\documentclass[a4paper,10pt
]{amsart}
\usepackage[english]{babel}

\usepackage{amsmath,amssymb,amsthm}
\usepackage[pagewise]{lineno}
\usepackage[bookmarksopen,final]{hyperref}
\usepackage[alphabetic, abbrev, msc-links,nobysame]{amsrefs}
\usepackage{enumerate}
\usepackage[all]{xy}
\newdir{ >}{{}*!/-5pt/@{>}} 
\SelectTips{cm}{} 

\numberwithin{equation}{section} 
\usepackage[utf8]{inputenc}
\newtheorem{theorem}[subsection]{Theorem}
\newtheorem{corollary}[subsection]{Corollary}
\newtheorem{lemma}[subsection]{Lemma}
\newtheorem{proposition}[subsection]{Proposition}
\theoremstyle{definition}
\newtheorem{definition}[subsection]{Definition}

\newtheorem{remark}[subsection]{Remark}
\newtheorem{example}[subsection]{Example}

\newtheorem{construction}[subsection]{Construction}

\newcommand{\bC}{\mathbb{C}}
\newcommand{\bE}{\mathbb{E}}
\newcommand{\bF}{\mathbb{F}}

\newcommand{\bR}{\mathbb{R}}
\newcommand{\bS}{\mathbb{S}}
\newcommand{\bT}{\mathbb{T}}
\newcommand{\bZ}{\mathbb{Z}}

\newcommand{\cC}{\mathcal{C}}
\newcommand{\cD}{\mathcal{D}}
\newcommand{\cE}{\mathcal{E}}
\newcommand{\cF}{\mathcal{F}}
\newcommand{\cG}{\mathcal{G}}
\newcommand{\cI}{\mathcal{I}}
\newcommand{\cJ}{\mathcal{J}}

\newcommand{\cO}{\mathcal{O}}

\newcommand{\cS}{\mathcal{S}}
\newcommand{\cV}{\mathcal{V}}
\newcommand{\cW}{\mathcal{W}}

\DeclareMathOperator*{\eq}{eq}

\DeclareMathOperator{\colim}{colim}

\DeclareMathOperator{\const}{const}

\DeclareMathOperator{\hocolim}{hocolim}

\renewcommand{\:}{\colon}

\newcommand{\<}{\langle}
\newcommand{\Alg}{\mathrm{Alg}}

\newcommand{\Ar}{\mathrm{Ar}}
\newcommand{\Assoc}{\mathrm{Assoc}}
\newcommand{\CAlg}{\mathrm{CAlg}}
\newcommand{\Cyc}{\mathrm{Cyc}}

\newcommand{\Env}{\mathrm{Env}}

\newcommand{\Fin}{\cF\mathrm{in}}
\newcommand{\Fun}{\mathrm{Fun}}
\newcommand{\GL}{\mathrm{GL}}
\newcommand{\HH}{\mathrm{HH}}

\newcommand{\Irep}{{I\text{-}\mathrm{rep}}}

\newcommand{\Map}{\mathrm{Map}}
\newcommand{\Mod}{\mathrm{Mod}}
\newcommand{\Pic}{\mathrm{Pic}}

\newcommand{\Sp}{\mathrm{Sp}}
\newcommand{\TC}{\mathrm{TC}}

\newcommand{\THH}{\mathrm{THH}}
\newcommand{\TP}{\mathrm{TP}}
\newcommand{\Th}{\mathrm{Th}}
\newcommand{\TwCyc}{\mathrm{TwCyc}}
\newcommand{\act}{\mathrm{act}}

\newcommand{\can}{\mathrm{can}}
\newcommand{\cof}{\mathrm{cof}}
\newcommand{\coll}{\mathrm{coll}}

\newcommand{\cy}{\mathrm{cy}}
\newcommand{\dlog}{d\kern-1pt\log}
\newcommand{\down}{{\textstyle (-\frac12)}}
\newcommand{\ev}{\mathrm{ev}}

\newcommand{\gp}{\mathrm{gp}}
\newcommand{\id}{\mathrm{id}}
\newcommand{\incl}{\mathrm{incl}}

\newcommand{\kup}{ku_{(p)}}

\newcommand{\longto}{\longrightarrow}
\newcommand{\modmod}[2]{#1/\!/#2}

\newcommand{\op}{\mathrm{op}}

\newcommand{\pr}{\mathrm{pr}}
\newcommand{\qqandqq}{\qquad\text{and}\qquad}
\newcommand{\rep}{\mathrm{rep}}
\newcommand{\res}{\mathrm{res}}
\newcommand{\sd}{\mathrm{sd}}

\newcommand{\trf}{\mathrm{trf}}

\newcommand{\up}{{\textstyle (+\frac12)}}
\newcommand{\zerorep}{{\pi_0\text{-}\mathrm{rep}}}
\renewcommand{\>}{\rangle}

\newcommand{\arxivlink}[1]{\href{https://arxiv.org/abs/#1}{\texttt{arXiv:#1}}}
\newcommand{\doilink}[1]{\href{http://dx.doi.org/#1}{doi:#1}}

\title[Localization sequences for log \texorpdfstring{$\TC$}{TC}]{Localization sequences for logarithmic topological cyclic homology}

\author{John Rognes} \address{Department of Mathematics, University of Oslo, 
  Norway}
\email{rognes@math.uio.no}

\author{Steffen Sagave}
\address{IMAPP, Radboud University Nijmegen, The Netherlands}
\email{s.sagave@math.ru.nl}

\author{Christian Schlichtkrull}
\address{Department of Mathematics, University of Bergen, 
  Norway}
\email{christian.schlichtkrull@math.uib.no}

\begin{document}
\date{June 25, 2026}

\begin{abstract}
  We introduce the notion of an $\bE_k$-ring with prelogarithmic structure, define logarithmic topological Hochschild homology and logarithmic topological cyclic homology in this context, and establish localization sequences for these theories. Our approach is based on Thom $R$-algebras. It extends and strengthens our earlier work on the subject in several regards. Our examples include the fraction field of topological $K$-theory, the existence of which was suggested by calculations by Ausoni and the first author. To illustrate the computational accessibility of log $\THH$ and log $\TC$, we determine these for non-negative even periodic sphere spectra, with their canonical prelogarithmic structures. 
\end{abstract}

\maketitle
\tableofcontents
\section{Introduction}
Trace maps from algebraic $K$-theory to topological Hochschild homology ($\THH$) and topological cyclic homology ($\TC$) provide powerful tools for the subtle computations of algebraic $K$-theory of discrete rings or structured ring spectra. Partly triggered by work of Nikolaus--Scholze~\cite{NS18} and Bhatt--Morrow--Scholze~\cite{BMS19}, these methods have seen a boost in recent years, witnessed for example by important results of Hahn--Wilson~\cite{HW22} and Burklund--Schlank--Yuan~\cite{BSY22} on chromatic redshift, of Burklund--Hahn--Levy--Schlank on Ravenel's telescope conjecture~\cite{BHLS23}, and of Antieau--Krause--Nikolaus on the algebraic $K$-theory of $\mathbb Z/p^n$~\cite{AKN24}.

However, there are interesting examples of discrete rings or structured ring spectra for which trace maps do not give good approximations to algebraic $K$-theory. This can be explained by the failure of $\THH$ and $\TC$ to admit the same localization sequences as algebraic $K$-theory. Motivated by a connection between relative $\THH$ and the de\ Rham complex with log poles~\cite{K89}, first appearing in work of Hesselholt--Madsen~\cite{HM03}, one strategy to overcome this is to extend the definition of $\THH$ and $\TC$ to appropriately defined prelog or log ring spectra, and to realize the terms missing from the localization sequences by the $\THH$ and $\TC$ of these objects. In this paper, we propose new notions of prelog and log $\bE_k$-rings, construct $\THH$ and $\TC$ for such objects, and establish localization sequences motivated by those in algebraic $K$-theory. Using simplicially enriched Waldhausen categories, Blumberg--Mandell~\cite{BM20} have in certain cases given a different construction of relative $\THH$ and $\TC$ terms fitting in such localization sequences. While their approach connects well to algebraic $K$-theory, it appears to have no direct relation to the de\ Rham complex with log poles. 

Our present approach extends and strengthens our earlier work on the subject~\cites{Rog09,RSS15,RSS18} in several regards: While our previous localization sequences and calculations were limited to $\bE_\infty$-rings with prelog structures generated by single homotopy classes (such as $ku$ with the Bott element $u \in \pi_2(ku)$), the scope of our new methods encompasses even $\bE_k$-rings with $k\geq 2$ and prelog structures generated by several different homotopy classes, such as the truncated Brown--Peterson spectrum $\mathrm{BP}\< n\>$ at a prime~$p$ with $p, v_1, \dots, v_n \in \pi_*\mathrm{BP}\< n\>$. More substantially, while our earlier work only addressed log $\THH$ as a cyclic spectrum, we now also construct cyclotomic structures using the graded refinement of the Nikolaus--Scholze approach that appears in work of Antieau--Mathew--Morrow--Nikolaus~\cite{AMMN22}. This allows us to define log $\TC$ for prelog $\bE_k$-rings, and we obtain the desired localization sequences both for log $\THH$ and log $\TC$. 
\subsection{\texorpdfstring{$R$}{R}-based prelog \texorpdfstring{$\bE_k$}{Ek}-rings and Thom spectra} A (discrete) prelog ring $(A,M,\bar\alpha)$ consists of a commutative ring $A$, a commutative monoid $M$, and a ring homomorphism $\bar\alpha \: \bZ[M] \to A$ from the integral monoid ring. It determines a localization $A[M^{-1}] = A \otimes_{A[M]}A[M^\gp]$, and one can interpret $(A,M,\bar\alpha)$ as an intermediate object between $A$ and $A[M^{-1}]$ that does not exist within commutative rings. (Making this factorization explicit requires discussing the log condition, which we skip in this introduction for simplicity.)

When generalizing prelog rings to the homotopy theoretic setup of structured ring spectra, we want to allow $A$ to be an $\bE_k$-ring for some $2 \leq k \leq \infty$. However, different choices suggest themselves when generalizing the remaining data. Most obviously, we can consider a discrete commutative monoid or an $\bE_k$-space $M$ together with a structure map $\bS[M] \to A$ from the spherical monoid ring. This works well when~$A$ is the Eilenberg--Mac\,Lane spectrum of a commutative ring~\cite{Rog09}, and also generalizes to the animated situation~\cites{SSV16,BLPO23-HKR}. However, it turns out to be too restrictive for capturing the prelog structure generated by a positive dimensional homotopy class in an $\bE_k$-ring. 

In~\cites{RSS15,RSS18}, we overcame this problem by working with a notion of graded $\bE_\infty$-spaces. However, the specific point set level model used there only allowed us to model $\bE_\infty$-spaces graded over the $\bE_\infty$-space $QS^0 = \Omega^\infty\Sigma^\infty S^0$, which in turn restricted the prelog structures we were able to construct. Our new solution to this problem is to use an $\infty$-categorical setup and work Picard graded relative to a base $\bE_{\infty}$-ring~$R$. That is, we consider $\bE_k$-spaces $M$ graded over the $\bE_{\infty}$-space $\Pic_R$, defined as the maximal subgroupoid spanned by the invertible $R$-modules in the symmetric monoidal $\infty$-category $\Mod_R$. There is an $\bE_k$ Thom $R$-algebra functor $\Th_R\: \Alg_{\bE_k}(\cS_{/\Pic_R}) \to \Alg_{\bE_k}(\Mod_R)$ that generalizes the spherical monoid ring  $\bS[-] \: \Alg_{\bE_k}(\cS) \to \Alg_{\bE_k}(\Sp)$.

We define an \emph{$R$-based prelog $\bE_k$-ring} to be a triple $(A,\xi,\bar\alpha)$ consisting of an $\bE_k$-ring $A$, an $\bE_k$-map $\xi \: M \to \Pic_R$, and an $\bE_k$-map ${\bar\alpha \: \Th_R(\xi) \to A}$ from the underlying $\bE_k$-ring of $\Th_R(\xi)$. (We also consider a more structured variant where~$A$ is an $\bE_k$ $R$-algebra and $\bar\alpha$ is an $\bE_k$ $R$-algebra map.) When $R = \bS$, we drop the term \emph{$\bS$-based} and simply speak of \emph{prelog $\bE_k$-rings}. While the examples in the present paper are mostly $\bS$-based, we develop a general $R$-based theory because this will allow us to consider $\bE_k$-maps $\xi \: M \to \Pic_R$ that do not factor through $\Pic_\bS \to \Pic_R$. 

One useful choice for the structure map $\xi$ of a prelog $\bE_2$-ring is the $\bE_2$-map $\xi_{2d}$ given by the composite
\[ \bZ_{\geq 0} \xrightarrow{1 \mapsto d} \bZ \xleftarrow{\simeq} \Omega^2(BU(1)) \xrightarrow{\Omega^2(\incl)} \Omega^2(BU) \xleftarrow{\simeq} \bZ\times BU  \xrightarrow{J_\bC}\Pic_\bS\,,\]
where $J_\bC$ denotes the complex $J$-homomorphism. It sends $1 \in \bZ_{\geq 0}$ to $2d \in \pi_0(\Pic_\bS)$. Using a minimal $\bE_2$-cell structure on $\Th_\bS(\xi_{2d})$ established by Galatius--Kupers--Randal-Williams~\cite{GKRW}, an obstruction theory argument allows us to build an $\bE_2$-map $\bar\alpha(a) \: \Th_\bS(\xi_{2d}) \to A$ hitting a prescribed homotopy class $a \in \pi_{2d}(A)$ whenever $\pi_*(A)$ is sufficiently sparse, e.g., even.  This defines a prelog $\bE_2$-ring
\begin{equation}\label{eq:S-based-E_2-prelog-intro}
(A,\xi_{2d},\bar\alpha(a)) = (A,\xi_{2d}\: \bZ_{\geq 0} \to \Pic_{\bS},\bar\alpha(a)\: \Th_\bS(\xi_{2d}) \to A). 
\end{equation}
For the connective complex topological $K$-theory spectrum $ku$ or the Adams summand $\ell$, with their respective Bott and Adams elements, we thus obtain the existence of prelog $\bE_2$-rings 
\[(ku,\< u \>) := (ku, \xi_{2}, \bar\alpha(u)) \qqandqq  (\ell,\< v_1\>) := (\ell, \xi_{2p-2}, \bar\alpha(v_1))\]
with prelog structures generated by these classes. 
The obstruction theory also applies to Thom spectra of  $\bE_2$-maps $(\bZ_{\geq 0})^r \to \Pic_\bS$ and allows us for example to define prelog $\bE_2$-rings $(\ell, \< p, v_1\>)$ and $(\mathrm{BP}\< n\>, \< p,v_1,\dots, v_n\>)$.

\subsection{Log \texorpdfstring{$\THH$}{THH}} The topological Hochschild homology $\THH(A)$ of an $\bE_k$-ring $A$ can be defined as the realization of the cyclic bar construction $B^\cy\: [q] \mapsto A^{\otimes 1+q}$. It comes equipped with a circle action because it is the realization of a cyclic object, and has an $\bE_{k-1}$-algebra structure because $B^\cy$ is symmetric monoidal. Moreover, there are $\bT$-equivariant cyclotomic structure maps $\varphi_p\: \THH(A) \to \THH(A)^{tC_p}$ to the $C_p$-Tate construction on $\THH(A)$ that equip $\THH(A)$ with the structure of a cyclotomic spectrum in the sense of~\cite{NS18}, and these data can be used to define the topological cyclic homology $\TC(A)$ of $A$. 

The key idea for extending $\THH$ to $R$-based prelog $\bE_k$-rings comes from a variant of the cyclic bar construction known as the \emph{replete bar construction}. Motivated by the algebraic structure of $\pi_*\THH(A|K)$ from~\cite{HM03}, the replete bar construction was introduced in~\cite{Rog09}. For a discrete monoid $M$, it is defined as the pullback of $M \to M^\gp \leftarrow B^\cy(M^\gp)$, i.e., of the group completion map of $M$ and the augmentation of the cyclic bar construction of $M^\gp$ induced by the iterated multiplication maps $(M^\gp)^{\times 1+q} \to M^\gp$. However, generalizing the last step to the present $\bE_k$-setup is delicate because there is no natural augmentation $B^\cy(M) \to M$ for a general $\bE_k$-space $M$. We resolve this issue by viewing the group completion~$M^\gp$ of an $\bE_k$-space $M$ with $k\geq 2$ as an object augmented over its commutative monoid of path components $\pi_0(M^\gp) \cong (\pi_0(M))^\gp$, so that there is an augmentation $B^\cy(M^\gp) \to B^\cy(\pi_0(M^\gp)) \to \pi_0(M^\gp)$ and we can base change along $\pi_0(M) \to \pi_0(M^\gp)$ to (re-)define $B^\rep(M)$.  

To implement this basic idea in the above setup of $R$-based prelog $\bE_k$-rings, we view the extension of an $\bE_k$-map $\xi\: M \to \Pic_R$ to the group completion $\xi^\gp \: M^\gp \to \Pic_R$ as an object $\xi^\gp_*$ graded over $\pi_0(M^\gp)$, or over a more general commutative monoid. It gives rise to a $\pi_0(M^\gp)$-graded $\bE_k$ Thom $R$-algebra $\Th_R(\xi^\gp_*)$ and we form the $\THH$ of (the underlying $\bE_k$-ring of)  $\Th_R(\xi^\gp_*)$  in $\pi_0(M^\gp)$-graded spectra. Restricting it along the group completion map $\pi_0(M) \to \pi_0(M^\gp)$ gives an $\bE_{k-1}$-ring $\THH(\Th_R(\xi), \xi)$. It comes with a canonical $\bE_{k-1}$ \emph{repletion map} $\rho \: \THH(\Th_R(\xi)) \to \THH(\Th_R(\xi), \xi)$ and we define the log $\THH$ of an $R$-based prelog $\bE_k$-ring $(A, \xi, \bar\alpha)$ to be the relative tensor product
\[
\THH(A, \xi, \bar\alpha) := \THH(A) \otimes_{\THH(\Th_R(\xi))}
\THH(\Th_R(\xi), \xi) \,. 
\]

While we show that $\THH(\Th_R(\xi), \xi)$ can also be obtained as the Thom spectrum over a suitable replete bar construction, the above approach has the advantage that we can use~\cite{AMMN22} to obtain the cyclotomic structure on $\THH(A, \xi, \bar\alpha)$ that we need to define $\TC(A, \xi, \bar\alpha)$.

When $k = 2$ and  $R = \bS$ or $R = \bS_{(p)}$, our version of log $\THH$ for $(A, \xi_{2d}, \bar\alpha(a))$, including the construction of its cyclotomic structure, agrees with the variant of log $\THH$ studied by Ausoni--Bay{\i}nd{\i}r--Moulinos~\cite{ABM23}. Moreover, we show that in examples including $(ku,\< u \>)$ and $ (\ell,\< v_1\>)$, it is equivalent to the version studied in~\cites{RSS15, RSS18}. Hence the calculations and log $\THH$-{\'e}taleness results from~\cite{RSS18} are also valid in the present model. 

It will be interesting to study the even filtration on $\THH(A, \xi, \bar\alpha)$ and the resulting log prismatic and syntomic cohomology using the motivic filtration~\cites{HRW, P23} (compare also~\cite{BLPO23-log-prismatic}). We intend to pursue this in future work.
  
\subsection{Repletion--residue sequences} In the situation of~\eqref{eq:S-based-E_2-prelog-intro}, there is an $\bE_2$-collapse map $\Th_\bS(\xi_{2d}) \to \bS$ that is a left inverse to the unit map $\bS \to \Th_\bS(\xi_{2d})$. Using it, we can define the residue $\bE_1$ $A$-algebra $\modmod{A}{\bar\alpha(a)} := A \otimes_{\Th_\bS(\xi_{2d})} \bS$. The next theorem is a simplified version of our main structural result about log $\THH$ and log~$\TC$: 
\begin{theorem}\label{thm:repletion-residue-intro} For a prelog $\bE_2$-ring as in~\eqref{eq:S-based-E_2-prelog-intro}, there are repletion--residue cofiber sequences
\[
\THH(A) \xrightarrow{\rho} \THH(A, \xi_{2d}, \bar\alpha(a))
\xrightarrow{\res} \Sigma \THH(\modmod{A}{\bar\alpha(a)})
\]
of cyclotomic $\THH(A)$-modules and
\[        
  \TC(A)_p
  \xrightarrow{\rho} \TC(A, \xi_{2d}, \bar\alpha(a))_p
  \xrightarrow{\res} \Sigma \TC(\modmod{A}{\bar\alpha(a)})_p \\
\]
of $\TC(A)_p$-modules.      
\end{theorem}
In the case of $(\ell,\< v_1\>)$, this specializes to a cofiber sequence
\[ \THH(\ell) \to \THH(\ell,\< v_1\>) \to \Sigma 
  \THH(H\bZ_{(p)})\]
of cyclotomic $\THH(\ell)$-modules and a similar sequence in log $\TC$. This matches the analogous localization sequence in algebraic $K$-theory, with the periodic theory $L$ instead of $(\ell,\< v_1\>)$, which was established by Blumberg--Mandell~\cite{BM08}. A comparison between our $\THH(\ell,\< v_1\>)$ and Blumberg--Mandell's relative term $\THH(\ell|L)$ from~\cite{BM20}, and therefore also a trace map from $K(L)$ to $\THH(\ell,\< v_1\>)$ or $\TC(\ell,\< v_1\>)$, appears in recent work by Lundemo~\cite{Lun25}. 

We give an example of the computational accessibility of our construction: For $d > 0$, the Thom spectrum $\bS[x] = \Th_\bS(\xi_{2d})$ is an $\bE_2$-model for the non-negative $2d$-periodic sphere spectrum. It is part of a canonical prelog $\bE_2$-ring $(\bS[x], \xi_{2d}, \id)$, for which we can explicitly identify the terms in the log $\THH$ and log $\TC$ repletion--residue sequence (see Proposition~\ref{prop:THHSP} and Theorem~\ref{thm:TCSP}). In the $\TC$-case, the sequence is $p$-equivalent to \[
\textstyle \TC(\bS)_p \vee
	\bigvee_{i>0} \Sigma ((\bS^{2d})^{\otimes i})_{hC_i} \longto \TC(\bS)[S^1]_p \vee
	\bigvee_{i>0} \Sigma ((\bS^{2d})^{\otimes i})_{hC_i} \longto \Sigma \TC(\bS)_p\,. \]

\subsection{The fraction field of topological \texorpdfstring{$K$}{K}-theory} Theorem~\ref{thm:repletion-residue-intro} generalizes to prelog $\bE_k$-rings $(A, \xi\: \<x_1,\dots,x_r\> \to \Pic_\bS,\bar\alpha)$ where $\<x_1,\dots,x_r\>$ is the free commutative monoid on $r$ generators. In this case, we obtain an $r$-dimensional cube of cofiber sequences with initial term $\THH(A)$ and final term $\Sigma^r\THH( \modmod{A}{\bar\alpha})$. The general case is spelled out in Theorem~\ref{thm:identification-of-hty-cofiber-cube}.

In the special case $(\ell, \< p, v_1\>)$ of the Adams summand with the $\bE_2$-prelog structure generated by $p \in \pi_0(\ell)$ and $v_1\in \pi_{2p-2}(\ell)$, we get a diagram of horizontal and vertical cofiber sequences of cyclotomic $\THH(\ell)$-modules 
\begin{equation}\label{eq:fraction-field-3x3-diagram-intro}
\xymatrix@-.8pc{
\THH(\ell) \ar[r]^-{\rho} \ar[d]_-{\rho}
	& \THH(\ell, \< p \>) \ar[r]^-{\res} \ar[d]_-{\rho}
	& \Sigma \THH(\ell/p) \ar[d]_-{\Sigma \rho} \\
\THH(\ell,\< v_1\>) \ar[r]^-{\rho} \ar[d]_-{\res}
	& \THH(\ell, \< p, v_1\>) \ar[r]^-{\res} \ar[d]_-{\res}
	& \Sigma \THH(\ell/p, \< v_1\>)
		\ar[d]_-{\Sigma \res} \\
\Sigma \THH(\bZ_{(p)}) \ar[r]^-{\Sigma \rho}
	& \Sigma \THH(\bZ_{(p)},  \< p \>)
		\ar[r]^-{\Sigma \res}
	& \Sigma^2 \THH(\bZ/p)\rlap{\,,}
}
\end{equation}
where the lower right-hand square anti-commutes and the three other squares commute. For this diagram, we extend our definition of log $\THH$ to cover the case of the $\bE_1$-ring $\ell/p$ by setting $\THH(\ell/p, \< v_1\>) = \THH(\ell/p)\otimes_{
  \THH(\ell)} \THH(\ell, \< v_1\>)$. Such a diagram of cofiber sequences was considered by Ausoni and the first author in the unpublished version~\cite{AR09} of~\cite{AR12}. There, the log $\THH$ terms constructed here were defined in terms of iterated homotopy cofibers. Based on the calculations made there, the authors hypothesized the existence of a \emph{fraction field of topological $K$-theory} whose $\THH$, $\TC$ and algebraic $K$-theory should participate in localization sequences of the type we construct. In view of this, the present results confirm that the prelog $\bE_2$-ring $(\ell, \< p, v_1\>)$ should be viewed as a model for this fraction field of (the Adams summand of) topological $K$-theory. We intend to demonstrate the connection between the calculations in~\cites{AR09,AR12} and our present constructions in more detail in future work. 
\subsection{Organization}
Section~\ref{sec:Thom-spectra} is about $R$-algebra Thom spectra. Section~\ref{sec:group-compl-repletion} discusses $\bE_k$-group completions and relative variants thereof. Section~\ref{sec:pic-graded-log-structures} contains the basic definitions and constructions of $R$-based prelog $\bE_k$-rings. Section~\ref{sec:logTHH} features the construction of log $\THH$ and in Section~\ref{sec:log-TC} we introduce log $\TC$. Section~\ref{sec:THH-of-Thom-rings} is about (log) $\THH$ of Thom spectra. In Section~\ref{sec:logifications}, we introduce the log condition and the logification construction. Sections~\ref{sec:repletion-residue-one-generator} and~\ref{sec:multiple-generator-case} contain the construction of the repletion--residue sequences, including a proof of Theorem~\ref{thm:repletion-residue-intro} and the construction of the diagram~\eqref{eq:fraction-field-3x3-diagram-intro}. The two appendices discuss the relations between the present approach and our previous results~in~\cites{RSS15,RSS18}. 
\subsection*{Notation and conventions} Except from in the appendices, we work $\infty$-cat\-e\-gor\-i\-cal\-ly and use~\cites{L:HTT,L21,L:HA} as primary references. We write $\cS$ for the $\infty$-category of spaces and $\Sp$ for the $\infty$-category of spectra, often viewed as the symmetric monoidal $\infty$-categories $(\cS,\times,*)$ and $(\Sp,\otimes,\bS)$. 

\subsection*{Acknowledgments} The authors would like to thank Christian Ausoni, Tobias Barthel, \"Ozg\"ur Bay{\i}nd{\i}r, S{\o}ren Galatius, Alexander Kupers, Tommy Lundemo, and Oscar Randal-Williams for helpful discussions related to this project. Moreover, the authors would like to thank the referee for useful comments. They would also like to thank the Isaac Newton Institute for Mathematical Sciences, Cambridge, for support and hospitality during the programme ``Equivariant homotopy theory in context'', where work on this paper was undertaken. This work was partially supported by EPSRC grant EP/Z000580/1, a grant from the Simons foundation, and the Dutch Research Council (NWO) ENW-XL grant ``Symmetry on the interface of topology and higher algebra'' \href{https://doi.org/10.61686/DGZBU69810}{doi:10.61686/DGZBU69810}.

\section{Thom spectra}\label{sec:Thom-spectra}
Let $R$ be an $\bE_{\infty}$-ring. The goal of this section is to review the symmetric monoidal Thom $R$-module functor $\Th_R$ and to study its right adjoint.

\subsection{A construction principle for Thom spectrum functors}
We follow the approach to $\Th_R$ in~\cite{HL17}*{\S 3} (see also~\cite{ABGHR14},~\cite{ABG18},~\cite{ACB19}) and work in a slightly more general setup that will become useful in Appendix~\ref{app:Thom-comparison}.

Let $P$ be an $\bE_{\infty}$-space. The slice $\infty$-category $\cS_{/ P}$ inherits a symmetric monoidal structure~\cite{L:HA}*{Rem.~2.2.2.5} and there is a canonical symmetric monoidal functor $y_P \: P \to \cS_{/ P}$ sending an object $p$ to the corresponding map $\Delta^0 \to P$ (compare~\cite{HL17}*{Not.~3.1.2}). The functor $\infty$-category $\Fun(P^{\op},\cS)$ has a symmetric monoidal convolution product, and under the symmetric monoidal equivalence
\begin{equation}\label{eq:str-unstr-for-P}
\Fun(P^{\op},\cS) \simeq \cS_{/ P}
\end{equation}
(see~\cite{HL17}*{Not.~3.1.2} or~\cite{Ram22}*{Cor.~5.9}), $y_P$ corresponds to the Yoneda embedding of $P$ into its cocompletion $\Fun(P^{\op},\cS)$. 

Let $(\cD,\otimes)$ be a presentably symmetric monoidal $\infty$-category, that is, $\cD$ admits all small colimits and $\otimes$ preserves small colimits separately in each variable. The next proposition is a minor generalization of~\cite{HL17}*{Prop.~3.1.3}.
\begin{proposition}\label{prop:extenstion-to-cocompletion}
  Composition with $y_P \: P \to \cS_{/ P}$ induces an equivalence of $\infty$-categories $\mathrm{LFun}^{\otimes}(\cS_{/ P},\cD) \to \Fun^{\otimes}(P,\cD)$ from the $\infty$-category of colimit-preserving symmetric monoidal functors $\cS_{/ P} \to \cD$ to the $\infty$-category of symmetric monoidal functors $P \to \cD$. \qed
\end{proposition}
\begin{proof}
The argument for~\cite{HL17}*{Prop.~3.1.3} applies also in this situation (see~\cite{L:HA}*{Prop.~4.8.1.10 and Cor.~4.8.1.12} for more details). 
\end{proof}
\begin{corollary}[\cite{HL17}*{Cor.~3.1.4}] \label{cor:construction-Fiota} Let $\iota \: P \to \cD$ be a symmetric monoidal functor. Then there is an essentially unique colimit-preserving symmetric monoidal functor $F_{\iota} \: \cS_{/ P} \to \cD$ such that $F_{\iota} \circ y_P \simeq \iota$.\qed
\end{corollary}
\begin{remark}
In the motivating example where $\iota$ is the inclusion $\Pic_R \to \Mod_R$, this $F_{\iota}$ is the Thom $R$-module functor $\Th_R$, see Definition~\ref{def:ThR} below. 
\end{remark}

Under the equivalence~\eqref{eq:str-unstr-for-P}, the functor $F_{\iota}$ can by construction be identified with the left Kan extension of $\iota \: P \to \cS$ along the Yoneda-embedding $y_P$. Since $y_P$ is fully faithful, we can use~\cite{L:HTT}*{Def.~4.3.2.2} to describe the value of this left Kan extension on a functor $Z \: P^{\op} \to \cS$: It is the colimit over the full subcategory of $\Fun(P^{\op},\cS)_{/ Z}$ spanned by the objects $y_P(p) \to Z$ with $p$ an object of $P$. Using again~\eqref{eq:str-unstr-for-P}, $Z$ corresponds to a map $\xi\: X \to P$ and $y_P$ to the map sending an object of $P$ to the map $\Delta^0 \to P$ it determines. Under this identification, the $\infty$-category over which we take the colimit corresponds to $X$, and the functor we are taking the colimit over is the composite of the projection to $P$ with $\iota$. 
\begin{corollary}
For an object $\xi \: X \to P$ of $\cS_{/ P}$, the object $F_{\iota}(\xi)$ is naturally equivalent to $\colim_{X}\iota\xi$.\qed 
\end{corollary}

The next construction is motivated by~\cite{ACB19}*{Def.~3.3}. Given an object $A$ of~$\cD$, we define $P_{\downarrow A}$ by the following pullback of $\infty$-categories:
\begin{equation}\label{eq:def-PdownarrawA}\xymatrix@-1pc{
    P_{\downarrow A} \ar[rr] \ar[d] && \cD_{/A} \ar[d] \\
    P \ar[rr]^{\iota} && \cD\rlap{\,.}
    }
  \end{equation}
  Here $\cD_{/A}$ is the slice of $\cD$ over $A$. Since $\cD_{/A} \to \cD$ is a right fibration~\cite{L:HTT}*{dual of Prop.~2.1.2.1} and $P$ is assumed to be an $\infty$-groupoid, \cite{L:HTT}*{Lem.~2.1.3.3} implies that $P_{\downarrow A}$ is an $\infty$-groupoid. In particular, we can view $P_{\downarrow A} \to P$ as an object of $\mathcal S_{/ P}$. This defines a functor \[G_{\iota}\: \cD \to \cS_{/P}, \quad A \mapsto (P_{\downarrow A} \to P).\]

  For an object $A$ in $\cD$, we notice that the functor $F_{\iota}$ sends $G_{\iota}(A) = (P_{\downarrow A} \to P)$ to the colimit of the composite $P_{\downarrow A} \to P \to \cD$. Since this map factors through the canonical map $\cD_{/A} \to \cD$, and the latter map preserves colimits~\cite{L:HTT}*{Prop.~1.2.13.8}, we can form the colimit in $\cD_{/A}$ and then project to $\cD$. This defines a map from the colimit of $P_{\downarrow A} \to  \cD$ to $A$, which we denote by $\varepsilon_A \: F_{\iota}(G_{\iota}(A)) \to A$. It is natural by the naturality of taking slices, pullbacks, and colimits.

\begin{proposition}\label{prop:additive-adjunction}
  The functor $F_{\iota}$ is left adjoint to $G_{\iota}$ with counit $\varepsilon \:  F_{\iota}\circ G_{\iota} \to \id_{\mathcal{D}}$. 
\end{proposition}
\begin{proof}
  We apply~\cite{L21}*{Prop.~5.1.13}, which is the dual of~\cite{L:kerodon}*{Tag 02FX}. For this we need to check that the composite 
  \begin{equation}\label{eq:adjunction-on-mapping-spaces} \Map_{\cS_{/ P}}(\xi, G_{\iota}(A)) \xrightarrow{F_{\iota}} \Map_{\cD}(F_{\iota}(\xi), F_{\iota}(G_{\iota}(A))) \xrightarrow{(\varepsilon_A)_*} \Map_{\cD}(F_{\iota}(\xi), A) 
  \end{equation}
  is an equivalence for all $\xi\: X \to P$ in $\cS_{/ P}$ and $A$ in $\cD$. Writing the $\infty$-groupoid $X$ as the colimit over $X$ of a terminal diagram reduces the claim to the case where $X$ is a point.

  To verify this case, recall that $\mathcal D_{/A}$ is cocomplete and that $\mathrm{pr}\: \mathcal D_{/A} \to \cD$ preserves colimits. Writing $C'$ for the colimit of the map $P_{\downarrow A} \to \cD_{/A}$ from~\eqref{eq:def-PdownarrawA}, it follows that $C = \mathrm{pr}(C') = F_\iota(G_\iota(A))$ is a colimit of the composite $P_{\downarrow A} \to P \to \cD$. The choice of a colimit cone for $C'$ provides a canonical map $P_{\downarrow A} \to (\cD_{/A})_{/ C'}$ that sends the objects of $P_{\downarrow A}$ to the values of the natural transformation from $P_{\downarrow A} \to \cD_{/A}$ to the constant $P_{\downarrow A}$-diagram with value~$C'$. We now consider the composite
  \begin{equation}\label{eq:composite-slice-categories} P_{\downarrow A} \to (\cD_{/A})_{/C'} \to \cD_{/C} \to \cD_{/A}
  \end{equation}
  induced by the map just discussed, the projection $\mathrm{pr}\: \cD_{/A} \to \cD$, and $\varepsilon_A \: C \to A$. Informally, the first map sends an object $(\xi(*),\iota\xi(*) \to A)$ of $P_{\downarrow A}$ to
  \[\xymatrix@-1pc{
      \iota\xi(*)\ar[rr] \ar[d] && C \ar[d]\\ 
    A \ar[rr]^{=} && A \rlap{\,,} 
      }\]
    the second map sends this object to $\iota\xi(*) \to C$, and the last map sends the latter to $\iota\xi(*) \to A$. More formally, the natural transformation associated with the colimit cone for $C'$ provides an equivalence between the composite in~\eqref{eq:composite-slice-categories} and the base change of $\iota$ in~\eqref{eq:def-PdownarrawA}. Next we consider the following diagram
  \begin{equation}\label{eq:mapping-eq-on-slice-cats}\xymatrix@-1pc{ \Delta^0 \ar[rr]^{\xi} \ar@{=}[d] && P \ar[d]^{\iota} && P_{\downarrow A} \ar[ll]\ar[d] \\
      \Delta^0 \ar[rr]^{\iota \xi}\ar@{=}[d]&& \mathcal D \ar@{=}[d]&& \mathcal D_{/C} \ar[ll] \ar[d]\\
      \Delta^0 \ar[rr]^{\iota \xi}&& \mathcal D && \mathcal D_{/A} \ar[ll]
      }
    \end{equation}
in which the upper right-hand vertical map is the composite of the first two maps in~\eqref{eq:composite-slice-categories}. The natural equivalences $\Map_{\cS}(\Delta^0,P) \simeq P$ and   $\Map_{\cS}(\Delta^0,P_{\downarrow A}) \simeq P_{\downarrow A}$ and \cite{L21}*{Prop.~3.3.18} imply that the horizontal homotopy pullback of the first row is a model for the mapping space $\Map_{\cS_{/P}}(\xi, P_{\downarrow A} \to P)$. By~\cite{L21}*{Cor.~2.5.34}, the horizontal homotopy pullbacks of the second and third rows are models for the mapping spaces $\Map_{\cD}(\iota\xi(*),C)$ and $\Map_{\cD}(\iota\xi(*),A)$, respectively. By construction, when passing to pullbacks the vertical maps in~\eqref{eq:mapping-eq-on-slice-cats} induce the maps in the composite~\eqref{eq:adjunction-on-mapping-spaces}. The composite of the two right-hand vertical squares in~\eqref{eq:mapping-eq-on-slice-cats} is equivalent to the homotopy pullback~\eqref{eq:def-PdownarrawA} defining $P_{\downarrow A}$. Hence the induced map of horizontal homotopy fibers is an equivalence, showing the claim. 
\end{proof}
\begin{corollary}\label{cor:Giota-lax-sym-mon}
The functor $G_{\iota}\:\cD \to \cS_{/ P}$ is lax symmetric monoidal, that is, it is a map of $\infty$-operads between the underlying $\infty$-operads of the symmetric monoidal categories $\cD$ and $\cS_{/ P}$. 
\end{corollary}
\begin{proof}
  This is a consequence of~\cite{L:HA}*{Cor.~7.3.2.7}.
\end{proof}

\subsection{ Thom \texorpdfstring{$R$}{R}-modules} 
Recall that $R$ denotes an $\bE_\infty$-ring. Let $\Mod_R$ be the symmetric monoidal $\infty$-category of $R$-modules and $\Pic_R$ the maximal subgroupoid of $\Mod_R$ spanned by the invertible $R$-modules, which inherits the structure of a grouplike $\bE_\infty$-space. 
\begin{definition}\label{def:ThR}
The Thom $R$-module functor is the colimit-preserving symmetric monoidal left adjoint in the adjunction 
\[ \Th_R \: \cS_{/ \Pic_R} \rightleftarrows \Mod_R \: \Omega^R \]
that arises from applying Corollary~\ref{cor:construction-Fiota} and Proposition~\ref{prop:additive-adjunction} to the inclusion $\iota \: \Pic_R \to \Mod_R$. In particular, $\Th_R(\xi\: M \to \Pic_R) \simeq \colim (\iota\xi\: M \to \Mod_R)$. We write $\Omega^R$ for the right adjoint, which is lax symmetric monoidal by Corollary~\ref{cor:Giota-lax-sym-mon}. Informally, $\Omega^R(A)$ is the space over $\Pic_R$ of $R$-module maps $P \to A$, with $P$ ranging through $\Pic_R$. 
\end{definition}
Let $1 \le k \le \infty$ count a number of coherently commuting multiplications.  As a consequence of~\cite{L:HA}*{Rem.~7.3.2.13}, we get an induced adjunction 
\begin{equation}
  \label{eq:R-mod-E-k-adjunction} \Th_R \: \Alg_{\bE_k}(\cS_{/ \Pic_{R}}) \rightleftarrows \Alg_{\bE_k}(\Mod_R) \: \Omega^R
\end{equation}
on the categories of $\bE_k$-algebras (or, more generally, for algebras over an $\infty$-operad). By abuse of notation, we often also write $\Omega^R(A)$ for the total space of $\Omega^R(A) = ((\Pic_R)_{\downarrow A} \to \Pic_R)$ and leave the $\bE_k$-structure map implicit.

We will often be interested in the composites
\begin{equation}\label{eq:def-uTh}
u\Th_R \: \cS_{/ \Pic_{R}}\to \Sp \qqandqq  u\Th_R \: \Alg_{\bE_k}(\cS_{/ \Pic_{R}}) \to \Alg_{\bE_k}(\Sp) 
\end{equation}
of $\Th_R$ with the lax symmetric monoidal functor $u = (\bS \to R)^* \: \Mod_R \to \Sp$. More explicitly, $u\Th_R(\xi\: M \to \Pic_R) \simeq \colim (u\iota\xi\: M \to \Sp)$. 

\begin{remark}
  The adjunction~\eqref{eq:R-mod-E-k-adjunction} induces an equivalence of mapping spaces \[\Map_{ \Alg_{\bE_k}(\Mod_R)}(\Th_R(\xi),A) \simeq \Map_{\Alg_{\bE_k}(\cS_{/ \Pic_{R}})}(\xi,\Omega^R(A)),\] which recovers the universal property of Thom spectra from~\cite{ACB19}*{Thm.~3.5}. 
\end{remark}

We will now explain why $\Omega^R$ generalizes the functor sending a ring to its underlying multiplicative monoid. For this we specialize to the case $R = \bS$ and $k=1$, consider an $\bE_1$-ring $A$, and define a $\bZ$-graded monoid $\pi_{0,*}(\Omega^\bS(A))$ by letting $\pi_{0,n}(\Omega^\bS(A))$ be the set of path components of the homotopy fiber of the $\Pic_\bS$-augmented $\bE_1$-space $\Omega^\bS(A)$ over $n \in \pi_0(\Pic_\bS)$. Then  $\pi_{0,*}(\Omega^\bS(A))$ inherits an associative and unital graded multiplication from $A$. Moreover, on each $\pi_{0,n}(\Omega^\bS(A))$ we get a $\{\pm 1\} = \pi_1(\Pic_\bS)$-action as in the long exact sequence of homotopy groups.

With this multiplication and sign action, $\pi_{0,*}(\Omega^\bS(A))$ is a \emph{graded signed monoid} in the sense of~\cite{SS12}*{Def.~4.15}, meaning that sign action and multiplication are suitably compatible. The $\bZ$-graded ring $\pi_*(R)$ also forms a graded signed monoid by forgetting all of the additive structure but the sign action. The following statement will be proved at the end of Appendix~\ref{app:Thom-comparison}:

\begin{proposition}\label{prop:Omega-S-multiplicative-monoid} Let $A$ be an $\bE_1$-ring. Then $\pi_{0,*}(\Omega^\bS(A))$ and $\pi_*(A)$ are isomorphic as graded signed monoids. \end{proposition}
This proposition motivates our approach to (pre-)log structures in the next sections, although only the comparison in Lemma~\ref{lem:map-to-direct-image} relies on it. 

\begin{remark}
While there are model categorical versions of the right adjoint to the Thom spectrum functor when working with $\bZ \times BO$ instead of $\Pic_\bS$ or $\Pic_R$ (see Appendix~\ref{app:Thom-comparison} or~\cite{Sch24}), we are not aware of a reference where this adjunction is established for general $R$. In line with this, we are not aware of an $\infty$-categorical proof of Proposition~\ref{prop:Omega-S-multiplicative-monoid}, but have not pursued such an $\infty$-categorical proof ourselves. 
\end{remark}

\subsection{Weight-graded Thom spectra}\label{subsec:weight-graded-Thom}
We will often consider functors with values in slice categories and begin with the following general observation:
\begin{lemma}\label{lem:functors-to-slice-categories} Let $K$ be a simplicial set, $\cC$ an $\infty$-category, and $X$ an object of~$\cC$. Then a lift of a functor $F \: K \to \cC$ to $\cC_{/X}$ is the same data as a natural transformation $F \to \const_{K}\!X = (K\to\Delta^0)^*X$ from $F$ to the constant $K$-diagram with value~$X$.
\end{lemma}
In fact, $\Fun(K,\cC_{/X}) \to \Fun(K,\cC)$ and $\Fun(K,\cC)_{/\const_K \!X} \to \Fun(K,\cC)$ are equivalent as right fibrations over $\Fun(K,\cC)$, but we will not use this stronger statement. 
\begin{proof}[Proof of Lemma~\ref{lem:functors-to-slice-categories}]
Using~\cite{L21}*{Prop.~2.5.27}, we may replace $\cC_{/X}$ by the fat slice $\cC^{/X}$. By~\cite{L21}*{Lem.~4.4.7 and Prop.~4.3.14}, the fiber of the right fibration $\Fun(K,\cC^{/X}) \to \Fun(K,\cC)$ over $F$ is $\mathrm{map}_{\Fun(K,\cC)}(F,\const_{K}X)$. 
\end{proof}
  
We now return to Thom spectra. In our applications, they will often come with an additional \emph{weight} grading. To explain this structure, let $(I, +, 0)$ be a discrete commutative monoid, such as $\{0\}$, $\bZ_{\ge0}$ or $\bZ$, viewed as a symmetric monoidal $\infty$-category. Moreover, let $\cC$ be a symmetric monoidal $\infty$-category with all small colimits. Then an $I$-graded object in $\cC$ is a functor $I\to \cC$, and we view the functor $\infty$-category $\Fun(I,\cC)$ as a symmetric monoidal $\infty$-category equipped with the Day convolution product \cite{L:HA}*{\S2.2.6}.

When $\cC = \cS$, we note that any map of spaces $\tau \: M \to I$ is automatically a coCartesian fibration, classified as in \cite{L:HTT}*{\S3.2} by the functor $M_* \: I \to \cS$ mapping $i \in I$ to $M_i = \tau^{-1}(i)$.   We refer to $M_i$ as the \emph{weight~$i$} graded piece of~$M$. The colimit of~$M_*$ recovers $M = \coprod_{i \in I} M_i$. We have the following characterization of multiplicative structure on the weight grading:
\begin{lemma}\label{lem:equivalent-formulations-gradings} Let $I$ be a discrete commutative monoid and let $\tau \: M \to I$ be a map of spaces. Then $\tau \: M \to I$ is an $\bE_k$-map of spaces if and only if the functor $M_*$ is lax $\bE_k$-monoidal, which holds if and only if $M_*$ is an $\bE_k$-algebra in~$\Fun(I, \cS)$. \qed
\end{lemma}
Suppose also given an $\bE_k$-map $\xi \: M \to \Pic_R$, and let $\xi_i = \xi|M_i \: M_i \to \Pic_R$ for each $i \in I$.  Then mapping $i \in I$ to $\xi_i$ defines a lax $\bE_k$-monoidal functor $\xi_*\: I \to \cS_{/\Pic_R}$, hence also an $\bE_k$-algebra in $\Fun(I, \cS_{/\Pic_R})$.  To emphasize the domains of the $\xi_i$, we often present $\xi_*$ as a map $M_* \to c^*\Pic_R$, where $c^*$ is the restriction along $c\: I \to \{0\}$ (and we implicitly apply~Lemma~\ref{lem:functors-to-slice-categories}). The colimit of $\xi_*$ in $\cS_{/\Pic_R}$ over~$I$ recovers~$\xi$ in $\cS_{/\Pic_R}$. The functor $\Th_R$ of Definition~\ref{def:ThR} induces a symmetric monoidal functor $\Fun(I, \cS_{/\Pic_R}) \to \Fun(I,\Mod_R)$ taking $\xi_*$ to the $\bE_k$-algebra $\Th_R(\xi_*)$ in $\Fun(I, \Mod_R)$, mapping $i \in I$ to $\Th_R(\xi_i)$ in~$\Mod_R$.  Its colimit over~$I$ recovers $\Th_R(\xi) = \bigvee_{i \in I} \Th_R(\xi_i)$.  Since the monoid pairing~$\cdot$ on~$M$ takes $M_i \times M_j$ to $M_{i+j}$ for all $i, j \in I$, the underlying $R$-algebra pairing on $\Th_R(\xi)$ takes $\Th_R(\xi_i) \otimes_R \Th_R(\xi_j)$ to $\Th_R(\xi_{i+j})$. Likewise, for the functor $u\Th_R$ of~\eqref{eq:def-uTh}, $u \Th_R(\xi_*)$ is an $\bE_k$-algebra in $\Fun(I, \Sp)$, with underlying $\bE_1$-ring pairing taking $u \Th_R(\xi_i) \otimes u \Th_R(\xi_j)$ to $u \Th_R(\xi_{i+j})$.

The following definition is our main source for weight gradings:
\begin{definition}\label{def:canonical-pi0grading}
Let $M$ be an $\bE_k$-space such that $\pi_0(M)$ is commutative. The $0$-th Postnikov truncation $\tau\: M \to \tau_{\leq 0}M \xrightarrow{\simeq} \pi_0(M)$ is an $\bE_k$-map and defines the \emph{canonical $\pi_0(M)$-grading} of $M$. 
\end{definition}
\begin{example}
For an $\bE_k$-ring $A$, Proposition~\ref{prop:Omega-S-multiplicative-monoid} implies that $\pi_0(\Omega^\bS(A)) \cong \pi_*(A)/\{\pm 1\}$. If  $k \geq 2$, then $\pi_*(A)$ is graded commutative, $\pi_*(A)/\{\pm 1\}$ is commutative, and  $\Omega^\bS(A)$ is canonically graded by $\pi_*(A)/\{\pm 1\}$. 
\end{example}

\subsection{Collapse maps}\label{subsec:collapse}
Let $f \: J \to I$ be a homomorphism of discrete commutative monoids. Then restriction and left Kan extension along $f$ induce an adjunction 
\begin{equation}\label{eq:restriction-extension-adjunction}
f_! \: \Fun(J, \Mod_R) \rightleftarrows \Fun(I, \Mod_R) \: f^*
\end{equation}
with a lax symmetric monoidal right adjoint $f^*$ and a symmetric monoidal left adjoint $f_!$. If we assume that $f\: J \to I$ is the inclusion of a submonoid, then the left and right Kan extensions
\[
f_! \simeq f_* \: \Fun(J, \Mod_R) \longto \Fun(I, \Mod_R)
\]
along~$f$ both map $Y_* \: J \to \Mod_R$ to $f_! Y_* \simeq f_* Y_*$
given by sending $i \in I$ to~$Y_i$ if $i \in J$ and to $0$ for $i \in I \setminus J$.  These naturally equivalent functors are both fully faithful and symmetric monoidal.

The following condition on a submonoid $J\subseteq I$ will be useful for controlling the multiplicative properties of collapsing the part of $\Th_R(\xi_*)$ complementary to $J$.  
\begin{definition}[\cite{Ogu18}*{\S1.4}]\label{def:face}
A \emph{face} of a discrete commutative monoid~$I$ is a submonoid~$J$
such that $i+j \in J$ for $i, j \in I$ only if $i \in J$ and $j \in J$.
\end{definition}
Equivalently, $J \subseteq I$ is a face if the complement $I \setminus J$ is  a prime ideal in~$I$ (again, compare~\cite{Ogu18}*{\S1.4}). 

If $J$ is a face of~$I$, then $f^*$ is symmetric monoidal because we have
\[ (f^*(X_*\otimes_R Y_*))_n = \bigvee_{\substack{i,j\in I\\ i+j=n}} X_i \otimes_R Y_j = \bigvee_{\substack{i,j\in J\\ i+j=n}} X_i \otimes_R Y_j = (f^*(X_*)\otimes_R f^*(Y_*))_n\] for $n \in J$, and $f^*$ preserves the unit. Moreover, the
localization functor $f_* f^*$ is compatible with the symmetric monoidal
structure on~$\Fun(I, \Mod_R)$, so by \cite{L:HA}*{Prop.~2.2.1.9}
it is a symmetric monoidal functor receiving a symmetric monoidal natural transformation from the identity. Applied to the $\bE_k$-algebra $\Th_R(\xi_*)$
we get a natural collapse $\bE_k$-map of $I$-graded $\bE_k$ $R$-algebras 
\[
\coll^{I\setminus J}_* \: \Th_R(\xi_*) \longto f_* f^* \Th_R(\xi_*)
\]
sending $i$ to $\Th_R(\xi_i)$ for $i \in J$ and to~$0$ for $i \in I \setminus J$.  Passing to total objects, i.e.,
colimits over~$I$, we obtain the $\bE_k$ $R$-algebra collapse map
\begin{equation}\label{eq:col-on-ThR}
  \coll^{I\setminus J} \: \Th_R(\xi) = \bigvee_{i \in I} \Th_R(\xi_i)
  \longto \bigvee_{i \in J} \Th_R(\xi_i) = \Th_R(\xi|M_J)
\end{equation}
collapsing the summands with $i \in I \setminus J$ to~$0$.  Here $M_J =
\tau^{-1}(J)$. We note that this construction is functorial for face inclusions, that is, for faces $K \subseteq J$ and $J \subseteq I$ we have $\coll^{I\setminus K} = \coll^{J\setminus K}\circ \coll^{I\setminus J}$.

These collapse maps will become relevant for the construction of the quotient objects appearing in our localization sequences (see Sections~\ref{subsec:trivial-locus-collapse}, \ref{sec:repletion-residue-one-generator} and \ref{sec:multiple-generator-case}). For these applications, it is useful to observe that $\coll^{I\setminus J}$ has structure not captured by the above construction. To see this, we consider the following more general construction also relevant in the next section.

Let $F\: \cC \to \cD$ be a symmetric monoidal left adjoint functor between presentably symmetric monoidal $\infty$-categories, with right adjoint $G$, and let $A$ be an $\bE_k$-algebra in $\cC$. Then $F(A)$ is an $\bE_k$-algebra in $\cD$, the categories $\mathrm{LMod}_A$ (of left $A$-modules) and $\mathrm{LMod}_{F(A)}$ inherit $\bE_{k-1}$-monoidal structures, and $F$ induces the left adjoint in an adjunction
\begin{equation}\label{eq:Ek-1-alg-adjunction-from-sym-mon-adj}
  \Alg_{\bE_{k-1}}(\mathrm{LMod}_A(\cC)) \rightleftarrows  \Alg_{\bE_{k-1}}(\mathrm{LMod}_{F(A)}(\cD))\,.
\end{equation}
Since $F$ is symmetric monoidal, its right adjoint $G$ is lax symmetric monoidal. So the right adjoint in~\eqref{eq:Ek-1-alg-adjunction-from-sym-mon-adj} can be identified with the composite of the functor \[\Alg_{\bE_{k-1}}(\mathrm{LMod}_{F(A)}(\cD)) \to \Alg_{\bE_{k-1}}(\mathrm{LMod}_{GF(A)}(\cC))\] induced by $G$ and the restriction along the adjunction unit $A \to G(F(A))$.

In the above situation where $J\subseteq I$ is a face, we can apply this to $F = f^*$ and $A = \Th_R(\xi_*)$. Evaluating the adjunction unit on the $\bE_{k-1}$ $\Th_R(\xi_*)$-algebra $\Th_R(\xi_*)$ and passing to total objects, we obtain the following result:
\begin{lemma}\label{lem:bEk-1-alg-structure-on-collapse}
  The underlying $\bE_{k-1}$-map of the collapse map \eqref{eq:col-on-ThR} extends to an $\bE_{k-1}$ $\Th_R(\xi)$-algebra structure on $\Th_R(\xi|M_J)$. \qed
\end{lemma}

\section{Group completion and repletion}\label{sec:group-compl-repletion}
As before, let $R$ be an $\bE_\infty$-ring, $k \geq 1$ a number of coherently commuting multiplications, and $I$ a discrete commutative monoid. 
\subsection{Group completion}
For any $\bE_k$-space $M$, the $\bE_k$-space $M^\gp = \Omega BM$ provides a group completion of $M$ that comes with a natural group completion $\bE_k$-map $\gamma \: M \to M^\gp$ (see~\cite{Leh24}*{Prop.~2.8} for a modern treatment).
\begin{lemma}\label{lem:detecting-group-completion-with-B}
A map $M \to M'$ in $\Alg_{\bE_k}(\cS)$ is equivalent to the group completion $M \to M^\gp$ if and only if $M'$ is grouplike and $BM \to BM'$ is an equivalence. 
\end{lemma}
\begin{proof}
This follows because~\cite{Leh24}*{Prop.~2.8} implies that $M' \to (M')^\gp$ is an equivalence if $M'$ is grouplike. 
\end{proof}
To lift $(-)^\gp$ to our relative and graded situations, it will be convenient to view it as a localization. For this, we write $\bE_k(*)$ for the free $\bE_k$-space on the one point space in (unbased) spaces and consider the set of morphisms
\[ W = \{ \gamma_{\bE_k(*)} \: \bE_k(*) \to \bE_k(*)^\gp, \id\:  \bE_k(*)^\gp \to \bE_k(*)^\gp\}\]
of $\Alg_{\bE_k}(\cS)$. By arguing with the path components, we see that $\Alg_{\bE_k}(\cS)^{\gp}$ is the collection of \emph{weakly $W$-local objects} in the sense of~\cite{L:kerodon}*{Tag~04LH}. Moreover, since the pushout of $\bE_k(*)^\gp \leftarrow \bE_k(*) \to \bE_k(*)^\gp$ in $\Alg_{\bE_k}(\cS)$ is $\bE_k(*)^\gp$, the set $W$ is closed under taking \emph{relative codiagonals} (\cite{L:kerodon}*{Tag~04K7}), i.e., the relative codiagonal of every morphism in $W$ belongs to $W$. Hence~\cite{L:kerodon}*{Tag~04LY} implies that $\Alg_{\bE_k}(\cS)^{\gp}$ also coincides with the collection of $W$-local objects in $\Alg_{\bE_k}(\cS)$.

\begin{lemma}\phantomsection\label{lem:grouplike-Ek-as-reflection}\begin{enumerate}[(i)] 
  \item  The subcategory  $\Alg_{\bE_k}(\cS)^{\gp}$ of grouplike objects is a reflective localization of $\Alg_{\bE_k}(\cS)$. 
  \item The $\Alg_{\bE_k}(\cS)^{\gp}$-reflection of an $\bE_k$-space $M$ can be obtained by applying the $\infty$-categorical small object argument~\cite{L:kerodon}*{Tag~04MA} to $W$ and $M$.
  \item The $\Alg_{\bE_k}(\cS)^{\gp}$-reflection of $M$ is equivalent to the group completion of $M$. 
  \end{enumerate}  
\end{lemma}
\begin{proof} Part (i) follows from ~\cite{L:kerodon}*{Tag~04MD} and part (ii) follows from ~\cite{L:kerodon}*{Tag~04MD}. The last part follows from (ii), Lemma~\ref{lem:detecting-group-completion-with-B}, and the observation that $B$ preserves colimits.
\end{proof}
Let $G$ be a grouplike $\bE_k$-space and let $W_{\!/G}$ be the set of morphisms
\[ \{ (M\to G) \to (M'\to G)\mid M\to M' \in W\} \]
in $\Alg_{\bE_k}(\cS)_{/G}$. 
\begin{lemma}\label{lem:grouplike-in-slice}
The set $W_{\!/G}$ is closed under taking relative codiagonals, and an object $M \to G$ in $\Alg_{\bE_k}(\cS)_{/G}$ is (weakly) $W_{\!/G}$-local if and only if $M$ is grouplike. 
\end{lemma}
\begin{proof}
  The first assertion follows from the definition of pushouts in $\Alg_{\bE_k}(\cS)_{/G}$ and the second from the definitions of $W_{\!/G}$ and of weakly local objects.
\end{proof}
Let $I^\gp$ be the algebraic group completion of $I$. In the rest of this subsection, we will use the shorthand notation 
\[\Alg^{(k)} = \Alg_{\bE_k}(\Fun(I^\gp,\cS_{/\Pic_R}))\]
for the category of $I^\gp$-graded $\bE_k$-algebras in $\cS_{/\Pic_R}$. We say that an object $\zeta_* \: N_* \to c^*\Pic_R$ in $\Alg^{(k)}$ is \emph{grouplike} if the underlying total $\bE_k$-space $N$ is. 

\begin{corollary}\label{cor:graded-group-completion-as-localization} The full subcategory $\Alg^{(k),\gp}\subseteq \Alg^{(k)}$ of grouplike objects is a reflective localization. It is spanned by the (weakly) $W'$-local objects for a set $W'$ of morphisms that is closed under taking relative codiagonals.
\end{corollary}
\begin{proof}
Identifying $\Alg^{(k)}$ with $(\Alg_{\bE_k}(\cS)_{/\Pic_R})_{/ I^\gp}$ as in Lem\-ma~\ref{lem:equivalent-formulations-gradings}, this follows from the previous lemma, taking $G$ to be $I^\gp\times \Pic_R$. 
\end{proof}
In analogy with Lemma~\ref{lem:grouplike-Ek-as-reflection}, this shows that the reflection provides a group completion that can be constructed using the small object argument with respect to~$W'$. 

Our next aim is to construct a relative group completion. For this we assume $k\geq 2$. As in adjunction~\eqref{eq:Ek-1-alg-adjunction-from-sym-mon-adj}, we consider $\bE_{k-1}$-algebras in left modules over an object $\zeta_* \: N_* \to c^*\Pic_R$ in $\Alg^{(k)}$ and denote this $\infty$-category by \[\Alg^{(k-1)}_{\zeta_*} = \Alg_{\bE_{k-1}}(\mathrm{LMod}_{\zeta_*}(\Fun(I^\gp,\cS_{/\Pic_R})))\,. \]  Writing $\iota_*$ for the initial object of $\Alg^{(k)}$, the $\infty$-category $\Alg^{(k-1)}_{\iota_*}$ is equivalent to $\Alg^{(k-1)}$ and the unit map $i\: \iota_* \to \zeta_*$ induces an adjunction
\begin{equation}\label{eq:free-forget-for-Ek-1-zeta} i_! \: \Alg^{(k-1)} \rightleftarrows \Alg^{(k-1)}_{\zeta_*}\: i^*\,\,.
\end{equation}
We say that an $\bE_{k-1}$ $\zeta_*$-algebra $\zeta'_*$ is \emph{grouplike} if the underlying $\bE_{k-1}$ algebra $i^*(\zeta'_*)$ is grouplike in $\Alg^{(k-1)}$.

\begin{corollary}\phantomsection\begin{enumerate}[(i)]
  \item The set $i_!(W')$ arising from the set $W'$ of Corollary~\ref{cor:graded-group-completion-as-localization} is closed under taking relative codiagonals. 
  \item The full subcategory $\Alg^{(k-1),\gp}_{\zeta_*} \subseteq \Alg^{(k-1)}_{\zeta_*}$ of grouplike objects is a reflective localization. It is spanned by the (weakly) $i_!(W')$-local objects and the group completion can be constructed using the small object argument for $i_!(W')$. 
  \end{enumerate}
\end{corollary}
\begin{proof}
Part (i) follows because $i_!$ preserves pushouts. Part (ii) is again analogous to Lemma~\ref{lem:grouplike-Ek-as-reflection} and uses the adjunction $(i_!,i^*)$.
\end{proof}
This corollary provides another way to group-complete an object $\zeta_*$ of $\Alg^{(k)}$: Viewing $\zeta_*$ as an object of $\Alg^{(k-1)}_{\zeta_*}$, we can apply the group completion functor from the last corollary to get an $\bE_{k-1}$ $\zeta_*$-algebra $\zeta_*^\gp$. This construction neither determines nor is determined by the group completion in $\Alg^{(k)}$ (compare~\cite{L:HA}*{Warning~7.1.3.9}). However, we have the following compatibility:
\begin{proposition}\label{prop:compatible-group-compl}
The unit map $\zeta_* \to \zeta_*^\gp$ of the  $\bE_{k-1}$ $\zeta_*$-algebra structure on $\zeta_*^\gp$ is a group completion of $\zeta_*$ in $\Alg^{(k-1)}$. In particular, it is equivalent to the $\bE_{k-1}$-map underlying the group completion $\zeta_* \to \zeta_*^\gp$ of $\zeta_*$ as an $\bE_k$-algebra. 
\end{proposition}

\begin{proof}
Lemma~\ref{lem:detecting-group-completion-with-B} implies that for the $\bE_k$-group completion $\zeta_* \to \zeta_*^\gp$, the underlying $\bE_{k-1}$-map is an $\bE_{k-1}$-group completion. 

For the remaining first claim, we need to show that that the unit map $\zeta_* \to \zeta_*^\gp$ is a group completion in $\Alg^{(k-1)}$. We may assume that $\zeta_*^\gp$ is obtained from $\zeta_*$ by the small object argument with respect to $i_!(W')$. Since $i^*$ preserves filtered (or, more generally, sifted) colimits, it is enough to show that that when $f$ in $\Alg^{(k-1)}_{\zeta_*}$ is a pushout of a morphism in $i_!(W')$, then $i^*(f)$ is the pushout of a morphism in $W'$. This holds because $i^*$ preserves pushouts (compare Remark~\ref{rem:equivalence-relative-tensor-products} below) and any map in $i^*i_!(W')$ is the pushout of a map in $W'$.
\end{proof}

\subsection{Repletion}\label{subsec:repletion} We now set up the relative group completion that will be a crucial ingredient for our constructions of log $\THH$.  
\begin{construction}\label{cons:xigp} Suppose given $\bE_k$-maps $\xi \: M \to \Pic_R$ and $\tau \: M \to I$. Writing $\gamma \: I \to I^\gp$ for the algebraic group completion of $I$, we can form $\xi_*^\gp = (\gamma_!\xi_*)^\gp$ in $\Alg_{\bE_k}(\Fun(I^\gp,\cS_{/\Pic_R}))$ and view it as an $\bE_{k-1}$ $\gamma_!\xi_*$-algebra or an $\bE_k$-algebra under $\gamma_!\xi_*$ (compare Proposition~\ref{prop:compatible-group-compl}). 
Encoding the gradings by maps to $I$ and $I^\gp$, these maps can be summarized by the commutative diagram
\begin{equation} \label{eq:xigp}
\xymatrix@-1pc{
I \ar[d]_-{\gamma}
	&& M \ar[ll]_-{\tau} \ar[rr]^-{\xi} \ar[d]^-{\gamma}
	&& \Pic_R  \\
I^\gp
	&& M^\gp \ar[ll]_-{\tau^\gp} \ar[rru]_-{\xi^\gp}
	\rlap{\,.}
}
\end{equation}
By the symmetric monoidality of $\Th_R$, the group completion map $\gamma \: M \to M^\gp$ induces a map $\Th_R(\gamma) \: \gamma_!\Th_R(\xi_*) \to \Th_R(\xi_*^\gp)$ of $\bE_k$ $R$-algebras. By Proposition~\ref{prop:compatible-group-compl}, its underlying $\bE_{k-1}$-map extends to an $\bE_{k-1}$ $\gamma_!\Th_R(\xi_*)$-algebra structure on $\Th_R(\xi^\gp_*)$. In particular, $\Th_R(\xi^\gp_*)$ is an $I^\gp$-graded $\bE_k$ $R$-algebra, i.e., an $\bE_k$-algebra in $\Fun(I^\gp, \Mod_R)$.
\end{construction}

\begin{definition}\label{def:repletion-on-Th}
In this situation, the \emph{repletion map} is the natural map 
\begin{equation}\label{eq:repletion-on-Th}
  \rho_* = \Th_R(\gamma_*) \: \Th_R(\xi_*)
  \longto \gamma^* \Th_R(\xi^\gp_*)
  = \Th_R(\xi^\gp_*) \circ \gamma
\end{equation}
of $I$-graded $\bE_k$ $R$-algebras, given in weight~$i \in I$ by the
map of Thom $R$-modules $\rho_i = \Th_R(\gamma_i) \: \Th_R(\xi_i) \to
\Th_R(\xi^\gp_{\gamma(i)})$ induced by the map $\gamma_i \: M_i \to
M^\gp_{\gamma(i)}$ over $\Pic_R$.  
\end{definition}

Arguing with the adjunction~\eqref{eq:Ek-1-alg-adjunction-from-sym-mon-adj} in the case of the symmetric monoidal left adjoint $\gamma_!$ implies the following statement:
\begin{corollary}\label{cor:repletion-as-alg-structure}
The underlying $\bE_{k-1}$-map of the repletion map~\eqref{eq:repletion-on-Th} extends to an $\bE_{k-1}$ $\Th_R(\xi_*)$-algebra structure on $\gamma^* \Th_R(\xi^\gp_*)$.\qed
\end{corollary}

\begin{definition}\label{def:zero-replete} Let $k\geq 2$. If $\tau \: M \to I$ is an $\bE_k$-map map to a discrete commutative monoid $I$, then $M$ is \emph{$I$-replete} if the canonical map $M \to I \times_{I^{\gp}} M^{\gp}$ is an equivalence. An $\bE_k$-space $M$ is \emph{$\pi_0$-replete} if it is $I$-replete for the canonical $\pi_0(M)$-grading, that is, if $M \to \pi_0(M)\times_{\pi_0(M^\gp)}M^\gp$ is an equivalence. 
\end{definition}
It is clear from the definitions that $\rho_*$ is an equivalence if $M$ is $I$-replete. 

\begin{remark}\label{rem:terms-replete-exact}
For discrete and commutative $M$, the condition that $M$ is $I$-replete
corresponds to $\tau \: M \to I$ being \emph{exact} in the
sense of~\cite{Ogu18}*{Def.~2.1.15}, and repletion corresponds to
\emph{exactification} in the terminology of~\cite{Ogu18}*{\S4.2}.  The
term `replete' was chosen in~\cite{Rog09} to avoid further overloading
the word `exact', and `repletion' seems to be easier to say than
`exactification'.
\end{remark}

\section{Picard-graded prelogarithmic structures}\label{sec:pic-graded-log-structures}
As before, let $R$ in $\CAlg(\Sp)$ be a fixed base $\bE_\infty$-ring, such as~$\bS$, $MU$ or~$H\bZ$, let $\Mod_R$ be the stable presentably symmetric monoidal $\infty$-category of $R$-modules, and let $u \: \Mod_R \to \Sp$ be the lax symmetric monoidal forgetful functor. Moreover, let $1 \le k \le \infty$ count a number of coherently commuting multiplications.

The following (Picard graded) notions generalize the prelog and log rings of logarithmic geometry (see e.g.~\cite{Ogu18}) and the versions of prelog and log ring spectra studied by the authors in earlier work (see e.g.~\cites{Rog09,RSS15}). Appendix~\ref{app:pointsetlevel-logTHH-comparision} contains a comparison with the latter objects. 

\begin{definition}
An \emph{$R$-based prelog $\bE_k$-ring} $(A, \xi, \bar\alpha)$
consists of an $\bE_k$-ring $A$, an $\bE_k$-map $\xi \: M \to
\Pic_R$, and an $\bE_k$-ring map $\bar\alpha \: u \Th_R(\xi) \to A$.

A \emph{prelog $\bE_k$ $R$-algebra} $(A, \xi, \bar\alpha)$ consists of
an $\bE_k$ $R$-algebra~$A$, an $\bE_k$-map $\xi \: M \to \Pic_R$,
and an $\bE_k$ $R$-algebra map $\bar\alpha \: \Th_R(\xi) \to A$.
\end{definition}
Prelog $\bE_k$ $R$-algebras are the objects of the $\infty$-category given by the pullback of
\[ \Alg_{\bE_k}(\cS_{/ \Pic_R}) \xrightarrow{\Th_R} \Alg_{\bE_k}(\Mod_R) \xleftarrow{\ev_0} \Alg_{\bE_k}(\Mod_R)^{\Delta^1}, \]
and analogously for $R$-based prelog $\bE_k$-rings. A map of prelog $\bE_k$ $R$-algebras $(f,f^\flat)\: (A,\xi,\bar\alpha) \to (B,\eta,\bar\beta)$ consists of maps $f\: A \to B$ in $\Alg_{\bE_k}(\Mod_R) $ and $f^\flat \: M \to N$ in $\Alg_{\bE_k}(\cS_{/ \Pic_R}) $ making the obvious square commutative, and analogously for  $R$-based prelog $\bE_k$-rings.

Each prelog $\bE_k$ $R$-algebra has an underlying $R$-based prelog
$\bE_k$-ring.  If $R = \bS$, these two notions coincide, and we speak of a
\emph{prelog $\bE_k$-ring}. When the context is clear, we refer to $(\xi,\bar\alpha)$ as a \emph{prelog structure} on $A$. 

\begin{example} The \emph{canonical} prelog $\bE_k$ $R$-algebra associated to
an $\bE_k$-map $\xi \: M \to \Pic_R$ is $(\Th_R(\xi), \xi, \id)$, where $\id \: \Th_R(\xi) \to \Th_R(\xi)$
is the identity map.
\end{example}

\subsection{\texorpdfstring{$\bE_2$ }{E2-}prelog structures} Our next goal is to build $\bE_2$ prelog structures on an $\bE_2$-ring determined by a given set of homotopy classes. 
\begin{construction}\label{cons:Sx1_xr}
  Let $r \geq 1$, let $d_1, \dots, d_r \in \bZ_{\geq 0}$, and let $\bZ_{\geq 0}^r = (\bZ_{\geq 0})^{\times r}$. We consider the composite
  \begin{equation}\label{eq:map-defining-Sx1xr}\xi_{2d_1, \dots, 2d_r}\:  \bZ_{\geq 0}^r\xrightarrow{(d_1,\dots, d_r)} \bZ^r \xrightarrow{(J_\bC\circ \Omega^2(\incl))^r} \Pic_\bS^r \to \Pic_\bS
\end{equation}
of the following maps: The first is the entry-wise multiplication by $d_i$ and the last is the iterated monoidal product of the $\bE_{\infty}$-space $\Pic_\bS$. The second is the $r$-fold product of the composite
\begin{equation}\label{eq:Z-to-PicS-map} \bZ \xleftarrow{\simeq} \Omega^2(BU(1)) \xrightarrow{\Omega^2(\incl)} \Omega^2(BU) \xleftarrow{\simeq} \bZ \times BU \xrightarrow{J_\bC} \Pic_\bS
\end{equation}
where $\Omega^2(\incl)$ is induced by the inclusion $\incl\: BU(1) \to BU$ and $J_\bC$ is the complex $J$-homomorphism, given as the composite of the realification $\bZ \times BU \to \bZ \times BO$ with the real $J$-homomorphism $J_\bR \: \bZ \times BO \to \Pic_\bS$ (see e.g.\ Appendix~\ref{app:Thom-comparison}). Since $\Omega^2(\incl)$ is $\bE_2$ and $J_\bC$ is $\bE_{\infty}$, the composite~\eqref{eq:Z-to-PicS-map} is $\bE_2$. The composite $\xi_{2d_1, \dots, 2d_r}$ in~\eqref{eq:map-defining-Sx1xr} is an $\bE_2$-map because the first and last map are $\bE_{\infty}$ and the middle one is $\bE_2$. The composite takes $(i_1,\dots, i_r) \in \bZ_{\geq 0}^r$ to a point in the path component of $2(d_1i_1+\dots+d_ri_r) \in \bZ \cong \pi_0(\Pic_\bS)$, whence the name $\xi_{2d_1, \dots, 2d_r}$. We write
\[ \bS[x_1,\dots, x_r] = \Th_{\bS}(\xi_{2d_1, \dots, 2d_r}) \]
for the associated Thom $\bE_2$-ring. The underlying $\bE_1$-ring of $\bS[x_1,\dots, x_r]$ is by construction the smash product $\bS[x_1]\otimes \dots \otimes \bS[x_r]$ of the $r$ free $\bE_1$-rings $\bS[x_1],\dots, \bS[x_r]$ on generators in degrees $|x_i| = 2d_i$. The identity of $\bZ_{\geq 0}^r$ provides a canonical weight grading in the sense of Subsection~\ref{subsec:weight-graded-Thom}.
\end{construction}
\begin{remark}\label{rem:E3-structure-on-Sx} This method cannot be improved to give an $\bE_3$-structure: Since $\widetilde{KU}^*(K(\bZ,3)) = 0$~\cite{AH68}*{Thm.~III}, we get $\widetilde{KO}^*(K(\bZ,3)) = 0$. So any map $K(\bZ,3) \to B^3(\mathbb Z\times BO)$ is nullhomotopic. Therefore, any $\bE_3$-map $\bZ \to {\bZ \times BO}$ is nullhomotopic, and we cannot realize $\bS[x]$ with $|x|>0$ as a Thom spectrum of an $\bE_3$-map $\bZ_{\geq 0} \to \Pic_\bS$ that factors through  $J_\bR \: {\bZ \times BO}\to \Pic_\bS$. There are also obstructions to other potential constructions of such $\bE_3$-structures: when $|x| = 2$ the $\bE_2$-structure on $\bS[x]$ (or its $2$-localization) cannot be lifted to an $\bE_3$-structure (see~\cite{DHLSW23}*{Rem.~3.11}). In Proposition~\ref{prop:obstructions-E3-structures-Sx} we show an analogous odd primary statement for $p\nmid |x|$.

However,  the further composites of $\bZ \to \Pic_\bS$ with the first or both of the $\bE_\infty$-maps
\[
\Pic_\bS \longto \Pic_{MU} \longto \Pic_{H\bZ}
\]
become $\bE_\infty$-maps, cf.~\cite{HL18}*{p.~85} or~\cite{D24}*{Lem.~2.1.3}. The composites
\[\bZ \to \Pic_{MU} \to \bZ\qqandqq \bZ \to \Pic_{H\bZ} \to \bZ\]
with the $0$-th Postnikov truncations $\Pic_{MU} \to \bZ$ and $\Pic_{H\bZ} \to \bZ$ send $1$ to~$2$, so there are factorizations $\Pic^{\ev}_{MU} \simeq \bZ \times BGL_1(MU)$ and $\Pic^{\ev}_{H\bZ} \simeq \bZ \times BGL_1(\bZ)$ as $\bE_\infty$-spaces, where $\Pic^{\ev}$ denotes the fiber over $2\bZ \subset \bZ$. Hence it is possible to consider $MU$-based prelog $\bE_\infty$-rings with $M = \bZ_{\geq 0}$ and $u\Th_{MU}(\xi_{2d}) \simeq \bigvee_{i\geq 0}\Sigma^{2di}MU$. When $d=1$, this defines the non-negative even periodic complex bordism spectrum $MUP_{\geq 0}$. 
\end{remark}

The inclusions of the generators $(0,\dots, 0,1,0, \dots 0)$ of $\bZ_{\geq 0}^r$ induce a map
\begin{equation}\label{eq:wedge-to-Sx_1x_r-map}\bS^{2d_1} \vee \dots \vee \bS^{2d_r} \to \bS[x_1,\dots, x_r]
\end{equation}
of spectra that is adjoint to a map of $\bE_2$-rings
\begin{equation}\label{eq:free-to-Sx_1x_r-map} \bE_2(\bS^{2d_1} \vee \dots \vee \bS^{2d_r}) \to \bS[x_1,\dots, x_r]\end{equation}
from the free $\bE_2$-ring on the source of~\eqref{eq:wedge-to-Sx_1x_r-map}.
\begin{proposition}
  The map~\eqref{eq:free-to-Sx_1x_r-map} can be obtained by attaching $\bE_2$-cells in positive even degrees, with equally many $d$-cells as the rank of the homogeneous degree~$d+2$ part of the cokernel of $\incl\: \bZ\{\sigma^2x_1,\dots,\sigma^2x_r\} \to \Gamma_{\bZ}(\sigma^2x_1,\dots,\sigma^2x_r)$, where $|\sigma^2x_i|=2d_i+2$. 
\end{proposition}
Here $\Gamma_{\bZ}(\sigma^2x_1,\dots,\sigma^2x_r) = \bZ\{\gamma_{k_1} \sigma^2 x_1 \cdots \gamma_{k_r} \sigma^2 x_r \mid k_1, \dots, k_r \ge 0\}$ denotes the divided power algebra over~$\bZ$ on the indicated generators. 
\begin{proof}
We will use~\cite{GKRW} and the claim that for $d\geq 0$, the map
\begin{equation}\label{eq:free-to-Sx_1x_r-map-in-E2} H_d^{\bE_2}(\bE_2(\bS^{2d_1} \vee \dots \vee \bS^{2d_r});\bZ) \to H_d^{\bE_2}(\bS[x_1,\dots, x_r];\bZ)
\end{equation}
in the $\bE_2$-homology of augmented $\bE_2$-rings induced by~\eqref{eq:free-to-Sx_1x_r-map} is isomorphic to the homogeneous degree $d+2$ part of $\incl\:\bZ\{\sigma^2x_1,\dots,\sigma^2x_r\} \to \Gamma_{\bZ}(\sigma^2x_1,\dots,\sigma^2x_r)$. Here $\bE_2$-homology is understood in the sense of~\cite{GKRW}*{Def.~10.7}.

We first explain how to deduce the proposition from this claim and~\cite{GKRW}*{Thm.~11.21}. In the notation of that theorem, we take $\mathsf{S}$ to be the category of spectra, $\mathsf{G}$ to be the discrete commutative monoid $\bZ_{\geq 0}^r$, and $\cO$ a non-unitary $\Sigma$-cofibrant $\bE_2$-operad (in the sense of~\cite{GKRW}*{Rem.~12.3}). The (canonical) rank functor $\bZ_{\geq 0}^r \to \bZ_{\geq 0}, (i_1,\dots,i_r) \mapsto i_1+\dots+i_r$ witnesses that $\mathsf{G}$ is Artinian in the sense of~\cite{GKRW}*{Def.~11.10}. In this situation, the map~\eqref{eq:free-to-Sx_1x_r-map} is a map of $\mathsf{G}$-graded objects, with the degree $(i_1,\dots,i_r)$-part of $\bS[x_1,\dots,x_r]$ a $2(d_1i_1+\dots+d_ri_r)$-fold suspension of $\bS$. Passing to the fibers of the augmentation to $\bS$ (viewed as a $\mathsf{G}$-graded object concentrated in degree $(0,\dots,0)$), we obtain a map $f \: \mathbf{R} \to \mathbf{S}$ of $\cO$-algebras in $\mathsf{S}^{\mathsf{G}}$. Our choice of grading ensures that $\mathbf{R}$ and $\mathbf{S}$ are reduced in the sense of~\cite{GKRW}*{Def.~11.11}. Finally, we choose the abstract connectivity in the sense of~\cite{GKRW}*{\S11.1} to be the map $c \: \mathsf{G} \to [-\infty, \infty]$ with $c(g) = 1$ if $g=(0,\dots,0)$ and $c(g) = 0$ otherwise. (An abstract connectivity leading to sharper estimates is not relevant for us since we will count cells in each dimension.)

Now our claim and the long exact sequence~\cite{GKRW}*{(10.1)} imply that 
\[ \bigoplus_{g \in \mathsf{G}} H_{g,d}^{\cO}(\mathbf{S},\mathbf{R};\bZ) \cong \left(\mathrm{coker}(\incl\:\bZ\{\sigma^2x_1,\dots,\sigma^2x_r\} \to \Gamma_{\bZ}(\sigma^2x_1,\dots,\sigma^2x_r))\right)_{d+2}.\]
Hence~\cite{GKRW}*{Thm.~11.21} applies, and it remains to determine the number of $d$-cells in the $\bE_2$-cell structure resulting from the theorem. We cannot apply the statement after that theorem directly because~\cite{GKRW}*{Ax.~11.19} requires us to work with $\bZ$-coefficients. However, the previous statement shows that each $H_{g,d}^{\cO}(\mathbf{S},\mathbf{R};\bZ)$ is finitely generated and free. With this, an inductive argument with the long exact sequences of the triples $\mathbf{Z}_{\epsilon -1} \to \mathbf{Z}_{\epsilon} \to \mathbf{S}$ appearing in the proof of the theorem shows that the sets of $(g,\epsilon)$-cells $I_{g,\epsilon}$ constructed there can be chosen so that the Hurewicz map and the inclusions induce isomorphisms \[\mathbb Z[I_{g,\varepsilon}] \xrightarrow{\cong} H_{g,\varepsilon}^{\cO}(\mathbf{S},\mathbf{Z}_{\epsilon -1} ;\bZ) \xleftarrow{\cong} H_{g,\varepsilon}^{\cO}(\mathbf{S},\mathbf{R};\bZ).\] 
This implies the desired statement about the number of $\bE_2$-cells. 

To verify our claim about~\eqref{eq:free-to-Sx_1x_r-map-in-E2}, we use the bar spectral sequence~\cite{GKRW}*{Thm.~14.1}. Firstly, we
observe that $H_*(\bS[x_1,\dots,x_r];\bZ) \cong \bZ[x_1,\dots,x_r]$. Next we apply the spectral sequence to compute $H_*(B^{\bE_1}(\bS[x_1,\dots,x_r]);\bZ)$. The $E^2$-page is the exterior algebra $\Lambda_\bZ(\sigma x_1,\dots, \sigma x_r)$ with $|\sigma x_i| = (1,2d_i)$. Then each  $\sigma x_i$ is an infinite cycle for filtration reasons, so there are no differentials from the $E^2$-page on, and no multiplicative extensions because
$(\sigma x_i)^2$ and $(\sigma x_i)(\sigma x_j) + (\sigma x_j)(\sigma x_i)$ cannot be detected in filtrations $0$ or~$1$. This implies that 
  \[ H_*(B^{\bE_1}(\bS[x_1,\dots,x_r]);\bZ) \cong \Lambda_\bZ(\sigma x_1,\dots, \sigma x_r) \]
  with $|\sigma x_i| = 2d_i+1$. We apply the spectral sequence a second time to obtain an additive isomorphism \[  H_*(B^{\bE_2}(\bS[x_1,\dots,x_r]);\bZ)  \cong \Gamma_{\bZ}(\sigma^2 x_1,\dots, \sigma^2 x_r).\]
  This time the spectral sequence collapses at the $E^2$-page because $\sigma^2 x_i$ has even total degree $2d_i + 2$. Using the isomorphism
  \[ H_*^{\bE_2}(\bS[x_1,\dots,x_r]; \bZ) \cong \tilde H_{*+2}(B^{\bE_2}(\bS[x_1,\dots,x_r]);\bZ) \] 
  from~\cite{GKRW}*{Thm.~13.7}, the claimed computation of $H_d^{\bE_2}(\bS[x_1,\dots, x_r];\bZ)$ follows. For the source of~\eqref{eq:free-to-Sx_1x_r-map-in-E2},~\cite{GKRW}*{Thm.~13.7 and 13.8} give that
  \[ H_d^{\bE_2}(\bE_2(\bS^{2d_1} \vee \dots \vee \bS^{2d_r});\bZ) \cong \tilde H_{d}(\bS^{2d_1} \vee \dots \vee \bS^{2d_r}; \bZ)\]
  for $d \geq 0$, with the generator corresponding to $\bS^{2d_i}$ matching the $\bE_2$-cell corresponding to $x_i$. Showing that the map in~\eqref{eq:free-to-Sx_1x_r-map-in-E2} is the claimed one reduces to the $r=1$ case by naturality, where it follows from~\cite{GKRW}*{Cor.~11.14} and the construction of~\eqref{eq:free-to-Sx_1x_r-map}.  
\end{proof}

\begin{remark}
  Analogous or related $\bE_2$-cell structures with $r = 1$ appear in~\cite{L:rotation}*{\S5.4} and~\cite{ABM23}*{Prop.~3.11}, while~\cite{HW22}*{Prop.~4.2.1} contains an $\bF_p$-version. 
\end{remark}

\begin{corollary}\label{cor:prelog-from-generators}
  Let $A$ be an $\bE_2$-ring with $\pi_*(A)$ concentrated in even degrees, and let $a_1,\dots,a_r$ be homotopy classes of degrees $|a_i| = 2d_i$ with $d_i\geq 0$. Then there exists an $\bE_2$-map $\bar \alpha = \bar \alpha(a_1,\dots,a_r) \: \bS[x_1,\dots,x_r] \to A$ such that $\pi_{2d_i}(\bar \alpha)$ sends the generator corresponding to $x_i$ to $a_i$.
\end{corollary}
In general, the map $\bar \alpha$ will not be unique. 

\begin{proof}
The chosen generators $a_1,\dots, a_r$ determine an $\bE_2$-map  $\bE_2(\bS^{2d_1} \vee \dots \vee \bS^{2d_r}) \to A$. The obstructions to extend it over the $\bE_2$-cells used to obtain $\bS[x_1,\dots,x_r]$ vanish by assumption, because we only need to attach even-dimensional cells. 
\end{proof}

\begin{remark}\label{rem:vanishing-for-one-gen-case} In the corollary, it is sufficient to ask that $\pi_{n-1}(A) \cong 0$ for all $n$ for which the relative cell complex~\eqref{eq:free-to-Sx_1x_r-map} has a positive number of relative $n$-cells. For example, when $r=1$, having $\pi_{j(2d_1+2)+2d_1-1}(A) \cong 0$ for all $j \geq 1$ is enough. When $r > 1$, the above weaker statement is sufficient for our examples and easier to work with. 
\end{remark}

\begin{definition}\label{def:ku-ell-BPn}
  When applied to the connective complex $K$-theory spectrum $ku$ and its $p$-local counterpart, with their Bott elements, this construction provides prelog $\bE_2$-rings
  \begin{linenomath*}\begin{align} & (ku,\< u \>) := (ku, \xi_{2}, \bar\alpha(u) \: \bS[x] \to
    ku)\qquad\text{and}\\
    &(\kup,\< u \>) := (\kup, \xi_{2}, \bar\alpha(u)\: \bS[x] \to \kup)  \end{align}\end{linenomath*}
  with $|x|=2$. In the case of the Adams summand,
  we similarly obtain
  \begin{linenomath*}\begin{equation}\label{eq:ell<v1>} (\ell,\< v_1\>) := (\ell, \xi_{2p-2}, \bar\alpha(v_1)\: \bS[x] \to \ell)
  \end{equation}\end{linenomath*}
with $|x|=2p-2$. These prelog $\bE_2$-rings admit logifications that are independent of the choices made in their constructions, see Lemma~\ref{lem:map-to-direct-image} and Example~\ref{ex:logification-of-E2-prelog-on-ell-ku} below. 

Taking $p$ instead of the Bott elements (with $|x|=0$) provides 
\begin{linenomath*}\begin{equation}\label{eq:ell<p>ku<p>} (ku,\< p \>), \quad (\kup,\< p \>), \qqandqq (\ell,\< p\>). 
  \end{equation}\end{linenomath*}
  The truncated Brown--Peterson spectrum $\mathrm{BP}\< n\>$ with $\pi_*\mathrm{BP}\< n\> \cong \bZ_{(p)}[v_1,\dots,v_n]$, $|v_i|=2p^{i}-2$, admits an $\bE_3$-structure by~\cite{HW22}*{Thm.~A}. Hence Corollary~\ref{cor:prelog-from-generators} applies to the elements $p, v_1, \dots, v_n$ and provides a prelog $\bE_2$-ring
  \begin{linenomath*}\begin{multline}\label{eq:BPnp-to-vn}(\mathrm{BP}\< n\>, \< p,v_1,\dots,v_n\>) :=\\  (\mathrm{BP}\< n\>, \xi_{0,2p-2,\dots,2p^n-2}, \bar\alpha(p, v_1,\dots,v_n)\: \bS[x_0,\dots,x_n] \to \mathrm{BP}\< n\>) 
  \end{multline}\end{linenomath*}
with $|x_i|=2p^i-2$, which depends both on the chosen $\bE_3$-structure and the chosen prelog structure map. This is especially interesting in the case $n=1$, where $\mathrm{BP}\< 1\> = \ell$ and this prelog $\bE_2$-ring specializes to
  \begin{equation}\label{eq:ell_p_v1} (\ell,\< p, v_1\>) =  (\ell, \xi_{0, 2p-2}, \bar\alpha(p, v_1)\: \bS[x_0,x_1] \to \ell).
  \end{equation}
Analogously, Corollary~\ref{cor:prelog-from-generators} provides $(ku,\< p, u \>)$ and $(\kup,\< p, u \>)$. 
\end{definition}

\subsection{Trivial locus and collapse}\label{subsec:trivial-locus-collapse} We introduce two constructions associating ring spectra to $R$-based prelog rings and prelog $R$-algebras.

\begin{definition} \label{def:triloc}
Let $(A, \xi \: M \to \Pic_R, \bar\alpha \: u \Th_R(\xi) \to A)$ be an
$R$-based prelog $\bE_k$-ring. By Proposition~\ref{prop:compatible-group-compl} and~\cite{L:HA}*{Cor.~7.1.3.4}, $u \Th_R(\xi^\gp)$ is an $\bE_{k-1}$ $u \Th_R(\xi)$-algebra. The \emph{trivial locus} of $(A,\xi,\bar\alpha)$ is the $\bE_{k-1}$ $A$-algebra 
\[
A[M^{-1}] := A \otimes_{u \Th_R(\xi)} u \Th_R(\xi^\gp)
\]
obtained as the base change of $u \Th_R(\xi^\gp)$ along $\bar\alpha$.
\end{definition}

\begin{remark}\label{rem:equivalence-relative-tensor-products}
If $(A, \xi, \bar\alpha)$ is a prelog $\bE_k$ $R$-algebra, then the relative tensor product in Definition~\ref{def:triloc} can be replaced by the relative tensor product
\[
A \otimes_{\Th_R(\xi)} \Th_R(\xi^\gp)
\]
 in $R$-modules  using the $\bE_k$ $R$-algebra maps $\bar\alpha$ and
 $\Th_R(\gamma) \: \Th_R(\xi) \to \Th_R(\xi^\gp)$. In this situation, \cite{L:HA}*{Thm.~4.5.2.1(2)} identifies the tensor product in~$\Mod_R$ with the relative tensor product over~$R$, which is introduced as a two-sided bar construction in~\cite{L:HA}*{\S4.4.2} and shown to satisfy associativity in~\cite{L:HA}*{\S4.4.3}.  It follows that the underlying spectrum of the construction in~$\Mod_R$ is equivalent to the construction in $\Sp$, so that the former provides an $R$-module enrichment of the latter.  In particular, we do not need to make~$R$ explicit in the notation~$A[M^{-1}]$.

 We also obtain an $\bE_{k-1}$ $R$-algebra structure on $A[M^{-1}]$ by restricting its $\bE_{k-1}$ $A$-algebra structure along $R \to A$. While this does not require $\bar\alpha$ to be an $\bE_k$ $R$-algebra map, we ask this property of $\bar\alpha$ in the definition of a prelog $\bE_k$ $R$-algebra to ensure the compatibility of this $\bE_{k-1}$ $R$-algebra structure with those on $\Th_R(\xi)$ and $\Th_R(\xi^\gp)$. A similar discussion applies to Definition~\ref{def:AMIJ} below.
\end{remark}

By construction, the prelog structure map $\bar\alpha \: u \Th_R(\xi)
\to A$ induces a map $u \Th_R(\xi^\gp) \to A[M^{-1}]$, where the map from
the Thom $R$-module of $\xi$ over~$M$ has been extended over~$M^\gp$, so that $M$ now acts invertibly.  This is intended to justify the notation~$A[M^{-1}]$.

\begin{example}\label{ex:non-neg-even-periodic-sphere-sp}
Taking $r = 1$ in Construction~\ref{cons:Sx1_xr}, we get the canonical prelog $\bE_2$-ring $(\bS[x],\<x\>) = (\bS[x],\xi_{2d}\: \< x\> \cong \bZ_{\geq 0} \to \Pic_\bS,\id)$. In this situation, the trivial locus $A[M^{-1}] = \Th_\bS(\xi_{2d}^\gp)$ is the $2d$-periodic sphere spectrum $\bigvee_{i\in\bZ} \bS^{2di}$. By Proposition~\ref{prop:compatible-group-compl}, we can view it as an $\bE_1$ $\bS[x]$-algebra or as an $\bE_2$-ring under $\bS[x]$, with equivalent underlying $\bE_1$-rings. 
\end{example}

\begin{example} For $(\ell,\< v_1\>)$ as in~\eqref{eq:ell<v1>}, $\ell[\< v_1 \>^{-1}] = \ell\otimes_{\Th_\bS(\xi_{2p-2})}\Th_\bS(\xi_{2p-2}^\gp)$ is the localization $L = E(1)$ of $\ell$ away from $v_1$, which has a preferred $\bE_{\infty}$ $\ell$-algebra structure. Analogously, we get ${KU}$ and ${KU}_{(p)}$ from $(ku,\< u \>)$ and $(\kup,\< u \>)$.\end{example}

\begin{definition} \label{def:AMIJ}
Let $(A, \xi \: M \to \Pic_R, \bar\alpha \: u \Th_R(\xi) \to A)$ be an
$R$-based prelog $\bE_k$-ring, let $\tau \: M \to I$ be an $\bE_k$-map
to a discrete commutative monoid, and let $J$ be a face of~$I$.  The \emph{collapse of $I\setminus J$} in $(A,\xi,\bar\alpha)$ is the $\bE_{k-1}$ $A$-algebra given by the base change \[
A/(M_{I \setminus J}) := A \otimes_{u \Th_R(\xi)} u \Th_R(\xi|M_J)
\]
of the $\bE_{k-1}$ $u \Th_R(\xi)$-algebra $u \Th_R(\xi|M_J)$ from Lemma~\ref{lem:bEk-1-alg-structure-on-collapse} along the $\bE_k$-map $\bar\alpha$. 

When $M = I = \bZ_{\geq 0}^r$, $\xi = \xi_{2d_1, \dots, 2d_r}$, and $J = \{0\}$, then we denote $A/(M_{I \setminus J})$ by $\modmod{A}{\bar\alpha}$, or by $\modmod{A}{\<x_1,\dots,x_r\>}$ if $\bar\alpha$ is understood from the context. 
\end{definition}
\begin{remark}
As explained in Remark~\ref{rem:equivalence-relative-tensor-products}, we can use the relative tensor product \[ A \otimes_{\Th_R(\xi)} \Th_R(\xi|M_J) \]
in $R$-modules instead when $(A, \xi, \bar\alpha \: \Th_R(\xi) \to A)$ is a prelog $\bE_k$ $R$-algebra.
\end{remark}

By construction, the collapse map $\coll^{I\setminus J} \: \Th_R(\xi) \to \Th_R(\xi|M_J)$
induces a map $A \to A/(M_{I \setminus J})$, where $M_{I \setminus J}
= \tau^{-1}(I \setminus J)$ now acts trivially on the target.  This is
intended to justify the notation $A/(M_{I \setminus J})$.  

\begin{example}\label{ex:collapse-canonical}
For $(\bS[x],\<x\>)$ as in Example~\ref{ex:non-neg-even-periodic-sphere-sp}, we get $\modmod{\bS[x]}{\<x\>} = \bS$ as $\bE_2$-algebras (since no base change is necessary).  
\end{example}

\begin{example}\label{ex:collapse-of-p-or-v1-in-ell} For $(\ell,\< v_1\>)$, we get $\modmod{\ell}{\<v_1\>} = H\bZ_{(p)}$ from applying $\ell\otimes_{\bS[x]}(-)$ to the cofiber sequence $\Sigma^{2p-2}\bS[x] \to \bS[x] \to \bS$ with $|x| = 2p-2$ as in~\eqref{eq:ell<v1>}.  Analogously, for $(ku,\< u \>)$ and $(\kup,\< u \>)$ we get $\modmod{ku}{\<u\>} = H\bZ$ and  $\modmod{\kup}{\<u\>} = H\bZ_{(p)}$. The same argument for the examples from~\eqref{eq:ell<p>ku<p>} gives the mod $p$ connective Adams summand $\modmod{\ell}{\<p\>} = \ell/p = k(1)$ and the mod $p$ connective complex $K$-theory spectrum $\modmod{ku}{\<p\>} = \modmod{\kup}{\<p\>} = ku/p$. 
\end{example}

\begin{example}\label{ex:collapse-of-p-and-v1-in-ell} For $(\ell,\< p, v_1\>)$, we get $\modmod{\ell}{\<p, v_1\>} \simeq \modmod{\ell}{\<p\>} \otimes_\ell \modmod{\ell}{\<v_1\>} \simeq H\bF_p$.
\end{example}

\begin{example}
  For $(\mathrm{BP}\< n\>, \< p,v_1,\dots,v_n\>)$ as in~\eqref{eq:BPnp-to-vn}, collapsing over the complement of the face $\< v_n\>$ gives the $n$-th connective Morava $K$-theory spectrum $k(n)$ at the prime $p$: This holds additively by arguing as in Example~\ref{ex:collapse-of-p-or-v1-in-ell}, and as $\bE_1$-rings because the $\bE_1$-structure on $k(n)$ is unique~\cite{A11}. Restricting the $\bE_2$-prelog structure on $\mathrm{BP}\< n\>$ to $\< v_n \>$ and composing with the $\bE_1$-map $\mathrm{BP}\< n\> \to k(n)$ from the collapse defines a prelog $\bE_1$-ring $(k(n),\< v_n \>)$. Its trivial locus is an $\bE_0$-ring that is equivalent as a spectrum to the $n$-th periodic Morava $K$-theory spectrum~$K(n)$. We can also get $K(n)$ as an $\bE_1$-ring from $(\mathrm{BP}\< n\>, \< p,v_1,\dots,v_n\>)$ by forming the relative $\otimes$-product $\mathrm{BP}\< n\> \otimes_{\bS[x_0,\dots,x_n]}\bS[x_n^{\pm 1}]\,$. This amounts to performing the collapse and trivial locus in one step. The identification of the $\bE_1$-structure again follows from its uniqueness~\cite{A11}.  
\end{example}
\section{Logarithmic topological Hochschild homology}\label{sec:logTHH}
Our next goal is to define the log $\THH$ of $R$-based prelog $\bE_k$-rings and prelog $\bE_k$ $R$-algebras with $k\geq 2$, including its cyclotomic structure. We begin by reviewing the cyclic bar construction and then turn to the construction of weight-graded $\THH$ from~\cite{AMMN22}*{App.~A}.

\subsection{The cyclic bar construction}
Recall that the forgetful functor from symmetric monoidal $\infty$-categories to $\infty$-operads has a left adjoint $\cC^\otimes \mapsto \Env(\cC)^\otimes$, given by the symmetric monoidal envelope from \cite{L:HA}*{\S2.2.4}. The underlying $\infty$-category $\Env(\cC)$ of $\Env(\cC)^\otimes$ agrees with the subcategory $\cC^\otimes_\act \subset \cC^\otimes$ spanned by all objects and the active morphisms between them \cite{L:HA}*{Def.~2.1.2.3, Rem.~2.2.4.3}.  The adjunction counit $\otimes \: \cC^\otimes_\act \to \cC$ is then informally given by $(X_1, \dots, X_n) \mapsto X_1 \otimes \dots \otimes X_n$ in the symmetric monoidal structure on~$\cC$, cf.~\cite{L:HA}*{Prop.~2.2.4.9} and \cite{NS18}*{Prop.~III.3.2}.

As in~\cite{L:HA}*{Def.~4.1.1.3}, we write $\Assoc$ for the associative $\infty$-operad. Let $\cC$ be a presentably symmetric monoidal $\infty$-category. Then an associative algebra $M \in \Alg(\cC)$ is given by an $\infty$-operad map $M^\otimes \: \Assoc^\otimes \to \cC^\otimes$ (compare~\cite{L:HA}*{Def.~4.1.1.6}) and restricts to a functor $M^\otimes_\act \: \Assoc^\otimes_\act \to \cC^\otimes_\act$ on the active morphisms. Let $\Delta$ denote the usual simplex category, let $\Lambda$ denote Connes' cyclic category, and let $\Lambda_\infty$ be the paracyclic category, as in \cite{NS18}*{p.~380}.  As explained in~\cite{NS18}*{Prop.~B.5} in the $\infty$-categorical context, the geometric realization $|X_\circ|$ of a cyclic object $X_\circ \: N(\Lambda^\op) \to \cC$ takes values in the category $\cC^{B\bT}$ of objects in $\cC$ with $\bT$-action. It is the colimit of the restriction of $X_\circ$ to $N(\Lambda_{\infty}^\op)$, and the underlying object in $\cC$ is equivalent to the colimit of the restriction $X_\bullet$ to $N(\Delta^\op)$.  Both of these colimits are sifted. Recall from \cite{NS18}*{Prop.~B.1} that there is a functor $V^\circ \: N(\Lambda^\op) \to \Assoc^\otimes_\act$ to the active part, whose restriction to $N(\Delta^\op)$ refines the functor $[q] \mapsto (\Delta^1/\partial\Delta^1)_q \cong \{0, 1, \dots, q\}$ represented by the simplicial circle $\Delta^1/\partial\Delta^1$.

\begin{definition}\label{def:Bcy} Let $\cC$ be a presentably symmetric monoidal $\infty$-category. The \emph{cyclic bar construction} is the composite \begin{multline}\label{eq:Bcy-composite}
    B^{\cy}\: \Alg(\cC) \xrightarrow{M \mapsto M^\otimes_\act} \Fun(\Assoc^\otimes_\act, \cC^\otimes_\act)\xrightarrow{(V^\circ)^* } \\
  \Fun(N(\Lambda^\op),\cC^\otimes_\act)  \xrightarrow{(\otimes)_*} \Fun(N(\Lambda^\op), \cC) \xrightarrow{|-|} \cC^{B\bT} 
\end{multline}
of the functors discussed above. 
\end{definition}
For $\cC = \Sp$, the first part of the next proposition is discussed in~\cite{NS18}*{\S{}IV.2}. 
\begin{proposition}\phantomsection\label{prop:Bcy-symm-mon-nat}\begin{enumerate}[(i)] 
  \item  The functor~\eqref{eq:Bcy-composite} is symmetric monoidal. 
  \item  If $F\: \cC \to \cD$ is a lax symmetric monoidal functor, then there is a natural comparison map $B^{\cy}\circ F \to F \circ B^{\cy}$ of lax symmetric monoidal functors $\Alg(\cC) \to \cD^{B\bT}$. If $G \: \cD \to \cE$ is another lax symmetric monoidal functor, then the composite
    \[ B^{\cy} \circ G \circ F \to G \circ B^{\cy} \circ F \to G \circ F \circ B^{\cy}\] 
    induced by the comparison maps for $F$ and $G$ is equivalent to the comparison map for $G\circ F$. 
  \item If $F\: \cC \to \cD$ is a symmetric monoidal functor commuting with sifted colimits, then the comparison map for $F$ is an equivalence.
  \end{enumerate}
\end{proposition}
\begin{proof}
For (i), we use that the symmetric monoidal structure on $\Alg(\cC)$ is defined in terms of the  pointwise symmetric  monoidal structure~\cite{L:HA}*{Rem.~2.1.3.4} on the functor category $\Fun(\Assoc^\otimes_\act, \cC^\otimes_\act)$. Therefore, $M \mapsto M^\otimes_\act$ is symmetric monoidal. The pointwise symmetric monoidal structure is compatible with precomposition by maps of simplicial sets, and with postcomposition by symmetric monoidal functors. Being a sifted colimit, geometric realization is symmetric monoidal. Hence it follows that $B^{\cy}$ is symmetric monoidal. 

For (ii), we observe that the first two functors in~\eqref{eq:Bcy-composite} commute with lax symmetric monoidal functors. To analyze the third functor, we note that the two composites in the square
\begin{equation}\label{eq:naturality-otimes-act} \xymatrix@-1pc{\cC^\otimes_\act \ar[rr]^{F^\otimes_\act} \ar[d]_{\otimes_\cC} && \cD^\otimes_\act \ar[d]^{\otimes_\cD} \\ \cC \ar[rr]^F && \cD }
\end{equation}
are not necessarily equivalent. However, using the characterization of the right adjoint of the adjunction $\incl\: \Fun_{\otimes}(\cC^\otimes_\act,\cD) \rightleftarrows  \Fun_{\mathrm{lax}}(\cC^\otimes_\act,\cD)\: R$ given in~\cite{NS18}*{Lem.~III.3.3}, we get an equivalence $\otimes_\cD \circ F^\otimes_\act \xrightarrow{\simeq} R(F\circ \otimes_\cC)$ of symmetric mo\-noid\-al functors because the precomposition of both functors with the lax symmetric mo\-noid\-al functor $\cC \to \cC^\otimes_\act$ is $F$.

The adjoint of this equivalence is the desired comparison map $\otimes_\cD \circ F^\otimes_\act \to F\circ \otimes_\cC$ of lax symmetric monoidal functors. Applying $\Fun(N(\Lambda^\op),-)$ this gives the comparison map for the third functor in~\eqref{eq:Bcy-composite}. For the last one, we use that the canonical comparison map $|F(-)| \to F(|-|)$ is a map of lax symmetric monoidal functors because $F$ is lax symmetric monoidal. The non-trivial part in the statement about $G\circ F$ is to prove that the composed map \[\otimes_\cE \circ G^\otimes_\act \circ F^\otimes_\act \to G \circ \otimes_\cD \circ F^\otimes_\act \to G \circ F\circ \otimes_\cC \]
is equivalent to the comparison map for $G\circ F$. By precomposing with $\cC \to \cC^\otimes_\act$, the adjoint of this composite can be checked to be equivalent to the equivalence $\otimes_\cE  \circ G^\otimes_\act \circ F^\otimes_\act \xrightarrow{\simeq} R(G\circ F\circ \otimes_\cC)$, from which the result follows. 

For (iii),~\cite{NS18}*{Prop.~III.3.2} implies that~\eqref{eq:naturality-otimes-act} commutes up to equivalence when $F$ is symmetric monoidal. Moreover, realization, being a sifted colimit, commutes with $F$ by assumption. 
\end{proof}
By Dunn additivity~\cite{L:HA}*{Thm.~5.1.2.2}, we can view $\bE_{k}$-algebras in $\cC$ as $\bE_{k-1}$-algebras in $\Alg(\cC)$. Hence the proposition has the following consequence: 
\begin{corollary}\phantomsection\label{cor:Bcy-on-E_k-alg}\begin{enumerate}[(i)]\item 
The cyclic bar construction induces a functor
\begin{equation}\label{eq:Bcy-on-E_k-alg} B^{\cy}\: \Alg_{\bE_k}(\cC) \to \Alg_{\bE_{k-1}}(\cC^{B\bT})\,.
\end{equation}
\item If $F\: \cC \to \cD$ is a lax symmetric monoidal functor and $A$ is any $\bE_k$-algebra in~$\cC$, then there is a natural map $B^{\cy}(F(A)) \to F(B^{\cy}(A))$ in $\Alg_{\bE_{k-1}}(\cD^{B\bT})$. This is compatible with compositions of lax symmetric monoidal functors.   
\item If $F\: \cC \to \cD$ is a symmetric monoidal functor commuting with sifted colimits, then the map in (ii) is an equivalence.\qed  
\end{enumerate}
\end{corollary}
We will also need a relative version of this corollary, which again follows from the symmetric monoidal structure of $B^\cy$:
\begin{corollary}\label{cor:Bcy-on-E_k-alg-relative} Let $R$ be an $\bE_\infty$-algebra in $\cC$, let $A$ be an $\bE_k$ $R$-algebra in $\cC$, and let $uA$ denote the underlying $\bE_k$-algebra of $A$ in $\cC$.
\begin{enumerate}[(i)]
\item The $\bE_{k-1}$-algebra with $\bT$-action $B^\cy(uA)$ resulting from Corollary~\ref{cor:Bcy-on-E_k-alg}(i) is naturally the underlying object of an $\bE_{k-1}$ $B^\cy(R)$-algebra with $\bT$-action. 
\item If $F\: \cC \to \cD$ is a lax symmetric monoidal functor, then there is a natural map $B^{\cy}(F(uA)) \to F(B^{\cy}(uA))$ of $\bE_{k-1}$-algebras with $\bT$-action over the $\bE_\infty$-algebra with $\bT$-action $B^{\cy}(F(R))$. This is compatible with compositions of lax symmetric monoidal functors.   
\item If $F\: \cC \to \cD$ is a symmetric monoidal functor commuting with sifted colimits, then the map in (ii) is an equivalence of $\bE_{k-1}$-algebras with $\bT$-action over the $\bE_\infty$-algebra with $\bT$-action $B^{\cy}(F(R)) \simeq F(B^{\cy}(R))$.
\end{enumerate}
\end{corollary}
\subsection{Weight-graded \texorpdfstring{$\THH$}{THH}}\label{subsec:weight-graded-THH}
Let $I$ be a discrete commutative monoid. An \emph{$I$-graded spectrum} is an object $X_*$ of $\Fun(I,\Sp)$, with total object $X = (c\: I \to \{0\})_! X_* = \bigvee_{i\in I}X_i$. We view $\Fun(I,\Sp)$ as a symmetric monoidal $\infty$-category equipped with the Day convolution product. An \emph{$I$-graded $\bE_k$-ring} is an $\bE_k$-algebra $A_*$ in $\Fun(I,\Sp)$, and has an $\bE_k$-ring  $A = c_!A_*$ as total object.

\begin{definition}\label{def:THHAstar} The \emph{$I$-graded topological
    Hochschild homology} \begin{equation}\label{eq:graded-THH-composite}
 \THH\: \Alg_{\bE_k}(\Fun(I, \Sp)) \to \Alg_{\bE_{k-1}}(\Fun(I, \Sp)^{B\bT})
\end{equation}
is the cyclic bar construction~\eqref{eq:Bcy-on-E_k-alg} for $\bE_k$-algebras in $\cC=\Fun(I,\Sp)$. It sends the underlying $I$-graded $\bE_1$-ring $A_* \in \Alg(\Fun(I, \Sp))$ of an $\bE_k$-algebra to the $I$-graded spectrum  with $\bT$-action $\THH(A_*) = \{i \mapsto \THH(A;i)\}$. 
\end{definition}

We now turn to cyclotomic structures. By~\cite{AMMN22}*{Def.~A.1, Lem.~A.7}, left Kan extension along $\ell_p \: I \to I, i \mapsto pi$ defines a symmetric monoidal functor $L_p = (\ell_p)_! \: \Fun(I, \Sp) \to \Fun(I, \Sp)$, and (for a class of symmetric monoidal $\infty$-categories including our discrete commutative monoid~$I$) there is a canonical \emph{$L_p$-twisted Tate diagonal}
\[
\Delta_p \: L_p(X_*) \longto (X_* \otimes \dots \otimes X_*)^{tC_p}
\]
on $\Fun(I, \Sp)$, with $p$ copies of $X_*$ in the target. Here $(-)^{tC_p}$ denotes the $C_p$-Tate construction (see e.g.~\cite{NS18}*{Ch.~I}). An \emph{$I$-graded twisted cyclotomic spectrum} is an $I$-graded spectrum $X_*$ with a $\bT$-action and a $\bT$-equivariant map $\varphi_p \: L_p(X_*) \to X_*^{tC_p}$ for every prime $p$. Here $X_*^{tC_p} = \{i \mapsto (X_i)^{tC_p}\}$ has the residual $\bT \cong \bT/C_p$-action for each $i \in I$. Analogously to~\cite{NS18}*{Def.~II.1.6}, we obtain a stable $\infty$-category $\TwCyc\Fun(I, \Sp)$ of $I$-graded twisted cyclotomic spectra. Applying~\cite{NS18}*{Con.~IV.2.1} provides a presentably symmetric monoidal structure on $\TwCyc\Fun(I, \Sp)$. 

\begin{proposition}[\cite{AMMN22}*{Cor.~A.9}]\label{prop:tw-cyc-structure-on-THH}
  The $I$-graded topological Hochschild homology canonically lifts to a symmetric monoidal functor \begin{equation*}\THH\: \Alg(\Fun(I, \Sp)) \to \TwCyc\Fun(I, \Sp)\,.\qed 
  \end{equation*} 
\end{proposition}
So, for an $I$-graded $\bE_1$-ring $A_*$ there are $L_p$-twisted cyclotomic structure maps
\[
\varphi_p \: L_p \THH(A_*) \longto \THH(A_*)^{tC_p} \,,  
\]
for every prime $p$. They can equivalently be described by their right adjoints $\varphi_{p;*} \:
\THH(A_*) \to \ell_p^* \THH(A_*)^{tC_p}$, having $\bT$-equivariant components
\[
\varphi_{p;i} \: \THH(A; i) \longto \THH(A; pi)^{tC_p}\,,
\]
for all $i \in I$.

If $I = \{0\}$, then an $I$-graded $L_p$-twisted cyclotomic spectrum is the same as a cyclotomic spectrum as in~\cite{NS18}*{Def.~II.1.1}, and $\TwCyc\Fun(I, \Sp)$ is the category $\Cyc\Sp$ of~\cite{NS18}*{Def.~II.1.6}.

\begin{lemma}\phantomsection\label{lem:left-right-kan-extension-vs-tw-cyc}
Let $f \: J \to I$ be a homomorphism of discrete commutative monoids.
  \begin{enumerate}[(i)]
\item Taking the left Kan extension along $f$ induces a symmetric monoidal functor $f_!\: \TwCyc\Fun(J, \Sp) \to \TwCyc\Fun(I, \Sp)$. If $B_*$ is a $J$-graded $\bE_k$-ring, then there is a natural equivalence $f_! \THH(B_*) \simeq \THH(f_! B_*)$ of $\bE_{k-1}$-algebras in $\TwCyc\Fun(I, \Sp)$. 
\item If $f$ is a face inclusion, then restriction along $f$ induces a symmetric monoidal functor  $f^* \: \TwCyc\Fun(I, \Sp) \to \TwCyc\Fun(J, \Sp)$. For an $I$-graded $\bE_k$-ring~$A_*$ there is a natural equivalence $ f^* \THH(A_*) \simeq \THH(f^*\! A_*)$ of $\bE_{k-1}$-algebras in $\TwCyc\Fun(J, \Sp)$.
\end{enumerate}
\end{lemma}
\begin{proof}
  For (i), the induced functor is that of~\cite{AMMN22}*{Lem.~A.14}. The $\bE_{k-1}$-algebra part of the equivalence follows from Corollary~\ref{cor:Bcy-on-E_k-alg}(iii) applied to the left adjoint in~\eqref{eq:restriction-extension-adjunction} and the twisted cyclotomic structure from~\cite{AMMN22}*{Cor.~A.15}.

 Part (ii) uses Corollary~\ref{cor:Bcy-on-E_k-alg}(iii) and a slight generalization of \cite{AMMN22}*{Ex.~A.17}, noting that $f^*$ is symmetric monoidal when $J$ is a face of~$I$ (compare the discussion after Definition~\ref{def:face}).
\end{proof}

\begin{corollary}
Let $A_*$ be an $I$-graded $\bE_k$-ring, with total $\bE_k$-ring $A =
c_! A_*$.  Then the total spectrum $c_! \THH(A_*) \simeq \THH(A)$
is an $\bE_{k-1}$-algebra in cyclotomic spectra. If $\{0\} \subseteq I$ is a face, then $\THH(A)_0 \simeq \THH(A_0)$ is an $\bE_{k-1}$-algebra in cyclotomic spectra. 
\end{corollary}

\begin{proof}
This follows from applying Lemma~\ref{lem:left-right-kan-extension-vs-tw-cyc} to $c\: I \to \{0\}$ and $\{0\} \subseteq I$. 
\end{proof}

\begin{definition}If $R$ is an $\bE_{\infty}$-ring, then $\THH(R)$ is an $\bE_{\infty}$-algebra in $\Cyc\Sp$. We can also view $R$ as an $I$-graded object concentrated in weight $0$, i.e., as the $I$-graded $E_{\infty}$-ring ${(\{0\}\to I)_!(R)}$, and $\THH(R)$ as an $E_{\infty}$-algebra in $\TwCyc\Fun(I, \Sp)$. We write $\Cyc\Mod_{\THH(R)}$ and $\TwCyc\Fun(I,\Mod_{\THH(R)})$ for the symmetric monoidal categories of $\THH(R)$-modules in  $\Cyc\Sp$ and in $\TwCyc\Fun(I, \Sp)$.
\end{definition}

The symmetric monoidal structure on $\THH$ implies the following statement:
\begin{corollary}
Let $A_*$ be an $I$-graded $\bE_k$ $R$-algebra with underlying $I$-graded $\bE_k$-ring~$uA_*$. Then $\THH(uA_*)$ is naturally an $\bE_{k-1}$-algebra in the symmetric monoidal category $\TwCyc\Fun(I, \Mod_{\THH(R)})$.\qed
\end{corollary}
The underlying $\bE_{k-1}$ $\THH(R)$-algebra with $\bT$-action of $\THH(uA_*)$ coincides by construction with the one resulting from Corollary~\ref{cor:Bcy-on-E_k-alg-relative}(i). 

\begin{remark} In the situation of the previous corollary, we shall generally omit `$u$' from the notation, as $\THH$ and its logarithmic and/or cyclic variants will always be applied to $\bE_k$-rings, not to $\bE_k$ $R$-algebras for $R \ne \bS$ (compare Remark~\ref{rem:relTHHcyclbase}).
\end{remark}

Next we discuss an operation that will become useful in Section~\ref{sec:multiple-generator-case}. Let $I_1$ and $I_2$ be discrete commutative monoids. Then the external product 
  \begin{multline}\label{eq:external-product}\TwCyc\Fun(I_1, \Mod_{\THH(R)}) \times \TwCyc\Fun(I_2, \Mod_{\THH(R)})\\ \longto \TwCyc\Fun(I_1\times I_2, \Mod_{\THH(R)})
  \end{multline}
is defined by first applying Lemma~\ref{lem:left-right-kan-extension-vs-tw-cyc}(i) to the inclusions of $I_1,I_2$ into $I_1\times I_2$ and then forming the symmetric monoidal product in $\TwCyc\Fun(I_1\times I_2, \Mod_{\THH(R)})$ of the resulting objects. It sends objects $X^1_*$ and $X^2_*$ to
\[X^1_* \otimes_{\THH(R)} X^2_* := (\incl\: I_1 \to I_1 \times I_2)_!(X^1_*) \otimes_{\THH(R)} (\incl\: I_2 \to I_1 \times I_2)_!(X^2_*).\]

If $A^1_*$ is an $I_1$-graded $\bE_k$ $R$-algebra and $A^2_*$ is an $I_2$-graded $\bE_k$ $R$-algebra, then we can analogously use the left adjoints from~\eqref{eq:restriction-extension-adjunction} and the monoidal structure on $\Alg_{\bE_k}(\Fun(I_1\times I_2,\Mod_R))$ to form an external product $A^1_* \otimes_R A^2_*$.
\begin{lemma}\label{lem:THH-commutes-w-external-product}
  In this situation, the universal property of $\otimes_{\THH(R)}$ induces an equivalence
  \[
    \THH(A^1_*) \otimes_{\THH(R)} \THH(A^2_*) \xrightarrow{\simeq} \THH(A^1_* \otimes_R A^2_*) \]
of $\bE_{k-1}$-algebras in $I_1\times I_2$-graded twisted cyclotomic $\THH(R)$-modules. \qed
  \end{lemma}
  
  \subsection{Log \texorpdfstring{$\THH$}{THH} of canonical prelog \texorpdfstring{$\bE_k$}{Ek}-algebras}\label{subsec:logTHH-of-canonical} Suppose given an $\bE_\infty$-ring~$R$, a discrete commutative monoid~$I$, and an $\bE_k$-algebra $\xi_*\: M_* \to c^*\Pic_R$ in the symmetric monoidal $\infty$-category $\Fun(I, \cS_{/\Pic_R})$. As in Subsection~\ref{subsec:weight-graded-Thom}, let $\tau \: {M\to I}$ be the $\bE_k$-map encoding the weight grading on the total object $\xi \: M \to \Pic_R$ of $\xi_*$. The $I$-graded $\bE_k$ $R$-algebra $\Th_R(\xi_*)$ from Subsection~\ref{subsec:weight-graded-Thom} has as total object the Thom $\bE_k$ $R$-algebra $\Th_R(\xi) = c_! \Th_R(\xi_*)$.  Its underlying $I$-graded
$\bE_k$-ring $u \Th_R(\xi)$ has $I$-graded topological Hochschild
homology $\THH(\Th_R(\xi_*))$, as in Definition~\ref{def:THHAstar},
which is an $\bE_{k-1}$-algebra in $\TwCyc\Fun(I, \Mod_{\THH(R)})$.  Similarly, the
$I^\gp$-graded $\bE_k$ $R$-algebra $\Th_R(\xi^\gp_*)$ has as total object
$\Th_R(\xi^\gp)$, and its underlying $I^\gp$-graded $\bE_k$-ring has
$I^\gp$-graded topological Hochschild homology
$\THH(\Th_R(\xi^\gp_*))$, which is an $\bE_{k-1}$-algebra in
$\TwCyc\Fun(I^\gp, \Mod_{\THH(R)})$.

\begin{remark}In examples, we will mostly consider the canonical $\pi_0(M)$-grading of Definition~\ref{def:canonical-pi0grading}. However, we allow for a general $I$-grading in our constructions because this helps studying induced maps on log $\THH$ (see Example~\ref{ex:maps-on-log-THH-w-can-grading} and Corollary~\ref{cor:Dx-comparison}) and enables us to vary the notion of repletion (compare Sections~\ref{subsec:repletion} and~\ref{subsec:comparison-rep-bar}).
\end{remark} 
\begin{definition} \label{def:logTHHcan}
The \emph{$I$-graded logarithmic topological Hochschild homology}
of the canonical prelog $\bE_k$ $R$-algebra $(\Th_R(\xi), \xi, \id)$
is the pullback
\[
\THH(\Th_R(\xi_*), \xi) := \gamma^* \THH(\Th_R(\xi^\gp_*))
= \THH(\Th_R(\xi^\gp_*)) \circ \gamma
\]
of the $I^\gp$-graded topological Hochschild homology of
$u\Th_R(\xi^\gp_*)$ along $\gamma \: I \to I^\gp$.  It is an
$\bE_{k-1}$-algebra in $\TwCyc\Fun(I, \Mod_{\THH(R)})$, with weight~$i$ component
\[
\THH(\Th_R(\xi_*), \xi; i)
	= \THH(\Th_R(\xi^\gp_*); \gamma(i))
\]
for each $i \in I$.  The $L_p$-twisted cyclotomic structure map has $\bT$-equivariant component
\[
\varphi_{p;\gamma(i)} \: \THH(\Th_R(\xi^\gp_*); \gamma(i))
        \longto \THH(\Th_R(\xi^\gp_*); p\gamma(i))^{tC_p}
\]
in weight $i$, keeping in mind that $p\gamma(i) = \gamma(pi)$.

The (ungraded) \emph{logarithmic topological Hochschild homology} of
$(\Th_R(\xi), \xi, \id)$ is
\[
\THH(\Th_R(\xi), \xi) := (c\: I \to \{0\})_! \THH(\Th_R(\xi_*), \xi) \,,
\]
which is an $\bE_{k-1}$-algebra in $\Cyc\Mod_{\THH(R)}$. 
 
Using Corollary~\ref{cor:Bcy-on-E_k-alg}(ii), the $I$-graded repletion map~\eqref{eq:repletion-on-Th} induces an \emph{$I$-graded
  repletion map}  \begin{equation}\label{eq:rho_for_canonical}\rho_* \: \THH(\Th_R(\xi_*)) \to \gamma^*\THH(\Th_R(\xi^\gp_*)) = \THH(\Th_R(\xi_*), \xi)\,,
\end{equation}
which is an $\bE_{k-1}$-algebra map in $\TwCyc\Fun(I, \Mod_{\THH(R)})$.
Its weight~$i$ component is
\begin{equation}\label{eq:rho_i_for_canonical}
\rho_i = \THH(\Th_R(\gamma); i) \: \THH(\Th_R(\xi_*); i)
	\longto \THH(\Th_R(\xi^\gp_*); \gamma(i)) \,.
\end{equation}
The corresponding map of total objects
\[
\rho \: \THH(\Th_R(\xi)) \longto \THH(\Th_R(\xi), \xi)
\]
is the (ungraded) \emph{repletion map}, which is an $\bE_{k-1}$-map
in $\Cyc\Mod_{\THH(R)}$.
\end{definition}
The repletion map~\eqref{eq:rho_for_canonical} is related to an $\bE_{k-2}$ $\THH(\Th_R(\xi_*))$-algebra structure:
\begin{corollary}\label{cor:canonical-log-thh-as-alg}
The underlying $\bE_{k-2}$-map of the repletion map~\eqref{eq:rho_for_canonical} extends to an $\bE_{k-2}$ $\THH(\Th_R(\xi_*))$-algebra structure on $\THH(\Th_R(\xi_*),\xi)$ in $\TwCyc\Fun(I, \Sp)$. 
\end{corollary}
\begin{proof}
  We can view $\xi^\gp_*$ as an $\bE_{k-1}$ $\gamma_!\xi_*$-algebra (see Proposition~\ref{prop:compatible-group-compl}), so that $u\Th_R(\xi^\gp_*)$ becomes an $\bE_{k-1}$ $u\Th_R(\gamma_!\xi_*)$-algebra by the symmetric monoidality of $\Th_R$ and~\cite{L:HA}*{Cor.~7.1.3.4}. Using $u\Th_R(\gamma_!\xi_*)\simeq \gamma_!u\Th_R(\xi_*)$, we can argue with adjunction~\eqref{eq:Ek-1-alg-adjunction-from-sym-mon-adj} as in Corollary~\ref{cor:repletion-as-alg-structure}.
\end{proof}
While not reflected in the notation, $\THH(\Th_R(\xi), \xi)$ depends on the $I$-grading of $\xi_*$ and, in general, changing the $I$-grading changes $\THH(\Th_R(\xi), \xi)$. We will always specify the grading in concrete examples. 

To study changes of the grading, let $\kappa\: I \to K$ be a homomorphism of discrete commutative monoids and consider the following commutative diagram extending~\eqref{eq:xigp}:
\begin{equation} \label{eq:xigp-with-kappa}
\xymatrix@-1pc{
K \ar[d]_-{\gamma_K}
	&& I \ar[d]_-{\gamma_I} \ar[ll]_-{\kappa}
	&& M \ar[ll]_-{\tau} \ar[rr]^-{\xi} \ar[d]^-{\gamma_M}
	&& \Pic_R \\
K^\gp
	&& I^\gp \ar[ll]_-{\kappa^\gp}
	&& M^\gp \rlap{\,.} \ar[ll]_-{\tau^\gp} \ar[rru]_-{\xi^\gp}
}
\end{equation}
Here $\tau$ specifies the grading of $\xi_* \in \Fun(I,\cS_{/\Pic_R})$, and the composite $\kappa\tau$ specifies the grading of $\kappa_!\xi_* \in  \Fun(K,\cS_{/\Pic_R})$. Analogously, $\tau^\gp$ and $\kappa^\gp\tau^\gp$ encode the gradings of $\xi_*^\gp\in \Fun(I^\gp,\cS_{/\Pic_R})$ and $\kappa_!^\gp \xi_*^\gp \simeq (\kappa_!\xi_*)^\gp\in \Fun(K^\gp,\cS_{/\Pic_R})$, respectively.
\begin{lemma}\label{lem:change-grading-of-canonical}
  In this situation, $\kappa$ induces a map of total objects 
  \[ (c_I\: I \to \{0\})_! \THH(\Th_R(\xi_*),\xi) \to (c_K\: K \to \{0\})_! \THH(\Th_R(\kappa_!\xi_*),\xi) \]
  of $\bE_{k-1}$-algebras in $\Cyc\Mod_{\THH(R)}$, which is natural in $\xi_*$. If $\kappa$ is an exact homomorphism of discrete commutative monoids, then this map is an equivalence. 
\end{lemma}
In other words, the ungraded log $\THH$ of canonical prelog $\bE_k$ $R$-algebras does not change if we modify the grading by extending along an exact homomorphism. 
\begin{proof}[Proof of Lemma~\ref{lem:change-grading-of-canonical}]
  Using that $\Th_R$ commutes with the colimit defining $\kappa_!^\gp$ and Lemma~\ref{lem:left-right-kan-extension-vs-tw-cyc}(i), we get equivalences 
  \begin{align*}
     (c_K\: K \to \{0\})_! \THH(\Th_R(\kappa_!\xi_*),\xi) 
    \simeq&\;  (c_K\: K \to \{0\})_! \gamma_K^* \THH(\Th_R(\kappa_!^\gp\xi_*^\gp))\\ \simeq&\; (c_K\: K \to \{0\})_! \gamma_K^* \kappa_!^\gp \THH(\Th_R(\xi_*^\gp))\\
    \simeq&\; \textstyle\coprod_{(k,i')\in K\times_{K^\gp}I^\gp} \THH(\Th_R(\xi_*^\gp);i')\,.
  \end{align*}
  The counit of the adjunction $(\nu_!,\nu^*)$ associated with the canonical map $\nu\: I \to K\times_{K^\gp}I^\gp$ induces a map from
\[ (c_I\: I \to \{0\})_! \THH(\Th_R(\xi_*),\xi) =  \textstyle\coprod_{i\in I} \THH(\Th_R(\xi_*^\gp);\gamma_I(i))\]
into the last term, which is an equivalence when $\kappa$ is exact (compare Remark~\ref{rem:terms-replete-exact} and~\cite{Ogu18}*{Def.~2.1.15}). \end{proof}

\begin{remark}
  When $M$ is an $\bE_k$-space and $A = \bS[M]$ is the associated spherical monoid ring, then the cyclotomic structure map $\varphi_p$ of $\THH(A)$ is determined by a $\bT$-equivariant structure map $B^\cy(M) \to (B^\cy(M))^{hC_p}$ on the space level~\cite{NS18}*{Lem.~IV.3.1}. If $\xi\: M \to \Pic_\bS$ is the constant map to basepoint, then $\Th_S(\xi) \simeq \bS[M]$. In view of this, one may wonder if for general $\xi\: M \to \Pic_\bS$ the cyclotomic structure map of $\THH(\Th_\bS(\xi))$ or of $\THH(\Th_\bS(\xi),\xi)$ arises from a similar structure on the cyclic bar construction. More generally, one can ask if these cyclotomic spectra admit a \emph{Frobenius lift} in the sense of~\cite{NS18}*{\S{}IV.3}, that is, a factorization of the $p$-cyclotomic structure map through the canonical map $(-)^{hC_p} \to (-)^{tC_p}$.

  It turns out that this is not the case: The existence of a Frobenius lift for $\THH(\bS[M])$ is related to the instability property of the (co-)homology of $\bS[M]$ as a module over the Steenrod algebra. Using~\cite{LNR12}*{(3.9) and Prop.~5.12}, one can show that for an $\bE_1$-ring $A$ with $H_*(A;\bF_p)$ of finite type, the existence of $\alpha \in H_q(A;\bF_p)$ and  $(\epsilon, j) \ne (0, 0)$ with $(\beta^\epsilon P^j)_*(\alpha) = 1$ prevent that $\varphi_p\: \THH(A) \to \THH(A)^{tC_p}$ admits a Frobenius lift. This applies for example to the Thom spectra $A = MU$ and $A = MUP_{\geq 0}$, as well as many other $\bE_1$-rings.

This also shows that the construction of log $\TC$ from~\cite{BLPO23-log-prismatic}*{\S{}3.2} for prelog structures involving only (spherical monoid rings of) discrete monoids does not generalize to our setup of prelog ring spectra. 
\end{remark}
\subsection{General log \texorpdfstring{$\THH$}{THH}}
To extend Definition~\ref{def:logTHHcan} to non-canonical prelog structures we assume that $2 \le k \le \infty$ for the rest of this section.  This ensures that $\THH(\Th_R(\xi))$ is at least an $\bE_1$-ring, so that we can form a relative tensor product over it. We also recall that $\THH(\Th_R(\xi), \xi)$ depends on an implicit $I$-grading $\tau\: M \to I$ of the domain of $\xi\: M \to \Pic_R$.

\begin{definition} \label{def:logTHHAxi}
The \emph{logarithmic topological Hochschild homology}
of an $R$-based prelog $\bE_k$-ring $(A, \xi, \bar\alpha)$ is the
$\bE_{k-2}$ $\THH(A)$-algebra in $\Cyc\Sp$ obtained as the base change 
\[
\THH(A, \xi, \bar\alpha) := \THH(A) \otimes_{\THH(\Th_R(\xi))}
	\THH(\Th_R(\xi), \xi) \,,
\]
along $\THH(\bar\alpha)$ of the $\bE_{k-2}$ $\THH(\Th_R(\xi))$-algebra $\THH(\Th_R(\xi), \xi)$ resulting from Corollary~\ref{cor:canonical-log-thh-as-alg}. We write
\begin{equation}\label{eq:repletion-relative-to-THHA}\rho \: \THH(A) \to \THH(A, \xi, \bar\alpha)
\end{equation}
for the unit of this $\bE_{k-2}$ $\THH(A)$-algebra structure and also refer to it as a \emph{repletion map}. In particular, it is both an $\bE_{k-2}$-algebra map in $\Cyc\Sp$ and a $\THH(A)$-module map in $\Cyc\Sp$. 

When $(A, \xi, \bar\alpha)$ is a prelog $\bE_k$ $R$-algebra, then this relative tensor product can be formed in $\Cyc\Mod_{\THH(R)}$ using the $\bE_{k-1}$-maps $\THH(\bar\alpha) \:
\THH(\Th_R(\xi)) \to \THH(A)$ and $\rho \: \THH(\Th_R(\xi)) \to
\THH(\Th_R(\xi), \xi)$ in that category (compare Remark~\ref{rem:equivalence-relative-tensor-products}). In this case,~\eqref{eq:repletion-relative-to-THHA} is both an $\bE_{k-2}$-algebra map in $\Cyc\Mod_{\THH(R)}$ and a $\THH(A)$-module map in $\Cyc\Mod_{\THH(R)}$. 
\end{definition}
We obtain a commutative diagram
\[
\xymatrix@-1pc{
\THH(\Th_R(\xi)) \ar[rr]^-{\THH(\bar\alpha)} \ar[d]_-{\rho}
	&& \THH(A) \ar[d]^-{\rho} \\
\THH(\Th_R(\xi), \xi) \ar[rr]
	&& \THH(A, \xi, \bar\alpha)
}
\]
in $\Cyc\Sp$ or $\Cyc\Mod_{\THH(R)}$, according to the case.
We stress that also $\THH(A, \xi, \bar\alpha)$ depends on the implicit $I$-grading $\tau\: M \to I$ of the domain of $\xi\: M \to \Pic_R$. 
This results in additional conditions for obtaining induced maps:

\begin{construction}\label{constr:maps-on-log-THH}
  Let $(f,f^\flat)\: (A, \xi\: M \to \Pic_R, \bar\alpha) \to (B,\eta\: N \to \Pic_R,\bar\beta)$ be a map of $R$-based prelog $\bE_k$-rings or of prelog $\bE_k$ $R$-algebras. Suppose $\xi$ is $I$-graded by $\tau_M \: M \to I$ and $\eta$ is $K$-graded by $\tau_N\: N \to K$. Let $\kappa \: I \to K$ be a commutative monoid homomorphism and let $f^\flat_* \: \kappa_!(\xi_*) \to \eta_*$ be an $\bE_k$-map in $\Fun(K,\cS_{/\Pic_R})$ whose underlying map of total objects in $\Alg_{\bE_k}(\cS_{/\Pic_R})$ is $f^\flat\: M \to N$. Then $(f^\flat_*,\kappa)$ induces a map $\THH(\Th_R(\xi),\xi) \to \THH(\Th_R(\eta),\eta)$ by Lemma~\ref{lem:change-grading-of-canonical}, and the compatibility of this map with $\THH(\Th_R(\xi)) \to \THH(\Th_R(\eta))$ implies that $(f,f^\flat,f^\flat_*,\kappa)$ induces a map
  \[ \THH(A,\xi,\bar\alpha) \to \THH(B,\eta,\bar\beta)\]
  of $\bE_{k-2}$-algebras in $\Cyc\Sp$ (respectively $\Cyc\Mod_{\THH(R)}$)
that is compatible with $\THH(A) \to \THH(B)$ via the repletion maps $\rho$ for the domain and codomain. 
\end{construction}
\begin{example}\label{ex:maps-on-log-THH-w-can-grading}
Given $(f,f^\flat)\: (A, \xi\: M \to \Pic_R, \bar\alpha) \to (B,\eta\: N \to \Pic_R,\bar\beta)$, we can equip $M$ with the canonical $\pi_0(M)$-grading of Definition~\ref{def:canonical-pi0grading}, $N$ with the canonical $\pi_0(N)$-grading, and set $\kappa = \pi_0(f^\flat)\: I = \pi_0(M) \to \pi_0(N) = K$. Then $f^\flat_*$ arises from the path component decomposition and we get an induced map on log $\THH$ formed with respect to the canonical $\pi_0$-gradings. 
\end{example}

We record another consequence of Lemma~\ref{lem:change-grading-of-canonical}.
\begin{corollary} Let $(A, \xi\: M \to \Pic_R, \bar\alpha)$ be an $R$-based prelog $\bE_k$-ring or a prelog $\bE_k$ $R$-algebra. Changing the grading $\tau \: M \to I$ by composing it with an exact homomorphism of discrete commutative monoids does not change log $\THH$. \qed
\end{corollary}
\begin{remark}
In the case where $R = \bS$ or $R = \bS_{(p)}$ and $k = 2$, our present definition of log $\THH$ coincides with the one used in~\cite{ABM23}*{Def.~6.6} (with the graded $\bE_2$-ring $\bS[\sigma_k]$ used in \emph{loc.\,cit.\ }corresponding to our $\Th_\bS(\xi_{2d})$). 
\end{remark}

\section{Logarithmic topological cyclic homology}\label{sec:log-TC}
In this section we use the previous constructions to define logarithmic topological cyclic homology and its relatives. We begin by recalling the basic definitions in the language of~\cite{NS18}. However, since we will only be interested in the $p$-complete topological cyclic homology of bounded below $p$-cyclotomic spectra, we deviate from their conventions by imposing $\bT$-actions on $p$-cyclotomic spectra (compare~\cite{NS18}*{Rem.~II.1.3}).

\begin{definition}[\cite{NS18}*{Prop.~II.1.9(i)}]
Let $X$ be a spectrum with $\bT$-action.  The \emph{negative topological cyclic homology} $\TC^-(X) := X^{h\bT}$ and the \emph{periodic topological cyclic homology} $\TP(X) := X^{t\bT}$ are the spectra given by the $\bT$-homotopy fixed points and the $\bT$-Tate fixed points, respectively. For any $\bE_1$-ring~$A$ we set $\TC^-(A) := \TC^-(\THH(A))$ and $\TP(A) := \TP(\THH(A))$.

Let $(X, \varphi_p)$ be a ($p$-)cyclotomic spectrum. The $p$-complete
\emph{topological cyclic homology}
\[
\TC(X)_p := \eq(G \circ \can, \varphi_p^{h\bT})
\]
is the equalizer of the maps in the diagram
\[
\xymatrix@-1pc{
\TC^-(X)_p = X^{h\bT}_p
	\ar[rr]^-{\can} \ar[drr]_-{\varphi_p^{h\bT}}
        && X^{t\bT}_p = \TP(X)_p
	\ar[d]^-{G} \\
&& (X^{tC_p})^{h\bT}_p
}
\]
where $G\: X^{t\bT} \to (X^{tC_p})^{h\bT}$ is the comparison map. For any $\bE_1$-ring~$A$ we set $\TC(A)_p = \TC(\THH(A))_p$. When the underlying spectrum of $\THH(A)$ is bounded below, $G$ is an equivalence, and, moreover, this definition of $\TC$ is equivalent to the original definition due to B\"{o}kstedt--Hsiang--Madsen~\cite{BHM93} (see~\cite{NS18}*{Thm.~1.2 and Thm.~II.3.8}). For a ($p$-)cyclotomic spectrum $(X, \varphi_p)$, we write 
\begin{equation}\label{eq:beta-TC-THH}
\beta = F \circ \pi \: \TC(X)_p
	\xrightarrow{\pi} \TC^{-}(X)_p = X^{h\bT}_p \xrightarrow{F} X_p
\end{equation}
for the composite forgetful map.

\end{definition}
The functors $\TC^-$, $\TP$ and~$\TC(-)_p$ are lax symmetric
monoidal, cf.~\cite{AMN18}*{Def.~4.4}, so if $X$ is an $\bE_{k-1}$
$\THH(R)$-algebra with $\bT$-action, then $\TC^-(X)$ is an $\bE_{k-1}$
$\TC^-(R)$-algebra and $\TP(X)$ is an $\bE_{k-1}$ $\TP(R)$-algebra, while
if $X$ is a cyclotomic $\bE_{k-1}$ $\THH(R)$-algebra, then $\TC(X)_p$
is an $\bE_{k-1}$ $\TC(R)_p$-algebra.  These observations apply, in
particular, for $X = \THH(A)$ when $A$ is an $\bE_k$ $R$-algebra.

\begin{definition}
Let $(A, \xi, \bar\alpha)$ be an $R$-based prelog $\bE_k$-ring.  Its \emph{logarithmic topological negative homology}
\[
\TC^-(A, \xi, \bar\alpha) := \THH(A, \xi, \bar\alpha)^{h\bT}
\]
is an $\bE_{k-2}$ $\TC^-(A)$-algebra, and its
\emph{logarithmic topological periodic homology}
\[
\TP(A, \xi, \bar\alpha) := \THH(A, \xi, \bar\alpha)^{t\bT}
\]
is an $\bE_{k-2}$ $\TP(A)$-algebra.
Its $p$-complete \emph{logarithmic topological cyclic homology} is the equalizer in the following diagram, and is an $\bE_{k-2}$ $\TC(A)_p$-algebra:
\[
\xymatrix@-1pc{
\TC(A, \xi, \bar\alpha)_p \ar[rr]^-{\pi}
        && \THH(A, \xi, \bar\alpha)^{h\bT}_p
        \ar[rr]^-{\can} \ar[drr]_-{\varphi_p^{h\bT}}
        && \THH(A, \xi, \bar\alpha)^{t\bT}_p
        \ar[d]^-{G} \\
&& && (\THH(A, \xi, \bar\alpha)^{tC_p})^{h\bT}_p\rlap{\,.}
}
\]
The \emph{repletion maps}
\begin{align*} \rho \: \TC^-(A) &\longto \TC^-(A, \xi, \bar\alpha), \quad
\rho \: \TP(A) \longto \TP(A, \xi, \bar\alpha), \qquad\text{and}\\  
\rho \: \TC(A)_p &\longto \TC(A, \xi, \bar\alpha)_p
\end{align*}
are induced by the map~\eqref{eq:repletion-relative-to-THHA}, and are in particular $\bE_{k-2}$-ring maps and left module maps in the respective module categories.

If $\bar\alpha \: \Th_R(\xi) \to A$ is an $\bE_k$ $R$-algebra map, so that $(A, \xi, \bar\alpha)$ is a prelog $\bE_k$ $R$-algebra, then the same constructions yield $\bE_{k-2}$-algebras and -maps over $\TC^-(R)$, $\TP(R)$ and~$\TC(R)_p$, respectively.
\end{definition}
\begin{example}
  We can in particular apply the previous definition to the canonical prelog $\bE_k$ $R$-algebra $(\Th_R(\xi), \xi, \id)$ to obtain $\bE_{k-1}$-algebras $\TC^-(\Th_R(\xi), \xi)$, $\TP(\Th_R(\xi), \xi)$, and $\TC(\Th_R(\xi), \xi)_p$ over  $\TC^-(R)$, $\TP(R)$ and~$\TC(R)_p$, respectively. This uses the $\bE_{k-1}$-structure on $\THH(\Th_R(\xi),\xi)$ from Definition~\ref{def:logTHHcan}, rather than the $\bE_{k-2}$-structure from its $\THH(\Th_R(\xi))$-algebra structure (compare Corollary~\ref{cor:canonical-log-thh-as-alg}). 
\end{example}

\subsection{Graded log \texorpdfstring{$\TC$}{TC}}
We now introduce parts of $\TC$ that can be built from graded $\THH$, and consider when the whole of~$\TC$ can be recovered from these pieces.

\begin{definition}
Given a discrete commutative monoid~$I$ let $I/{\sim_p}$ be the set of equivalence classes~$[i]_p$ with respect to the equivalence relation on the underlying set of~$I$ generated by $i \sim_p pi$ for all $i \in I$.
\end{definition}

For example, the elements of~$\bZ_{\ge0}/{\sim_p}$ are~$[0]_p$ together with the $[i]_p$ for $i \in \bZ_{>0}$ with $p \nmid i$.
\begin{definition}[\cite{AMMN22}*{Rem.~5.39}]
Let $X_*$ be an $I$-graded spectrum with $\bT$-action.  The \emph{$I$-graded topological negative homology} $\TC^-(X_*) := (X_*)^{h\bT}$ and the \emph{$I$-graded topological periodic homology} $\TP(X_*) := (X_*)^{t\bT}$ are the $I$-graded spectra given by the $\bT$-homotopy fixed points $(X_i)^{h\bT}$ and $\bT$-Tate fixed points $(X_i)^{t\bT}$, respectively, in each weight $i \in I$.  For any $I$-graded $\bE_1$-ring~$A_*$ we set $\TC^-(A_*) := \TC^-(\THH(A_*))$ and $\TP(A_*) := \TP(\THH(A_*))$.
\end{definition}

The $I$-graded functors $\TC^-$ and $\TP$ are lax symmetric monoidal, so if $X_*$ is an $I$-graded $\bE_{k-1}$ $\THH(R)$-algebra with $\bT$-action, then $\TC^-(X_*)$ is an $I$-graded $\bE_{k-1}$ $\TC^-(R)$-algebra and $\TP(X_*)$ is an $I$-graded $\bE_{k-1}$ $\TP(R)$-algebra.  This applies, in particular, for $X_* = \THH(A_*)$ when $A_*$ is an $I$-graded $\bE_k$ $R$-algebra.

\begin{remark}
Let $c \: I \to \{0\}$ and $X = c_! X_*$.
Beware that the canonical maps
\begin{align*}
\textstyle\bigvee_{i \in I} \TC^-(X_i) = c_! \TC^-(X_*)
	&\longto \TC^-(c_! X_*) = \TC^-(X) \\
\textstyle\bigvee_{i \in I} \TP(X_i) = c_! \TP(X_*)
	&\longto \TP(c_! X_*) = \TP(X)
\end{align*}
are generally not equivalences, since the colimit defining total objects
does not generally commute with the limits defining $\TC^-$ and $\TP$.
\end{remark}

\begin{definition}
Let $(X_*, \varphi_p)$ be an $I$-graded $L_p$-twisted cyclotomic spectrum.
For each $\sim_p$-equivalence class $[i]_p \in I/{\sim_p}$, let
\[
\textstyle X_{[i]_p} = \bigvee_{j \in [i]_p} X_j \,.
\]
The composite map
\[
\textstyle \varphi_p \: X_{[i]_p} = \bigvee_{j \in [i]_p} X_j
\xrightarrow{\bigvee_{j \in [i]_p} \varphi_{p;j}}
\bigvee_{j \in [i]_p} \bigl( X_j^{tC_p} \bigr)
\longto \bigl( \bigvee_{j \in [i]_p} X_j \bigr)^{tC_p}
	= X_{[i]_p}^{tC_p}
\]       
makes $X_{[i]_p}$ a cyclotomic subspectrum of $X = c_! X_*$.  The
$p$-completed \emph{$I/{\sim_p}$-graded topological cyclic homology}
is the functor
\[
[i]_p \longmapsto \TC(X_{[i]_p})_p \,.
\]
\end{definition}
\begin{definition} \label{def:propconn}
A spectrum~$Y$ is $m$-connective if $\pi_*(Y) = 0$ for all $* < m$.
Let $X_*$ be an $I$-graded spectrum.  We say that $X_*$ is \emph{properly
connective} if for each $m \in \bZ$ the set
\[
\{ i \in I \mid \text{$X_i$ is $m$-connective} \}
\]
has finite complement in~$I$.
\end{definition}

\begin{example}
For $R$ connective and $\xi\: \<x\> \to \Pic_R$, the $\bZ_{\ge0}$-graded spectra $X_*
= \THH(\Th_R(\xi_*))$ and $\THH(\Th_R(\xi_*), \xi)$ will be properly
connective if $\xi(x)$ is at least $1$-connective, since then
$\xi(x^i) = \xi(x)^{\otimes_R i}$ will be at least $i$-connective.  They will
typically not be properly connective if $\xi(x)$ is $0$-connective,
such as for $\xi = 0$ appearing in the examples  $(ku,\< p \>)$,  $(\kup,\< p \>)$, and $(\ell,\< p\>)$ from~\eqref{eq:ell<p>ku<p>}.

Similarly, for $\xi\: \<x, y\> \to \Pic_R$ the corresponding $(\bZ_{\ge0})^2$-graded spectra will be properly connective if $\xi(x)$ and $\xi(y)$ are both
$1$-connective, but usually not if one of them is just $0$-connective as in the examples $(\mathrm{BP}\< n\>, \< p,v_1,\dots,v_n\>)$ from~\eqref{eq:BPnp-to-vn} or $(\ell,\< p, v_1\>)$ from~\eqref{eq:ell_p_v1}.   
\end{example}

\begin{lemma} \label{lem:bigvee-prod}
Let $X_*$ be an $I$-graded spectrum with $\bT$-action.  Suppose that $X_*$
is properly connective.
\begin{enumerate}[(i)]
  \item 
The canonical map
\[
\textstyle X = \bigvee_{i \in I} X_i \xrightarrow{\simeq} \prod_{i \in I} X_i
\]
is an equivalence.
\item 
The middle and right-hand canonical maps
\[
\textstyle \bigvee_{i \in I} (X_i^{tC_p})
	\longto \bigl( \bigvee_{i \in I} X_i \bigr)^{tC_p}
	\xrightarrow{\simeq} \bigl( \prod_{i \in I} X_i \bigr)^{tC_p}
	\xrightarrow{\simeq} \prod_{i \in I} (X_i^{tC_p})
\]
are equivalences.

\item  The middle and right-hand canonical maps
\[
\textstyle \bigvee_{i \in I} (X_i^{t\bT})
	\longto \bigl( \bigvee_{i \in I} X_i \bigr)^{t\bT}
	\xrightarrow{\simeq} \bigl( \prod_{i \in I} X_i \bigr)^{t\bT}
	\xrightarrow{\simeq} \prod_{i \in I} (X_i^{t\bT})
\]
are equivalences.
\end{enumerate}
\end{lemma}

\begin{proof}
The first displayed map induces the inclusion $\bigoplus_{i \in I}
\pi_*(X_i) \to \prod_{i \in I} \pi_*(X_i)$ in homotopy, and we are
assuming that in each fixed degree~$*$ there are only finitely many
nonzero factors in this product.  This implies that the map in~(i)
is an equivalence.  It follows that the middle maps in the second and
third displays are also equivalences.

To show that the right-hand map in~(ii) is an equivalence, the map of
norm--restriction cofiber sequences
\[
\xymatrix@-1pc{
\bigl( \prod_{i \in I} X_i \bigr)_{hC_p} \ar[rr]^-{N^h} \ar[d]
&& \bigl( \prod_{i \in I} X_i \bigr)^{hC_p} \ar[rr]^-{R^h} \ar[d]^-{\simeq}
&& \bigl( \prod_{i \in I} X_i \bigr)^{tC_p} \ar[d] \\
\prod_{i \in I} ((X_i)_{hC_p}) \ar[rr]^-{N^h}
&& \prod_{i \in I} (X_i^{hC_p}) \ar[rr]^-{R^h}
&& \prod_{i \in I} (X_i^{tC_p})
}
\]
shows that it suffices to prove that the left-hand vertical map is
an equivalence.  Forming homotopy orbits preserves connectivity,
so for each fixed degree~$*$ only finitely many of the factors $X_i$
contribute to $\pi_*$ of the two left-hand spectra.  The claim then
follows, because $(\bigvee_{i \in I} X_i)_{hC_p} \to \bigvee_{i \in I}
((X_i)_{hC_p})$ is an equivalence.

The proof that the right-hand map in~(iii) is an equivalence is
similar, using the norm--restriction cofiber sequence
\begin{equation} \label{eq:T-NRseq}
\xymatrix{
\Sigma (X_i)_{h\bT} \ar[r]^-{N^h} & X_i^{h\bT} \ar[r]^-{R^h} & X_i^{t\bT}
}
\end{equation}
and its analogue for $\prod_{i \in I} X_i$.
\end{proof}

\begin{proposition} \label{prop:TC-prodIp}
Let $(X_*, \varphi_p)$ be an $I$-graded $L_p$-twisted cyclotomic spectrum.
The $p$-cyclotomic structure map for $X = c_! X_*$ factors
as
\[
\textstyle{}X = \bigvee_{i \in I} X_i
	\xrightarrow{\bigvee_{i \in I} \varphi_{p;i}}
	\bigvee_{i \in I} (X_i^{tC_p})
	\longto \bigl( \bigvee_{i \in I} X_i \bigr)^{tC_p} \,.
\]
If $X_*$ is properly connective, then this composite is equivalent to
the product map
\[
\textstyle{}\prod_{i \in I} X_i
	\xrightarrow{\prod_{i \in I} \varphi_{p;i}}
	\prod_{i \in I} (X_i^{tC_p}) \,,
\]
and there is a natural equivalence
\[
\textstyle{}\TC(X)_p \simeq \prod_{[i]_p \in I/{\sim_p}} \TC(X_{[i]_p})_p \,.
\]
\end{proposition}

\begin{proof}
The first claims are immediate from Lemma~\ref{lem:bigvee-prod}.
Moreover, there are canonical equivalences
\begin{align*}
\textstyle{}X^{h\bT} = \bigl( \bigvee_{i\in I} X_i \bigr)^{h\bT}
        &\xrightarrow{\simeq} \bigl( \prod_{i\in I} X_i \bigr)^{h\bT}
        \simeq \prod_{i\in I} \bigl( X_i^{h\bT} \bigr) \\
\textstyle{}X^{t\bT} = \bigl( \bigvee_{i\in I} X_i \bigr)^{t\bT}
        &\xrightarrow{\simeq} \bigl( \prod_{i\in I} X_i \bigr)^{t\bT}
        \simeq \prod_{i\in I} \bigl( X_i^{t\bT} \bigr) \\
\textstyle{}(X^{tC_p})^{h\bT} = \bigl( \bigl( \bigvee_{i\in I} X_i \bigr)^{tC_p} \bigr)^{h\bT}
        &\xrightarrow{\simeq} \bigl( \bigl( \prod_{i\in I} X_i \bigr)^{tC_p} \bigr)^{h\bT}
        \simeq \prod_{i\in I} \bigl( \bigl( X_i^{tC_p} \bigr)^{h\bT} \bigr)
\,.
\end{align*}
Under these identifications, the maps $\can \: X^{h\bT} \to X^{t\bT}$
and $G \: X^{t\bT} \to (X^{tC_p})^{h\bT}$ correspond to the product
over $i\in I$ of the maps
\[
\can_i \: X_i^{h\bT} \longto X_i^{t\bT}
\qqandqq
G_i \: X_i^{t\bT} \longto (X_i^{tC_p})^{h\bT} \,,
\]
while $\varphi_p^{h\bT} \: X^{h\bT} \to (X^{tC_p})^{h\bT}$ corresponds
to the product over~$i\in I$ of the maps
\[
\varphi_{p;i}^{h\bT} \: X_i^{h\bT}
        \longto \bigl( X_{pi}^{tC_p} \bigr)^{h\bT} \,.
\]
Hence $\TC(X)_p$ is the product over the equivalence classes~$[i]_p \in
I/{\sim_p}$ of the equalizer of the two maps
\begin{equation} \label{eq:TCXipp}
\xymatrix@C+2pc{
\prod_{j \in [i]_p} (X_j)^{h\bT}_p
	\ar@<1ex>[r]^-{\prod_j G_j \circ \can_j}
	\ar@<-1ex>[r]_-{\prod_j \varphi_{p;j}^{h\bT}}
	& \prod_{j \in [i]_p} (X_j^{tC_p})^{h\bT}_p \,.
}
\end{equation}
Since $X_*$ restricted to $[i]_p \subset I$ remains properly connective,
this factor in the product is also the equalizer~$\TC(X_{[i]_p})_p$
of the two maps
\[
\xymatrix@C+2pc{
(X_{[i]_p})^{h\bT}_p =
\bigl( \bigvee_{j \in [i]_p} X_j \bigr)^{h\bT}_p
        \ar@<1ex>[r]^-{G \circ \can}
        \ar@<-1ex>[r]_-{\varphi_p^{h\bT}}
        &
        \bigl( \bigl( \bigvee_{j \in [i]_p} X_j \bigr)^{tC_p} \bigr)^{h\bT}_p
	= (X_{[i]_p}^{tC_p})^{h\bT}_p \,.
}
\]
\end{proof}

\section{Topological Hochschild homology of Thom rings and algebras}\label{sec:THH-of-Thom-rings}
The goal of this section is to show that the topological Hochschild homology spectra of (weight-graded) Thom rings and algebras are equivalent to (weight-graded) Thom spectra of suitable cyclic bar constructions, and that the log $\THH$ spectra of canonical prelog algebras are equivalent to Thom spectra of suitable replete bar constructions. These will be equivalences of $\bE_{k-1}$-rings or $\bE_{k-1}$ $\THH(R)$-algebras with $\bT$-action, not generally fully capturing the cyclotomic structures constructed in Section~\ref{sec:logTHH}.

\subsection{\texorpdfstring{$\THH$}{THH} of Thom rings}
We first handle the case $R = \bS$. Let $I$ be a discrete commutative monoid. As in Subsection~\ref{subsec:weight-graded-Thom}, we consider the symmetric monoidal category $\Fun(I,\cS_{/\Pic_\bS})$ with the Day convolution product structure. Then Corollary~\ref{cor:Bcy-on-E_k-alg}(i) provides a functor
  \begin{equation}\label{eq:graded-Bcy-composite}
 B^\cy\: \Alg_{\bE_k}(\Fun(I, \cS_{/\Pic_\bS})) \to \Alg_{\bE_{k-1}}(\Fun(I, \cS_{/ \Pic_\bS})^{B\bT}).
\end{equation}
Using Lemma~\ref{lem:functors-to-slice-categories}, we note that it sends an $I$-graded $\bE_k$-algebra $\xi_*\: M_* \to c^*\Pic_\bS$ in $\cS_{/\Pic_\bS}$ to the $I$-graded  $\bE_{k-1}$-algebra with $\bT$-action $B^\cy(M_*) = \{i \mapsto B^\cy(M_*;i)\}$ with augmentation $\xi^\cy_*\: B^\cy(M_*) \to c^*\Pic_\bS$, where $c^*\Pic_\bS$ is the constant $I$-diagram on $\Pic_\bS$ with trivial $\bT$-action. The map $\xi^\cy_*$ is the composite of $B^\cy(\xi_*)\: B^\cy(M_*) \to B^\cy(c^*\Pic_\bS)$ and the map $B^\cy(c^*\Pic_\bS) \to c^*\Pic_\bS$ induced by the iterated multiplications of the $\bE_{\infty}$-object $c^*\Pic_\bS$. When $I = \{0\}$, this specializes to an absolute cyclic bar construction sending an $\bE_k$-map $\xi\: M \to \Pic_\bS$ to an $\bE_{k-1}$-algebra $\xi^\cy\: B^\cy(M) \to \Pic_\bS$ in $(\cS_{/ \Pic_\bS})^{B\bT}$. 

We record some compatibility statements, of which part~\eqref{eq:ungraded-THH-Bcy-over-S} is closely related to \cite{BCS10}*{Thm.~2.1}:

\begin{proposition}\label{prop:Bcy-THH-comparison-over-PicS} Let $c\: I \to \{0\}$ be the collapse map and let $\xi_*\: M_* \to c^*\Pic_\bS$ be an $\bE_k$-algebra in $\Fun(I, \cS_{/\Pic_\bS})$ with total object $\xi\: M \to \Pic_\bS$.
\begin{enumerate}[(i)]
\item There is a natural equivalence of $I$-graded $\bE_{k-1}$-rings with $\bT$-action
\begin{equation}\label{eq:graded-THH-Bcy-over-S}
\THH(\Th_\bS(\xi_*)) \simeq \Th_\bS(\xi^\cy_* \: B^\cy(M_*) \to c^*\Pic_\bS)
\end{equation}
and an analogous equivalence of total objects
\begin{equation}\label{eq:ungraded-THH-Bcy-over-S}
\THH(\Th_\bS(\xi)) \simeq \Th_\bS(\xi^\cy \: B^\cy(M) \to \Pic_\bS)
\,.
\end{equation}
\item There is a natural equivalence of $\bE_{k-1}$-algebras with $\bT$-action in  $\cS_{/ \Pic_{\bS}}$ between
  \[ c_!(\xi^\cy_*\: B^\cy(M_*) \to c^*\Pic_\bS) \qqandqq \xi^\cy\: B^\cy(M) \to \Pic_\bS. \]
\item Passing to total objects, the equivalence~\eqref{eq:graded-THH-Bcy-over-S} recovers~\eqref{eq:ungraded-THH-Bcy-over-S}.
\end{enumerate}
\end{proposition}
\begin{proof}
All statements follow from Corollary~\ref{cor:Bcy-on-E_k-alg}(iii) applied to suitable parts of the commutative square 
\[\xymatrix@-1pc{\Fun(I,\cS_{/ \Pic_\bS}) \ar[rr]^{\Th_\bS} \ar[d]_{c_!} && \Fun(I,\Sp)\ar[d]^{c_!} \\ 
    \cS_{/ \Pic_\bS} \ar[rr]^{\Th_\bS}  && \Sp}
\]
of colimit-preserving symmetric monoidal functors.
\end{proof}

\subsection{\texorpdfstring{$\THH$}{THH} of Thom algebras}
To analyze~$\THH$ in the case of Thom $R$-algebras for a general base $\bE_\infty$-ring~$R$, we must take into account that the forgetful functor $u \: \Mod_R \to \Sp$ is usually only lax symmetric monoidal. This implies that the statement of Proposition~\ref{prop:Bcy-THH-comparison-over-PicS}(i) does in general not hold with $R$ in place of $\bS$.  To alleviate this problem, we pass to a more general setting where the $\bE_\infty$-ring $R$ is allowed to vary.

\begin{remark} \label{rem:relTHHcyclbase}
If we were concerned with the \emph{relative} topological Hochschild homology $\THH^R(A) = | [q] \mapsto A^{\otimes_R 1+q}|$ of an $R$-algebra~$A$, then we could work in $\Mod_R$ throughout.  However, $\THH^R(A)$ will only admit a cyclotomic structure under additional hypotheses on~$R$, such as the existence of a $\bT$-map $R \to R^{tC_p}$ factoring $\epsilon^{tC_p} \varphi_p \: \THH(R) \to \THH(R)^{tC_p} \to R^{tC_p}$ over $\epsilon \: \THH(R) \to R$, i.e., that~$R$ is a \emph{cyclotomic base} (compare~\cite{AMMN22}*{Ex.~A.3}, \cite{HRW}*{Def.~3.2.1},~\cite{BMY23}, and work in progress by Devalapurkar, Hahn, Raksit, and Yuan in the $R = MU$ case). 
\end{remark}
\begin{lemma}\label{lem:Pic-lax-sym-mon}
The functor $\Pic \: \CAlg(\Sp) \to \cS$ is lax symmetric monoidal and preserves filtered colimits.
\end{lemma}
Informally, for $\bE_\infty$-rings $R^1, R^2 \in \CAlg(\Sp)$ the lax
structure map
\[
\hat\otimes \: \Pic_{R^1} \times \Pic_{R^2}
	\longto \Pic_{R^1 \otimes R^2}
\]
takes an invertible $R^1$-module $P^1$ and an invertible $R^2$-module
$P^2$ to the invertible $R^1 \otimes R^2$-module $P^1 \otimes P^2$.
\begin{proof}[Proof of Lemma~\ref{lem:Pic-lax-sym-mon}]
This pairing is obtained from the composite
\[
\Mod_{R^1} \times \Mod_{R^2}
\longto \Mod_{R^1} \otimes \Mod_{R^2}
	\simeq \Mod_{R^1 \otimes R^2}
\]
from \cite{L:HA}*{Prop.~4.8.1.17} and \cite{L:HA}*{Thm.~4.8.5.16} by
passage to the maximal grouplike (with respect to $\otimes$) $\infty$-groupoids inside these
symmetric monoidal $\infty$-categories.

The functor $R \mapsto \Mod_R$ preserves $K$-indexed colimits for weakly contractible simplicial sets~$K$ by \cite{L:HA}*{Cor.~4.8.5.13}, and this includes all filtered~$K$ by \cite{L:HTT}*{Lem.~5.3.1.20}.  Moreover, $\Mod_R \mapsto \Pic_R$ preserves filtered colimits, since passage to the groupoid core $\Mod_R^\simeq$, and to the full subcategory of invertible objects, both have this property.
\end{proof}

\begin{definition}\label{def:FunIcsovercstarPic}
Let $I$ be a discrete commutative monoid. We define 
\[
\Fun(I,\cS)_{/c^*\Pic} = \CAlg(\Sp) \times_{\Fun(\{1\}, \Fun(I,\cS))} \Fun(\Delta^1, \Fun(I,\cS)) 
\]
to be $\infty$-category with objects pairs $(R, \xi_*\: M_* \to c^*\Pic_R)$ where $c^*$ is the restriction along $c\: I \to \{0\}$. To equip it with a symmetric monoidal structure, we define $(\Fun(I,\cS)_{/c^*\Pic})^\otimes$ as the pullback in $\infty$-operads of
\[ \CAlg(\Sp)^\otimes \xrightarrow{(c^*\Pic)^\otimes}  \Fun(I,\cS)^\otimes \xleftarrow{\ev_1^\otimes} \Fun(\Delta^1, \Fun(I,\cS)^\otimes), \] where
$\ev_1^\otimes$ is a coCartesian fibration
by~\cite{L:HTT}*{Cor.~2.4.7.12}.  The composite of its pullback
$(\Fun(I,\cS)_{/c^*\Pic})^\otimes\to \CAlg(\Sp)^\otimes$ with the coCartesian
fibration $\CAlg(\Sp)^\otimes \to N(\Fin_*)$ that defines the
symmetric monoidal structure of~$\CAlg(\Sp)$ is the coCartesian
fibration that exhibits the symmetric monoidal structure of
$\Fun(I,\cS)_{/c^*\Pic}$.
\end{definition}

Informally, the symmetric monoidal pairing of $\Fun(I,\cS)_{/c^*\Pic}$ takes two objects 
$(R^1, \xi_*^1 \: M_*^1 \to c^*\Pic_{R^1})$ and
$(R^2, \xi_*^2 \: M_*^2 \to c^*\Pic_{R^2})$ to
\[
  (R^1 \otimes R^2, \bar\otimes (\xi^1_* \times \xi^2_*) \: M^1_* \times M^2_*
  \to c^*\Pic_{R^1 \otimes R^2})\,,
\]
where $\bar\otimes$ is the lax structure map for the functor $c^*\Pic$ (that results as the composite of the lax structure map $\hat\otimes$ of $\Pic$ with that of $c^*$).
\begin{lemma}
The $\infty$-category $\Fun(I,\cS)_{/c^*\Pic}$ is presentably symmetric monoidal. 
\end{lemma}
\begin{proof}
We consider the functors $F,G\: \Fun(I,\cS) \times \CAlg(\Sp) \to \Fun(I,\cS)$ given by the projection to the first factor and the projection to the second factor followed by $c^*\Pic$, respectively. 
Then by construction, $\Fun(I,\cS)_{/c^*\Pic}$ is equivalent to the lax equalizer $\mathrm{LEq}(F,G)$ in the sense of~\cite{NS18}*{Def.~II.1.4}. Since $F$ preserves colimits and $G$ preserves filtered colimits by Lemma~\ref{lem:Pic-lax-sym-mon},~\cite{NS18}*{Prop.~II.1.5(iv)} and its proof show that $\Fun(I,\cS)_{/c^*\Pic}$ is presentable and that the projection $\mathrm{LEq}(F,G) \to \Fun(I,\cS) \times \CAlg(\Sp)$ preserves and detects colimits. Using that  $\CAlg(\Sp)$ and $\Fun(I,\cS)$ are presentably symmetric monoidal, it follows from the construction of the symmetric monoidal product of $\Fun(I,\cS)_{/c^*\Pic}$ that it preserves colimits separately in each variable.
\end{proof}

We define the \emph{underlying Thom spectrum} functor to be
\begin{equation}\label{eq:uTh_onFunIS_cPic}
u\Th \: \Fun(I,\cS)_{/c^*\Pic} \to \Fun(I,\Sp), \quad (R, \xi_*\: M_* \to c^*\Pic_R) \mapsto u \Th_R(\xi_*). 
\end{equation}
When $I = \{0\}$, this specializes to a functor $u\Th \: \cS_{/\Pic} \to \Sp$. 

\begin{lemma} \label{lem:uThsifted}
The functor~\eqref{eq:uTh_onFunIS_cPic} is symmetric monoidal and
preserves sifted colimits.
\end{lemma}
\begin{proof}Recall that $u\Th_R(\xi_*) \simeq \{i \mapsto \colim(u\iota\xi_i\: M_i \to \Mod_R)\}$. We have 
\begin{align*}
&\;u \Th_{R_1}(\xi_*^1) \otimes u \Th_{R_2}(\xi_*^2)\\
&\;=\; \{i\mapsto \textstyle\coprod_{i_1+i_2=i}\colim (u \iota \xi_{i_1}^1 \: M_{i_1}^1 \to \Sp)	\otimes \colim (u \iota \xi_{i_2}^2 \: M_{i_2}^2 \to \Sp)\} \\
&\;\simeq\; \{i\mapsto \colim (\textstyle\coprod_{i_1+i_2=i} M_{i_1}^1 \times M_{i_2}^2 \xrightarrow{u\iota\xi^1_{i_1} \times u\iota\xi^2_{i_2}} \Sp \times \Sp
	\xrightarrow{\otimes} \Sp)\} \\
&\;= \;\colim (u \iota \bar\otimes (\xi^1_* \times \xi^2_*)
	\: M^1_* \times M^2_* \to \Sp) \\
&\; =\; u \Th_{R_1 \otimes R_2}(\bar\otimes (\xi^1_* \times \xi^2_*))\,,
\end{align*}
since the smash product of spectra preserves colimits in each variable.

For the rest of this proof, we will follow the terminology of~\cites{L:HTT, L21} and say that a functor between $\infty$-categories is \emph{cofinal} if restriction along it preserves colimits. As a first preparatory step, let $\cD$ be an $\infty$-category and let $\Ar(\cD) = \Fun(\Delta^1, \cD)$ denote the category of arrows $d \to e$ in~$\cD$. We claim that the $\infty$-functor $\cD \to \Ar(\cD)$ sending $d$ to $\id_d$ is cofinal. By the $\infty$-categorical version of Quillen's Theorem A (see e.g.~\cite{L21}*{Thm. 4.4.20}), it is sufficient to show that for any object $\alpha \: d \to e$ in $\Ar(\cD)$, the slice $\cD_{/\alpha}$ is weakly contractible. This slice is weakly equivalent to the slice $\cD_{e/}$, which is weakly contractible because it has an initial object.

Let $\cD$ be a sifted $\infty$-category. As a second preparatory step, we claim that for any object $d$ of $\cD$, the forgetful functor $U\: \cD_{d/} \to \cD$ is cofinal. Arguing as above, it is sufficient to show that for any object $e$ of $\cD$ the slice $(\cD_{d/})_{e/}$ is weakly contractible. Identifying the latter with the slice of the diagonal functor $\cD \to \cD\times \cD$ under $(d,e)$, the conclusion follows from the definition of a sifted $\infty$-category.

Let  $d \mapsto (R_d, (\xi_*)_d \: (M_*)_d \to c^*\Pic_{R_d})$ be a $\cD$-diagram in $\Fun(I,\cS)_{/c^*\Pic}$, with colimit $(R, \xi_*\: M_* \to c^*\Pic_R)$. We can extend the functor $d \mapsto u \Th_{R_d}((\xi_*)_d)$ over $\Ar(\cD)$ by sending $\alpha \: d \to e$ to $u
\Th_{R_e}((c^*\Pic_\alpha) \circ (\xi_*)_d)$, where $\Pic_\alpha \: \Pic_{R_d}
\to \Pic_{R_e}$ is given by base change along $\alpha \: R_d \to R_e$.
We then obtain equivalences
\begin{align*}
\colim_{d \in \cD} \ u \Th_{R_d}((\xi_*)_d)
	&\simeq \colim_{\alpha\: d \to e \in \Ar(\cD)}\ %
		u \Th_{R_e}((c^*\Pic_\alpha)\circ(\xi_*)_d) \\
	&\simeq \colim_{d \in \cD} \ %
		\colim_{\alpha \: d \to e \in \cD_{d/}} \ %
		u \Th_{R_e}((c^*\Pic_\alpha) \circ (\xi_*)_d) \\
	&\simeq \colim_{d \in \cD} \ u \Th_R((c^*\Pic_{R_d \to R}) \circ (\xi_*)_d) \\
	&\simeq u \Th_R(\xi_*) \,.
\end{align*}
The first equivalence arises from the first preparatory step and the last one from $u\Th_R$ commuting with colimits.  The second one arises as in~\cite{HV92}*{Prop.~6.2}, viewing the arrow category $\Ar(\cD)$ as the Grothendieck construction on the undercategories $\cD_{d/}$ for $d \in \cD$. For the third equivalence, the map we need to show to be an equivalence is induced by the maps $u(R_e\otimes_{R_d}\Th_{R_d}((\xi_*)_d)) \to u(R\otimes_{R_d}\Th_{R_d}((\xi_*)_d)$ 
via the diagram of $R$-modules underlying the restriction of a colimit cone for the equivalence of $\bE_\infty$-rings $\colim_{d \to e \in \cD_{d/}}R_e \to R$ resulting from the second preparatory step. By writing $(M_*)_d$ as a colimit of its points, this reduces to the case where $\Th_{R_d}((\xi_*)_d)$ is $\otimes$-invertible over $R_d$. So we need to show that $\colim_{d \to e \in \cD_{d/}}uR_e \to uR$ is an equivalence. Because the colimit system arises by restriction along the projection $\cD_{d/} \to \cD$, we can apply the second preparatory step once more so that the claim reduces to showing that $\colim_{e \in \cD}uR_e \to uR$ is an equivalence. This holds because sifted colimits of $\bE_\infty$-rings are sifted colimits of the underlying spectra~\cite{L:HA}*{Cor.~3.2.3.2}.
\end{proof}

The passage to total objects $c_! \: \Fun(I,\cS) \to \cS$ induces a functor 
\begin{equation}\label{eq:cshriek_onFunIS_cPic}
c_! \: \Fun(I,\cS)_{/c^*\Pic} \to \cS_{/\Pic}
\end{equation}
sending $(R, \xi_*\: M_* \to c^*\Pic_R)$ to the adjoint $\xi\: M = c_!(M_*) \to \Pic_R$ of $\xi_*$. 
\begin{lemma} 
The functor~\eqref{eq:cshriek_onFunIS_cPic} is symmetric monoidal and preserves colimits.
\end{lemma}
\begin{proof}
The functor preserves colimits because it has a right adjoint induced by~$c^*$. It is symmetric monoidal because the underlying total object of the Day convolution product in $\Fun(\cI,\cS)$ is equivalent to the product of the underlying total objects.  
\end{proof}
We also need the category $(\Fun(I,\cS)_{/c^*\Pic})^{B\bT}$ of $\bT$-objects in $\Fun(I,\cS)_{/c^*\Pic}$.  Since $\Fun(\bT,-)$ commutes with pullbacks, Definition~\ref{def:FunIcsovercstarPic} implies that such a $\bT$-object consists of an $\bE_\infty$-ring with $\bT$-action and $\bT$-map in $\Fun(I,\cS)$ to $c^*\Pic$ of this $\bE_\infty$-ring with $\bT$-action. We extend $u\Th$ to this category of $\bT$-objects by simply appealing to functoriality.

Applying Corollary~\ref{cor:Bcy-on-E_k-alg}(i) provides a cyclic bar construction  
\[  B^{\cy}\: \Alg_{\bE_k}(\Fun(I,\cS)_{/c^*\Pic}) \to \Alg_{\bE_{k-1}}((\Fun(I,\cS)_{/c^*\Pic})^{B\bT})\, \]
which can be described as follows:

\begin{lemma}\label{lem:Bcy-on-FunIS_cPic}
For an object $(R,\xi_*\: M_* \to c^*\Pic_R)$ of $\Alg_{\bE_k}(\Fun(I,\cS)_{/c^*\Pic})$, the cyclic bar construction $B^\cy(R,\xi_*)$ is, as an object of $\Alg_{\bE_{k-1}}((\Fun(I,\cS)_{/c^*\Pic})^{B\bT})$, equivalent to
  \[ (\THH(R), \xi_*^\cy := \bar\epsilon \circ B^\cy(\xi_*) \: B^\cy(M_*)
	\longto B^\cy(c^*\Pic_R) \longto c^*\Pic_{\THH(R)}) \,,\] 
      where $\bar\epsilon \: B^\cy(c^*\Pic_R) \to c^*\Pic_{\THH(R)}$ is the $\bT$-equivariant $\bE_\infty$-map obtained by applying Corollary~\ref{cor:Bcy-on-E_k-alg}(ii) to the lax symmetric monoidal structure of $c^*\Pic$.
\end{lemma}
\begin{proof}
The forgetful functors from $\Fun(I,\cS)_{/c^*\Pic}$ to $\CAlg(\Sp)$ and~$\Fun(I,\cS)$ are symmetric monoidal and create the colimit in $\Fun(I,\cS)_{/c^*\Pic}$ defining geometric realization, using the lax symmetric monoidal structures on $c^*$ and~$\Pic$.  This explains the appearance of $\THH(R)$, $B^\cy(M_*)$, and $\xi^\cy_*$, respectively. 
\end{proof}

An $\bE_k$-algebra $\xi_*\: M_* \to c^*\Pic_R$ in  $\Fun(I,\cS_{/\Pic_R})$ can be viewed as an $\bE_k$-algebra in $\Fun(I,\cS)_{/c^*\Pic}$ by letting the $\bE_k$-algebra structure in the $\CAlg(\Sp)$-component be the restricted $\bE_\infty$-algebra structure. More generally, an $\bE_k$-algebra over an $\bE_\infty$-algebra in $\Fun(I,\cS_{/\Pic_R})$ gives rise to an $\bE_k$-algebra over an $\bE_\infty$-algebra in $\Fun(I,\cS)_{/c^*\Pic}$. 

\begin{proposition}\label{prop:Bcy-THH-comparison-over-PicR}
  Let $\xi_*\: M_* \to c^*\Pic_R$ be an $\bE_k$-algebra in $\Fun(I,\cS_{/\Pic_R})$, viewed as an  $\bE_k$-algebra in $\Fun(I,\cS)_{/c^*\Pic}$ as explained above, and let $\xi\: M \to \Pic_R$ be its total object. 
  \begin{enumerate}[(i)]
\item There is a natural equivalence of $I$-graded $\bE_{k-1}$ $\THH(R)$-algebras with $\bT$-action 
\begin{equation}\label{eq:graded-THH-Bcy-over-R}
\THH(\Th_R(\xi_*)) \simeq  \Th_{\THH(R)}(\xi^\cy_* \: B^\cy(M_*) \to c^*\Pic_{\THH(R)})
\end{equation}
and an analogous equivalence of total $\bE_{k-1}$ $\THH(R)$-algebras with $\bT$-action
\begin{equation}\label{eq:ungraded-THH-Bcy-over-R}
\THH(\Th_R(\xi)) \simeq  \Th_{\THH(R)}(\xi^\cy \: B^\cy(M) \to \Pic_{\THH(R)})
\,.
\end{equation}
\item There is a natural equivalence of $\bE_{k-1}$-algebras in  $(\cS_{/ \Pic_{\THH(R)}})^{B\bT}$ between
  \[ c_!(\xi^\cy_*) \qqandqq \xi^\cy\: B^\cy(M) \to \Pic_{\THH(R)}. \]
\item Passing to total objects, the equivalence~\eqref{eq:graded-THH-Bcy-over-R} recovers~\eqref{eq:ungraded-THH-Bcy-over-R}.
\end{enumerate}
\end{proposition}
\begin{proof}
  Let $\iota_R \: \Delta^0 \to \Pic_R$ be the unit of $\cS_{/\Pic_R}$. Then $\iota_R$ can be viewed as an $\bE_{\infty}$-algebra in $\cS_{/\Pic_R}$ (resp. in $\Fun(I,\cS_{/\Pic_R})$) and $\xi$ (resp. $\xi_*$) can be viewed as an $\bE_k$ $\iota_R$-algebra. Viewing $\xi$ and $\xi_*$ as $\iota_R$-algebras in $\cS_{/\Pic_R}$ and $\Fun(I,\cS)_{/c^*\Pic}$ in the way explained before the proposition, the claims now follow from applying Corollary~\ref{cor:Bcy-on-E_k-alg-relative}(iii) to suitable parts of the commutative square 
\[\xymatrix@-1pc{\Fun(I,\cS)_{/ c^*\Pic} \ar[rr]^{u \Th} \ar[d]_{c_!} && \Fun(I,\Sp)\ar[d]^{c_!} \\ 
    \cS_{/ \Pic} \ar[rr]^{u \Th}  && \Sp }
\]
  of symmetric monoidal functors that preserve sifted colimits.
\end{proof}

\begin{remark}
If $R = \bS$, then the equivalences in Proposition~\ref{prop:Bcy-THH-comparison-over-PicS} agree with the equivalences in Proposition~\ref{prop:Bcy-THH-comparison-over-PicR}, keeping in mind that $\THH(\bS) = \bS$.  This boils down to the observation that as a symmetric monoidal functor, $\Th_\bS \: \cS_{/\Pic_\bS} \to \Sp$ factors as the composite of $\cS_{/\Pic_\bS} \to \cS_{/\Pic}$, mapping $\xi$ to $(\bS, \xi)$, followed by~$u\Th$.
\end{remark}

\subsection{Log  \texorpdfstring{$\THH$}{THH} of Thom algebras} As in Proposition~\ref{prop:Bcy-THH-comparison-over-PicR}, we consider an $\bE_k$-map $\xi_*\: M_* \to c^*\Pic_R$. By the discussion in Construction~\ref{cons:xigp}, the group completion $M^\gp$ of $M$ gives rise to a map $\xi_*^\gp \: M_*^\gp \to c^*\Pic_R$.
\begin{definition}\label{def:BIrep}
  The \emph{$I$-replete bar construction} of $\xi_*$ is the $\bT$-equivariant $\bE_{k-1}$-map
  \[\xi^\rep_*\: B^{\Irep}(M_*)\to c^*\Pic_{\THH(R)}\]
  obtained by first applying $B^{\cy}$ to $(R,\xi^{\gp}_*\: M_*^{\gp} \to c^*\Pic_R)$ and then applying the restriction $\Alg_{\bE_{k-1}}((\Fun(I^\gp,\cS)_{/c^*\Pic})^{B\bT})\to \Alg_{\bE_{k-1}}((\Fun(I,\cS)_{/c^*\Pic})^{B\bT})$ along $\gamma \: I \to I^\gp$. Passing to total objects defines the ungraded map $\xi^\rep\: B^{\Irep}(M) \to \Pic_{\THH(R)}$. 

Using Corollary~\ref{cor:Bcy-on-E_k-alg}(ii), the group completion map $M_* \to \gamma^*(M_*^\gp)$ induces a natural \emph{$I$-graded repletion map} $\rho_*\: B^\cy(M_*) \to B^\Irep(M_*)$, which is a $\bT$-equivariant $\bE_{k-1}$-map over $c^*\Pic_{\THH(R)}$. Passage to total objects defines the (ungraded) \emph{repletion map} $\rho\: B^\cy(M) \to B^\Irep(M)$. 
\end{definition}
Encoding the gradings of $M_*$ and $M_*^\gp$ by $\bE_k$-maps $\tau\: M \to I$ and $\tau^\gp\: M^\gp \to I^{\gp}$ as in~\eqref{eq:xigp}, the construction of $B^{\Irep}$ can be rephrased by requiring the lower left-hand and the lower middle square in the following diagram to be pullbacks:
\begin{equation} \label{eq:repbar}
\xymatrix@C+.5pc{
I \ar[d]_-{=}
	& B^\cy(I) \ar[d]^-{\rho} \ar[l]_-{\epsilon}
	& B^\cy(M) \ar[d]^-{\rho} \ar[l]_-{B^\cy(\tau)} \ar[r]^-{B^\cy(\xi)}
	& B^\cy(\Pic_R) \ar[dd]^-{B^\cy(\gamma)}_-{\simeq}  \ar[r]^-{\hat\epsilon} & \Pic_{\THH(R)} \\
I \ar[d]_\gamma
	& B^\rep(I) \ar[d]^-{\tilde\gamma} \ar[l]_-{\tilde\epsilon}
	& B^\Irep(M) \ar[d]^-{\tilde\gamma} \ar[l]_-{B^\rep(\tau)} \ar@{-->}[ur]_-{B^\rep(\xi)} \\
I^\gp
	& B^\cy(I^\gp) \ar[l]_-{\epsilon}
	& B^\cy(M^\gp) \ar[l]_-{B^\cy(\tau^\gp)} \ar[r]^-{B^\cy(\xi^\gp)}
	& B^\cy(\Pic_R^\gp)\rlap{\,.}
}
\end{equation}
Here the two maps labeled $\epsilon$ result from $I$ being commutative and $\hat\epsilon$ arises from applying Corollary~\ref{cor:Bcy-on-E_k-alg}(ii) to the lax symmetric monoidal structures of~$\Pic$. The augmentation map $\xi^\rep\: B^{\Irep}(M)\to \Pic_{\THH(R)}$ from Definition~\ref{def:BIrep} is the composite of $B^\rep(\xi)$ and $\hat\epsilon$. This description of $B^{\Irep}(M)$ highlights its (implicit) dependence on the choice of the $I$-grading $\tau\: M \to I$. In contrast, $B^\cy(M)$ is also defined for an ungraded $\bE_k$-space.

\begin{proposition}\label{prop:THHThRxistarxi}
In the situation of Proposition~\ref{prop:Bcy-THH-comparison-over-PicR}, there
is a natural equivalence of $I$-graded $\bE_{k-1}$ $\THH(R)$-algebras with $\bT$-action 
\[
\THH(\Th_R(\xi_*), \xi) \simeq \Th_{\THH(R)}(\xi^\rep_*\: B^\Irep(M_*) \to c^*\Pic_{\THH(R)}) \,.
\]
Passing to total objects, this provides an equivalence
\[
\THH(\Th_R(\xi), \xi) \simeq  \Th_{\THH(R)}(\xi^\rep
	\: B^\Irep(M) \to \Pic_{\THH(R)})
\]
of $\bE_{k-1}$ $\THH(R)$-algebras with $\bT$-action.

These equivalences are compatible with those
of Proposition~\ref{prop:Bcy-THH-comparison-over-PicR} via the
repletion maps $\rho \: \THH(\Th_R(\xi)) \to \THH(\Th_R(\xi), \xi)$, $\rho
\: B^\cy(M) \to B^\Irep(M)$ and their $I$-graded variants.
\end{proposition}
\begin{proof}
Arguing as in the proof of Proposition~\ref{prop:Bcy-THH-comparison-over-PicR} with $\xi^{\gp}_*\: M_*^{\gp} \to c^*\Pic_R$ (in place of $\xi_*$) and forming the restriction gives the first equivalence, and since $u\Th_{\THH(R)}$ commutes with coproducts passing to total objects gives the second. For the compatibility with the equivalences from  Proposition~\ref{prop:Bcy-THH-comparison-over-PicR}, we first notice that the two lax symmetric monoidal functors $\gamma^*\circ u\Th$ and $u\Th \circ \gamma^*$ from $\Fun(I^\gp,\cS)_{/c^*\Pic}$ to $\Fun(I,\Sp)$ are equivalent. The claim then follows from Corollary~\ref{cor:Bcy-on-E_k-alg-relative}(ii), making use of the compatibility of the comparison map with composites of lax symmetric monoidal functors. 
\end{proof}
\subsection{Comparison of replete bar constructions}\label{subsec:comparison-rep-bar} For the rest of this section, we will allow $I$ to be any $\bE_{\infty}$-space. If $M \to I$ is an $\bE_{\infty}$-map, we define $B^{\Irep}(M)$ to be the pullback $I\times_{I^\gp}B^\cy(M^\gp)$ as in~\eqref{eq:repbar}. For $M \to I$ the identity, we denote the resulting replete bar construction by $B^\rep(M)$. It coincides with the version in~\cite{Rog09}*{Cor.~8.10}.

Recall from Definition~\ref{def:zero-replete} that $M$ is $\pi_0$-replete if $M \to \pi_0(M)\times_{\pi_0(M^\gp)}M^\gp$ is an equivalence.
\begin{lemma}
  Let $M$ be a $\pi_0$-replete $\bE_{\infty}$-space and let $M \to I = \pi_0(M)$ be the canonical $\pi_0(M)$-grading from Definition~\ref{def:canonical-pi0grading}. Then there is a natural equivalence $B^{\Irep}(M) \to B^\rep(M)$ of $I$-graded $\bE_{\infty}$-spaces. 
\end{lemma}
\begin{proof}
Writing $B^{\Irep}(M)$ as an iterated pullback, this is immediate. 
\end{proof}
\begin{corollary} When $R = \bS$ and $M$ is a $\pi_0$-replete $\bE_\infty$-space, then the notion of log $\THH$ from Definition~\ref{def:logTHHAxi} is equivalent to that of~\cite{Rog09}*{Def.~8.11}. 
\end{corollary}
\begin{proof}
This follows from the last corollary and Proposition~\ref{prop:THHThRxistarxi}. 
\end{proof}  
An analogous compatibility result with the notion of log $\THH$ from~\cite{RSS15}*{Def.~4.6}, \cite{RSS18}*{Def.~4.1} is discussed in Appendix~\ref{app:pointsetlevel-logTHH-comparision}. The $\pi_0$-repleteness hypothesis is satisfied in all the cases for which specific calculations were made in those papers.

\begin{remark}
  If $M$ is an $\bE_k$-space, then~\cite{L:HA}*{Cor.~5.1.1.7} implies that the Postnikov truncation $\tau_{\leq n}M$ has a preferred $\bE_{\infty}$-structure if $ n+1<k$. In this situation, the truncation $\tau \: M \to \tau_{\leq n}M$ provides a $\tau_{\leq n}M$-grading on $M$ that gives rise to a $\tau_{\leq n}M$-replete bar construction. For $n=\infty$, this gives the $M$-replete bar construction. For $n=0$ and $k \geq 2$, we get the $I$-replete bar construction for the canonical $I=\pi_0(M)$-grading. We will not study the versions for $0 < n < \infty$ here.  
\end{remark}

\section{Logarithmic structures and logification}\label{sec:logifications}
The purpose of this section is to generalize the log condition and the logification process for discrete prelog rings to prelog $\bE_k$ $R$-algebras. Based on this, we compare different notions of log $\THH$ (see Appendix~\ref{app:pointsetlevel-logTHH-comparision}) and use the results of Section~\ref{sec:THH-of-Thom-rings} to show that in some cases of interest, log $\THH$ is independent of the choices made in the construction of $\bE_2$-prelog structures. 

\subsection{Log \texorpdfstring{$\bE_k$ R}{Ek R}-algebras} 
If $(A, \xi, \bar\alpha)$ is a prelog $\bE_k$ $R$-algebra, then the adjunction~\eqref{eq:R-mod-E-k-adjunction} allows us to equivalently describe the structure map $\bar\alpha$ by its right adjoint $\alpha\: M \to \Omega^R(A)$ in  $\Alg_{\bE_k}(\cS_{/ \Pic_{R}})$. Since $u\: \Mod_R \to \Sp$ is in general only lax symmetric monoidal, $u\: \Alg_{\bE_k}(\Mod_R) \to \Alg_{\bE_k}(\Sp)$ has no right adjoint, and so there is no adjoint structure map for an $R$-based prelog $\bE_k$-ring.

\begin{example}\label{ex:directimageprelog}
  Let $(B, \eta\: N \to \Pic_R, \bar\beta\: \Th_R(\eta) \to B)$ be a prelog $\bE_k$ $R$-algebra and let $f\: A \to B$ be an $\bE_k$ $R$-algebra map. The \emph{direct image} prelog $\bE_k$ $R$-algebra $(A,f_*\eta\: f_*N \to \Pic_R, f_*\bar\beta\: \Th_R(f_*\eta) \to A)$ has as structure map $f_*\bar\beta$ the left adjoint of the map $f_*\beta$ defined by the pullback 
  \[\xymatrix@-1pc{
      f_*N \ar[d]^{\tilde f} \ar[rr]^{f_*\beta} && \Omega^R(A) \ar[d]^{\Omega^R(f)} \\ 
      N \ar[rr]^{\beta} && \Omega^R(B)
    }\]
in $\Alg_{\bE_k}(\cS_{/ \Pic_{R}})$, while $f_*\eta = \eta \circ \tilde f $. 
\end{example}
\begin{definition}
  Let $A$ be an $\bE_k$ $R$-algebra. Its \emph{space of units} $\GL_1^R(A)$ consists of the connected components of $\Omega^R(A)$ corresponding to units in $\pi_0(\Omega^R(A))$, and is itself an $\bE_k$-space by~\cite{ABG18}*{Lem.~7.2}. The composite $\xi^{\GL_1} \: \GL_1^R(A) \to \Pic_R$ of the inclusion~$\iota\: \GL_1^R(A) \to \Omega^R(A)$ and the structure map of $\Omega^R(A)$ provides an object of $\Alg_{\bE_k}(\cS_{/ \Pic_{R}})$ sometimes also denoted by $\GL_1^R(A)$.
\end{definition}
Proposition~\ref{prop:Omega-S-multiplicative-monoid} implies that $\GL_1^\bS(A)$ (as a space over $\Pic_\bS$) captures the units of the graded ring $\pi_*(A)$. 

\begin{definition}
  A \emph{prelog $\bE_k$ $R$-algebra} $(A, \xi, \bar\alpha)$ is a \emph{log $\bE_k$ $R$-algebra} if the map $\tilde\alpha$ in the following pullback square in $\Alg_{\bE_k}(\cS_{/ \Pic_{R}})$ is an equivalence:
  \begin{equation}\label{eq:pullback-log-condition}\xymatrix@-1pc{
      \alpha^{-1}(\GL_1^R(A)) \ar[rr]^-{\tilde\alpha}\ar[d] && \GL_1^R(A) \ar[d]^\iota \\ 
      M \ar[rr]^-{\alpha} && \Omega^R(A)\rlap{\,.}
    }
  \end{equation}
  A prelog $\bE_k$-ring is a \emph{log $\bE_k$-ring} if it is a log $\bE_k$ $\bS$-algebra.
\end{definition}
\begin{example}
  Let $A$ be an $\bE_k$ $R$-algebra and let $\bar\iota$ be the adjoint of the inclusion~$\iota$. Then  $(A, \xi^{\GL_1}, \bar\iota)$ is the \emph{trivial log $\bE_k$ $R$-algebra}. 
\end{example}

\begin{example}\label{ex:directimagelog}
If $f\: B \to A$ is an $\bE_k$ $R$-algebra map and $(A, \xi, \bar\alpha)$ is a log $\bE_k$ $R$-algebra, then the direct image prelog $\bE_k$ $R$-algebra $(B,f_*\xi, f_*\bar\alpha)$ of Example~\ref{ex:directimageprelog} is in fact a log $\bE_k$ $R$-algebra. (This can be seen by identifying the pullback of $f_*\alpha$ along $\GL_1^R(B) \to \Omega^R(B)$ with the pullback of $\alpha$ along $\GL_1^R(B) \to \GL_1^R(A) \to \Omega^R(A)$.) In particular, any morphism of $\bE_k$ $R$-algebras $f\: B \to A$ gives rise to a log $\bE_k$ $R$-algebra $(B,f_* \xi^{\GL_1}, f_*\bar\iota)$, the direct image log structure associated with the trivial log structure on $A$.  
\end{example}

\subsection{Logification} We can build log structures from prelog structures:
\begin{construction}\label{cons:logification}
  Let $(A, \xi, \bar\alpha)$ be a prelog $\bE_k$ $R$-algebra with $k\geq 2$. We view the maps $M \leftarrow \alpha^{-1}(\GL_1^R(A)) \to \GL_1^R(A)$ in $\Alg_{\bE_k}(\cS_{/ \Pic_{R}})$ from the pullback square~\eqref{eq:pullback-log-condition} as an $\bE_{k-1}$-algebra in the span category of~$\Alg_{\bE_1}(\cS/_{\Pic_R}))$ with the pointwise symmetric monoidal structure~\cite{L:HA}*{Rem.~2.1.3.4}. Viewing $M$ as a left $\alpha^{-1}(\GL_1^R(A))$-module and $\GL_1^R(A)$ as a right $\alpha^{-1}(\GL_1^R(A))$, we can form the two-sided bar construction $M^a := B(M,\alpha^{-1}(\GL_1^R(A)),\GL_1^R(A))$ in the symmetric monoidal category $\cS/_{\Pic_R}$. Since the two-sided bar construction is symmetric monoidal with respect to the pointwise symmetric monoidal structure, $M^a$ inherits an $\bE_{k-1}$-algebra structure. The unit of $\alpha^{-1}(\GL_1^R(A))$ induces an $\bE_{k-1}$-map of spans of algebras, yielding  $\bE_{k-1}$-maps $M \to M^a$ and $\GL_1^R(A) \to M^a$. The $\bE_{k-1}$-map to $\Omega^R(A) \leftarrow \Omega^R(A) \to \Omega^R(A)$ in the span category induces the structure map $\alpha^a\: M^a \to \Omega^R(A)$. This factorization induces the \emph{logification} map
  \begin{equation}\label{eq:logification}
 (A,\xi, \bar\alpha)  \longto (A,\xi^a\: M^a \to \Pic_R, \bar\alpha^a\: \Th_R(\xi^a) \to A) =: (A,\xi^a, \bar\alpha^a) 
\end{equation}
of prelog $\bE_{k-1}$ $R$-algebras whose codomain is the \emph{associated log $\bE_{k-1}$ $R$-algebra} of $(A, \xi, \bar\alpha)$. 
\end{construction}

\begin{lemma}
The prelog $\bE_{k-1}$ $R$-algebra $(A,\xi^a, \bar\alpha^a)$ resulting from Construction~\ref{cons:logification} is a log $\bE_{k-1}$ $R$-algebra. Any map of prelog $\bE_k$ $R$-algebras $(A,\xi, \bar\alpha) \to (B,\eta,\bar\beta)$ to a log $\bE_k$ $R$-algebra extends canonically as an $\bE_{k-1}$ map over the logification map~\eqref{eq:logification}. 
\end{lemma}
\begin{proof}
  We first show that $(A,\xi^a, \bar\alpha^a)$ is log and argue as in the case of discrete prelog rings or $\cJ$-space prelog ring spectra (see e.g.~\cite{S14}*{Lem.~3.12}). For a map $W \to \Omega^R(A)$, the canonical $\pi_0(\Omega^R(A))$-grading of $\Omega^R(A)$ induces a $\pi_0(\Omega^R(A))$-grading on $W$, which induces a coproduct decomposition $\tilde W\amalg \hat W$ in $\alpha^{-1}(\GL_1^R(A))$-modules in $\cS_{/\Pic_R}$ of the parts of $W$ over the units $\pi_0(\Omega^R(A))^{\times}$ and its complement. For $\alpha \: M \to \Omega^R(A)$ we have $\tilde M = \alpha^{-1}(\GL_1^R(A))$, and there are equivalences
\[ M^a = M\otimes_{\tilde M} \GL_1^R(A) \simeq (\tilde M\textstyle\amalg \hat M)\otimes_{\tilde M} \GL_1^R(A) \simeq \GL_1^R(A) \amalg (\hat M\otimes_{\tilde M} \GL_1^R(A)) \] 
where the latter coproduct is the decomposition for $M^a \to \Omega^R(A)$. This shows that $M^a \to \Omega^R(A)$ is log. 
The second part follows by naturality and the observation that $N \to N^a$ is an equivalence if $(B,\eta\colon N \to \Pic_R,\bar\beta \colon \Th_R(\eta) \to B)$ is log. 
\end{proof}
When $k=2$, then $(A,\xi^a, \bar\alpha^a)$ is only $\bE_1$ so that its log $\THH$ is not defined. To remedy this, we use the following tweak.
\begin{definition}\label{def:E_k-logification}
A map of prelog $\bE_k$ $R$-algebras $(A,\xi, \bar\alpha) \to (B,\eta,\bar\beta)$ is an \emph{$\bE_k$-logification} if $(B,\eta,\bar\beta)$ is log and the map induces an equivalence of prelog $\bE_{k-1}$ $R$-algebras $(A,\xi^a, \bar\alpha^a) \to (B,\eta^a,\bar\beta^a)$. 
\end{definition}
Since logification does not change the underlying $\bE_k$ $R$-algebra, the map $A \to B$ in the definition is an equivalence. 

\begin{construction}\label{constr:map-to-direct-image} Suppose that $k = 2$, that $d > 0$, that $f\: A \to B$ is a map of $\bE_k$-rings that exhibits $A$ as a connective cover of $B$, and that $a \in \pi_{2d}(A)$ is a homotopy class such that $f$ induces an isomorphism $\pi_*(A)[1/a] \to \pi_*(B)$. (In other words, $A$ is the connective cover of the $2d$-periodic $\bE_k$-ring $B$.) Suppose further that $\pi_{i(2d+2)+2d-1}(A) \cong 0$ for all $i \geq 1$. Then by Remark~\ref{rem:vanishing-for-one-gen-case}, there exists an $\bE_2$-map $\bar\alpha \: \Th_{\bS}(\xi_{2d}) = \bS[x] \to A$ sending the generator $x$ of degree $2d$ to $a$. This defines a prelog $\bE_2$-ring $(A,\xi_{2d},\bar\alpha)$. The adjoint of the composite $\bS[x] \to A \to B$ factors through the inclusion $\GL_1^\bS(B) \to \Omega^\bS(B)$ and therefore induces a map $(A,\xi_{2d},\bar\alpha) \to (A,f_* \xi^{\GL_1}, f_*\bar\iota)$ to the direct image log $\bE_k$-ring determined by~$f$ (compare Example~\ref{ex:directimagelog}).
\end{construction}

\begin{lemma}\label{lem:map-to-direct-image}
  In this situation, $(A,\xi_{2d},\bar\alpha) \to (A,f_* \xi^{\GL_1}, f_*\bar\iota)$ is an $\bE_2$-logification of $(A,\xi_{2d},\bar\alpha)$, and any $\bE_2$-logification of $(A,\xi_{2d},\bar\alpha)$ is equivalent as a prelog $\bE_2$-ring to $(A,f_* \xi^{\GL_1}, f_*\bar\iota)$. 
\end{lemma}
\begin{proof}
  The first part is analogous to~\cite{S14}*{Lem.~4.7}. For $\alpha \: \bZ_{\geq 0} \to \Omega^\bS(A)$ we have $\alpha^{-1}(\GL_1^\bS(A)) \simeq \{0\}$ since $\pi_*(A)$ has no units of positive degree. Hence $(\bZ_{\geq 0})^a \simeq \bZ_{\geq 0} \times \GL_1^\bS(A)$ and it remains to show that $\bZ_{\geq 0} \times \GL_1^\bS(A) \to f_*\GL_1^\bS(B)$ is an equivalence. For this it is convenient to use the $\bZ$-grading of $\Omega^\bS(B)$ induced by the canonical $\pi_*(B)/\{\pm 1\}$ grading and the map $\pi_*(B)/\{\pm 1\} \to \bZ$ sending a homotopy class to its degree. Proposition~\ref{prop:Omega-S-multiplicative-monoid} implies that after restricting along $(\bZ_{\geq 0} \to \bZ)$, the map $\Omega^{\bS}(A) \to \Omega^{\bS}(B)$ becomes an equivalence. Hence the map $\GL_1^\bS(A) \to (\{0\}\to \bZ)^*(\GL_1^\bS(B))$ is an equivalence, and so the induced map
  \[ \bZ_{\geq 0} \times \GL_1^\bS(A) \to f_*(\GL_1^\bS(B)) \simeq (\bZ_{\geq 0} \to \bZ)^* (\GL_1^\bS(B))\]
is an equivalence because $\GL_1^\bS(B)$ is a grouplike $\bE_2$-space and all non-empty components in the $\bZ$-grading are therefore equivalent. 

  For the second statement, let $(A,\xi_{2d},\bar\alpha) \to (A,\eta\: N \to \Pic_\bS,\bar\beta)$ be any $\bE_2$-logification. Then the image of $\pi_0(N)$ in $\pi_0(\Omega^\bS(B))$ is contained in the units $\pi_0(\Omega^\bS(B))^\times = (\pi_*(B))^\times /\{\pm 1\}$ by the first part. Hence there is an induced map $ (A,\eta,\bar\beta) \to (A,f_* \xi^{\GL_1}, f_*\bar\iota)$, which is an equivalence since both are logifications of $(A,\xi_{2d},\bar\alpha)$. 
\end{proof}
\begin{example}\label{ex:logification-of-E2-prelog-on-ell-ku}
The previous lemma applies in particular to the prelog $\bE_2$-rings $(ku,\< u \>)$, $(\kup,\< u \>)$, and $(\ell,\< v_1\>)$ of Definition~\ref{def:ku-ell-BPn} and shows that they admit $\bE_2$-logifications that do not depend on the choices made in the construction of the $\bE_2$-prelog structures. Moreover, these logifications are log $\bE_\infty$-rings since they arise as direct image prelog $\bE_\infty$-rings from the $\bE_\infty$-maps $A \to A[1/a]$ and the trivial log $\bE_\infty$-rings $(A[1/a], \xi^{\GL_1}, \bar\iota)$ in the way explained in Example~\ref{ex:directimagelog}. This argument does not apply to the other examples introduced in Definition~\ref{def:ku-ell-BPn}, since they all involve degree $0$ elements $a$ for which $A$ is not the connective cover of $A[1/a]$. We do not know how to construct the logifications as direct image log structures in these cases. 
\end{example}

\subsection{Logification invariance of log \texorpdfstring{$\THH$}{THH}} The next result uses the induced map of Example~\ref{ex:maps-on-log-THH-w-can-grading}. It generalizes~\cite{RSS15}*{Thm.~4.24}.  
\begin{theorem}\label{thm:logification-invariance}
  Let $(A,\xi,\bar\alpha)$ be a prelog $\bE_{k}$ $R$-algebra with an $\bE_k$-logification $(f,f^\flat) \: (A,\xi, \bar\alpha) \to (B,\eta,\bar\beta)$. Then $(f,f^\flat)$ induces an equivalence $\THH(A,\xi,\bar\alpha) \to \THH(B,\eta,\bar\beta)$ on log $\THH$ formed with respect to the canonical $\pi_0$-grading. 
\end{theorem}
\begin{proof}
  The argument for~\cite{RSS15}*{Thm.~4.24} applies also in the present context and we only explain the necessary modifications. Firstly, using Propositions~\ref{prop:Bcy-THH-comparison-over-PicR} and~\ref{prop:THHThRxistarxi} we can argue with $B^\cy$ and $B^{\zerorep}$, the replete bar construction for the $\pi_0$-repletion. As in \emph{loc.\,cit.}, we write $G = \GL_1^\bS$ and $P = \alpha^{-1}(G)$, while $N$ replaces the $M^a$ there. We then get a commutative square
  \begin{equation}\label{eq:P-G-M-N-square}
\xymatrix@-1pc{P \ar[rr] \ar[d] && G \ar[d] \\ M \ar[rr] && N } 
\end{equation}
in $\Alg_{\bE_k}(\cS/_{\Pic_R})$ inducing an equivalence of $\bE_{k-1}$-algebras $M\otimes_P G \to N$. 

As in \emph{loc.\,cit.}, the canonical map $B^\cy(P^\gp) \to B^\zerorep(P)\otimes_{B^\cy(P)}B^\cy(P^\gp)$ is an equivalence, its domain and codomain are grouplike, and it induces an equivalence after group completion (as one can check by applying $B$ to the counterpart of the first square in the proof of \emph{loc.\,cit.}). Replacing the $\bS^\cJ$ in the reference by $\Th_R$, the same argument implies that for the repletion map $\rho \: B^\cy(P) \to  B^\zerorep(P)$, the cobase change of $\Th_R(\rho)$ along $\Th_R(B^\cy(P)\to \Pic_R) \to \THH(A)$ is an equivalence. 

Since $B^\cy$ is symmetric monoidal, \eqref{eq:P-G-M-N-square} induces an equivalence of $\bE_{k-2}$-algebras $B^\cy(M)\otimes_{B^\cy(P)}B^\cy(G) \to B^\cy(N)$. Similarly, it induces an equivalence of discrete commutative monoids $\pi_0(M) \times_{\pi_0(P)}\pi_0(G) \to \pi_0(N)$. Arguing as in the proof of~\cite{RSS15}*{Lem.~4.26}, this implies that also  $B^\zerorep(M)\otimes_{B^\zerorep(P)}B^\zerorep(G) \to B^\zerorep(N)$ is an equivalence. Applying $\Th_R$, we get the counterpart of the cubical diagram in the proof of \emph{loc.\,cit.}, and the conclusion follows analogously. 
\end{proof}

\begin{corollary}The log $\THH$ of the prelog $\bE_2$-rings discussed in Example~\ref{ex:logification-of-E2-prelog-on-ell-ku} is independent of the choices made in the construction of the $\bE_2$-prelog structures. \qed
\end{corollary}
Theorem~\ref{thm:logification-invariance} will also be used in Appendix~\ref{app:pointsetlevel-logTHH-comparision} to compare the present notion of log $\THH$ to the one used in \cite{RSS15}*{Def.~4.6} and \cite{RSS18}*{Def.~4.1}.

\section{The repletion--residue sequence in the one generator case}\label{sec:repletion-residue-one-generator}
In this section, we specialize to $\bE_k$ prelog structures $\xi\: M \to \Pic_R$ with $k\geq 2$ where $M = \<x\> = \{x^i \mid i\ge0\}$ is the free commutative monoid on one generator~$x$,  graded by $I = \bZ_{\ge0}$ via the $\bE_\infty$-isomorphism $\tau \: M \to I, x^i \mapsto i$. As in Construction~\ref{cons:xigp}, we write $\xi^\gp\: M^\gp = \<x^{\pm1}\> \to \Pic_R$ for the preferred $\bE_k$-extension of $\xi$ along the group completion map of $M$.

Our goal is to identify the fiber of the repletion map $\rho \: \THH(A) \to \THH(A, \xi, \bar\alpha)$ from~\eqref{eq:repletion-relative-to-THHA} for an $R$-based prelog $\bE_k$-ring  $(A, \xi \: \<x\> \to \Pic_R, \bar\alpha \: \Th_R(\xi) \to A)$. To do so, we first consider the weight~$0$ component
\begin{equation} 
\rho_0 \: \THH(\Th_R(\xi_*); 0) \longto \THH(\Th_R(\xi_*), \xi; 0) = \THH(\Th_R(\xi^\gp_*); 0)
\end{equation}
of the repletion map for the canonical prelog structure, introduced in~\eqref{eq:rho_i_for_canonical}. As discussed in Subsection~\ref{subsec:logTHH-of-canonical}, this is a map of cyclotomic $\THH(R)$-modules. 
\begin{definition}
The \emph{repletion fiber} $\THH^\xi(R)$ is the fiber of $\rho_0$ in cyclotomic $\THH(R)$-modules. We say that the pair $(R,\xi_*)$ is \emph{cyclotomically good} if $\THH(R)$ and $\THH^\xi(R)$ are equivalent as cyclotomic $\THH(R)$-modules.
\end{definition}
We will see in Lemma~\ref{lem:repletion-fiber-as-THHR-module} below that the underlying $\THH(R)$-modules of $\THH(R)$ and $\THH^\xi(R)$ are equivalent without further assumptions on~$R$, and refine this to a $\bT$-equivariant statement in Proposition~\ref{prop:THHThRxixi0-Tact}. However, we only know how to control the full cyclotomic structure of $\THH^\xi(R)$ under additional assumptions: 
\begin{proposition}\label{prop:sufficient-criteria-cycl-good}
If $R = \bS$, $R = MU$ or $R = H\bF$ with $\bF$ a perfect field of characteristic $p$ or a perfectoid ring, then $(R,\xi_*)$ is cyclotomically good. 
\end{proposition}
We prove this proposition in Subsection~\ref{subsec:p-cyc-good-criteria}, but first explain its relevance for understanding the fiber of the repletion map of an $R$-based prelog $\bE_k$-ring $(A, \xi \: \<x\> \to \Pic_R, \bar\alpha \: \Th_R(\xi) \to A)$. The $\bE_k$ $R$-algebra collapse map $\coll \: \Th_R(\xi) \to R$ obtained from the face inclusion $\{1\} = \{x^0\} \subseteq \<x\>$ as an instance of~\eqref{eq:col-on-ThR} allows us to form the $\bE_{k-1}$ $A$-algebra
\begin{equation}\label{eq:Amodbaralphax}
  \modmod{A}{\bar\alpha} :=  A \otimes_{\Th_R(\xi)} R
\end{equation}
as in Definition~\ref{def:AMIJ}. Restricting the $\THH(R)$-module structure on $\THH^\xi(R)$ along $\THH(\coll) \: \THH(\Th_R(\xi)) \to \THH(R)$, we can form the balanced smash product
\[
\THH^\xi(\modmod{A}{\bar\alpha}) := \THH(A) \otimes_{\THH(\Th_R(\xi))} \THH^\xi(R) 
\]
in  cyclotomic $\THH(A)$-modules. 

The next theorem, proved in Subsection~\ref{subsec:constr-repletion-residue} below, generalizes~\cite{RSS15}*{Thm.~1.1}, where we only considered certain $\bE_{\infty}$ prelog structures over $R=\bS$ (compare Appendix~\ref{app:pointsetlevel-logTHH-comparision}) and did not record any cyclotomic structure. Together with Proposition~\ref{prop:sufficient-criteria-cycl-good} and Corollary~\ref{cor:rep-res-TC}, it implies Theorem~\ref{thm:repletion-residue-intro} from the introduction. 
\begin{theorem} \label{thm:rep-res-Rbas}
Let $(A, \xi \: \<x\> \to \Pic_R, \bar\alpha
\: \Th_R(\xi) \to A)$ be an $R$-based prelog $\bE_k$-ring with $k\geq 2$.  There is a natural repletion--residue cofiber sequence
\[
\THH^\xi(\modmod{A}{\bar\alpha})
	\xrightarrow{\trf} \THH(A)
	\xrightarrow{\rho} \THH(A, \xi, \bar\alpha)
	\xrightarrow{\res} \Sigma \THH^\xi(\modmod{A}{\bar\alpha})
\]
of cyclotomic $\THH(A)$-modules. If the pair $(R,\xi_*)$ is cyclotomically good, then $\THH^\xi(\modmod{A}{\bar\alpha})$ is equivalent to $\THH(\modmod{A}{\bar\alpha})$ as a cyclotomic $\THH(A)$-module.
\end{theorem}
Informally,  when $A=\Th_R(\xi)$, the map $\res \: \THH(\Th_R(\xi), \xi) \to \Sigma \THH^\xi(R) \simeq
\Sigma \THH(R)$ is given by first collapsing all weight $i>0$ summands to the zero spectrum~$0$, then viewing the quotient as $\Th_{\THH(R)}(\xi^\rep_0)$ formed over $B^\rep(\<x\>; 0) \simeq S^1(0)$, then collapsing the copy of $\THH(R)$ over the $0$--simplex $(1)$ to $0$, and finally identifying this quotient as a suspension of $\THH(R)$.
We write $\trf$ to suggest transfer, but have not verified compatibility of this map under the (cyclotomic) trace map with the algebraic $K$-theory transfer map $K(R) \to K(\Th_R(\xi))$. Partial results in this direction can be derived from~\cite{Lun25}. \begin{remark} When $(R,\xi_*)$ is not known to be cyclotomically good, then Proposition~\ref{prop:THHThRxixi0-Tact} still implies that $\THH^\xi(\modmod{A}{\bar\alpha})$ is equivalent to $\THH(\modmod{A}{\bar\alpha})$ as a $\THH(A)$-module with $\bT$-action. 
\end{remark}

\begin{corollary} \label{cor:rep-res-TC}
  Let $(A, \xi\: \<x\> \to \Pic_R, \bar\alpha)$ be an $R$-based prelog $\bE_k$-ring with $k\geq 2$ such that $(R,\xi_*)$ is cyclotomically good. Then there is a repletion--residue cofiber sequence
\[
\TC(\modmod{A}{\bar\alpha})_p
	\xrightarrow{\trf} \TC(A)_p
	\xrightarrow{\rho} \TC(A, \xi, \bar\alpha)_p
	\xrightarrow{\res} \Sigma \TC(\modmod{A}{\bar\alpha})_p
\]
of $\TC(A)_p$-modules. 
\end{corollary}
\subsection{Construction of the repletion--residue sequence}\label{subsec:constr-repletion-residue} Next we recall the cyclic and replete bar constructions for the free commutative monoid on one generator. Here and later, weak equivalences of spaces with $\bT$-action or $\bT$-equivariant weak equivalences refer to $\bT$-equivariant maps that are underlying weak equivalences, also known as naive equivalences.  
\begin{proposition}[\cite{Hes96}*{Lem.~2.2.3} and~\cite{Rog09}*{Prop.~3.20, Prop.~3.21}] \label{prop:BrepZge0}
  There are weight decompositions \[ B^\cy(\<x\>) = \textstyle \coprod_{i\ge0} B^\cy(\<x\>; i) \qqandqq B^\rep(\<x\>) = \coprod_{i\ge0} B^\rep(\<x\>; i)\]
  and weak equivalences of spaces with $\bT$-action
  \[ B^\cy(\<x\>;0) = {*},\quad B^\cy(\<x\>; i) \simeq S^1(i) \textrm{ if }i\geq 1, \text{and}\quad B^\rep(\<x\>; i) \simeq S^1(i) \textrm{ if }i\geq 0.\] 
Here $S^1(i)$ denotes the unit circle with the $\bT$-action given by
$z \cdot w = z^i w$ for $z \in \bT$ and $w \in S^1(i) \subset
\bC$. The repletion map $\rho \: B^\cy(\<x\>) \to B^\rep(\<x\>)$ is compatible with this decomposition. It induces a $\bT$-equivariant equivalence in positive weights and the inclusion of the basepoint in weight $0$.
\end{proposition}
\begin{corollary}\label{cor:THH-ThR-xi-zero}
The unit map $R \to \Th_{R}(\xi)$ induces an equivalence $\THH(R) \to \THH(\Th_{R}(\xi_*);0)$ of cyclotomic $\bE_{k-1}$ $\THH(R)$-algebras. 
\end{corollary}
\begin{proof}
It is sufficient to show that it is an equivalence of underlying $\bE_{k-1}$ $\THH(R)$-algebras with $\bT$-action. This follows from Proposition~\ref{prop:Bcy-THH-comparison-over-PicR} and the equivalence $* \to B^{\cy}(\<x\>; 0)$ from the last proposition. 
\end{proof}

The next proposition is the key ingredient for the proof of Theorem~\ref{thm:rep-res-Rbas}. It generalizes~\cite{RSS15}*{Prop.~6.11}.

\begin{proposition} \label{prop:rep-res-can}
Let $\xi \: \<x\> \to \Pic_R$ be an $\bE_k$-map with $k\geq 2$.  There is a natural repletion--residue cofiber sequence
\[
\THH^\xi(R) \xrightarrow{\trf} \THH(\Th_R(\xi_*))
	\xrightarrow{\rho} \THH(\Th_R(\xi_*), \xi)
	\xrightarrow{\res} \Sigma \THH^\xi(R)
\]
of $\bZ_{\ge 0}$-graded twisted cyclotomic $\THH(\Th_R(\xi_*))$-modules, with $\THH^\xi(R)$ viewed as an object in this category by restriction along $\THH(\coll)$.
\end{proposition}
\begin{proof}
  As explained in Subsection~\ref{subsec:collapse}, the face inclusion $f \: \{0\} \to \bZ_{\ge 0}$ induces a natural transformation $\id \to f_*f^*$. Applied to the $\bZ_{\geq 0}$-graded repletion map~\eqref{eq:rho_for_canonical}, this provides the left-hand commutative square of $\bZ_{\ge 0}$-graded twisted cyclotomic $\bE_{k-1}$ $\THH(R)$-algebras in the following diagram:
\begin{equation} \label{eq:cofrho}
\xymatrix@-1pc{
\THH(\Th_R(\xi_*)) \ar[rr]^-{\rho} \ar[d]
	&& \THH(\Th_R(\xi_*), \xi) \ar[rr]^-{\res} \ar[d]
	&& \cof(\rho) \ar[d] \\
\THH(\Th_R(\xi_*); 0) \ar[rr]^-{\rho_0}
	&& \THH(\Th_R(\xi_*), \xi; 0) \ar[rr]^-{\res_0}
	&& \cof(\rho_0)
}
\end{equation}
Using these $\bE_{k-1}$-algebra structures, we view the left-hand square as one of $\bZ_{\ge 0}$-graded twisted cyclotomic $\THH(\Th_R(\xi_*))$-modules and obtain the right-hand square by forming horizontal cofibers in that category. By construction, the left-hand vertical map in~\eqref{eq:cofrho} is equivalent to the composite of $\THH(\coll)$ and the equivalence of Corollary~\ref{cor:THH-ThR-xi-zero}.  Therefore, $\cof(\rho_0)$ is equivalent to $\Sigma \THH^\xi(R)$ as a $\bZ_{\ge 0}$-graded twisted cyclotomic $\THH(\Th_R(\xi_*))$-module, where $\Sigma\THH^\xi(R)$ carries the module structure obtained by restriction along $\THH(\coll)$.

With this identification, it suffices to show that the right-hand vertical map $\cof(\rho) \to \cof(\rho_0)$ is an equivalence. In weight $0$, this holds by definition. In weights $i>0$, the homotopy equivalences $B^{\cy}(\<x\>; i) \simeq B^\rep(\<x\>; i)$ from Proposition~\ref{prop:BrepZge0} and the compatible equivalences from Propositions~\ref{prop:Bcy-THH-comparison-over-PicR} and~\ref{prop:THHThRxistarxi} imply that the maps
\begin{equation} \label{eq:rhoi}
\rho_i \: \THH(\Th_R(\xi_*); i)
	\xrightarrow{\simeq} \THH(\Th_R(\xi_*), \xi; i)
\end{equation}
are equivalences of $\THH(R)$-modules. So both $\cof(\rho)$ and $\cof(\rho_0)$ are contractible in positive weights. 
\end{proof}

\begin{proof}[Proof of Theorem~\ref{thm:rep-res-Rbas}]
  We apply $\THH(A) \otimes_{\THH(\Th_R(\xi))}(-)$ to the repletion--residue cofiber sequence of Proposition~\ref{prop:rep-res-can}. This provides the desired cofiber sequence. Now suppose that $(R,\xi_*)$ is cyclotomically good, so that there is an equivalence $\THH(R) \simeq \THH^\xi(R)$ of cyclotomic $\THH(R)$-modules. Restricting along $\THH(\coll)$ and applying $\THH(A) \otimes_{\THH(\Th_R(\xi))}(-)$ gives an equivalence relating
  \[\THH(A) \otimes_{\THH(\Th_R(\xi))} \THH(R) \simeq \THH(A\otimes_{\Th_R(\xi)} R) = \THH(\modmod{A}{\bar\alpha}) \]
and $\THH^\xi(\modmod{A}{\bar\alpha}) = \THH(A) \otimes_{\THH(\Th_R(\xi))} \THH^\xi(R)$.
\end{proof}
\subsection{Identifying the repletion fiber as a \texorpdfstring{$\THH(R)$}{THH(R)}-module}
We turn to identifying $\THH^\xi(R)$ as a $\THH(R)$-module. While this result will be superseded by Proposition~\ref{prop:THHThRxixi0-Tact} below, the present analysis will be useful for the identification of the cyclotomic structure map of $\THH^\xi(R)$ and future calculations.

As before, we write $\xi^\gp\: M^\gp = \<x^{\pm1}\> \to \Pic_R$ for the  preferred $\bE_k$-extension of $\xi \: M = \<x\>  \to \Pic_R$. We write $\xi^{\gp}_* \: \<x^{\pm1}\> \to c^*\Pic_R$ for the graded version of our preferred $\bE_k$-extension of $\xi$ to $\<x^{\pm1}\>$ and $\xi^{\gp,\cy}_*\: B^\cy(\<x^{\pm1}\>) \to c^*\Pic_{\THH(R)}$ for the map induced by $\xi_*^\gp$ and the comparison map from Corollary~\ref{cor:Bcy-on-E_k-alg}(ii). Moreover, we write $P := \xi^\gp(x)$ and $P^* := \xi^\gp(x^{-1})$.

The coherence data for the $\bE_k$-map $\xi^\gp$ include a homotopy between two composites $M \times M \to \Pic_R$, which we can interpret as part of the coherence data for the map of cyclic bar constructions induced by $\xi_*^\gp$, namely a choice of a homotopy in the diagram 
\begin{equation}\label{eq:coherence-for-xi-gp} \xymatrix@-1pc{ B^\cy_1(\<x^{\pm 1}\>) \ar[rr] \ar[d]_{d_0} && B^\cy_1(c^*\Pic_R)  \ar[d]_{d_0}  \\  B^\cy_0(\<x^{\pm 1}\>) \ar[rr] && B^\cy_0(c^*\Pic_R)\,.}
\end{equation}
Likewise, the composite of the homotopy with the twist map determines the analogous coherence datum for the simplicial structure map $d_1$. Evaluated at $(x,x^{-1})$ and $(x^{-1},x)$, the homotopy provides two equivalences
\begin{align*}
\lambda_{x,x^{-1}} \:
&P \otimes_R P^* = \xi^\gp(x) \otimes_R \xi^\gp(x^{-1})
	\xrightarrow{\simeq} \xi^\gp(x \cdot x^{-1}) = \xi^\gp(1) = R \qquad\text{and}\\
\lambda_{x^{-1},x} \: 
&P^* \otimes_R P = \xi^\gp(x^{-1}) \otimes_R \xi^\gp(x)
	\xrightarrow{\simeq} \xi^\gp(x^{-1} \cdot x) = \xi^\gp(1) = R 
\end{align*}
witnessing that $P$ and $P^*$ are mutually inverse $R$-modules. 

\begin{lemma}\label{lem:repletion-fiber-as-THHR-module} The weight $0$ component \[\rho_0\: \THH(\Th_R(\xi_*); 0) \longto \THH(\Th_R(\xi_*), \xi; 0)\] of the repletion map is equivalent as a map of $\THH(R)$-modules to the composite
  \begin{equation}\label{eq:model-for-0-repl} \THH(R) \xrightarrow{\simeq} \THH(R)\wedge \{*\}_+ \longto \THH(R)\wedge S^1_+
  \end{equation}
  induced by the basepoint inclusion $\{*\} \to S^1$, and the repletion fiber  $\THH^\xi(R)$ is equivalent to $\THH(R)$ as a $\THH(R)$-module. 
\end{lemma}
\begin{proof}
  Proposition~\ref{prop:THHThRxistarxi} provides an equivalence
  \begin{equation}\label{eq:logTHHThRxi0-as-ThTHHR} \THH(\Th_R(\xi_*),\xi;0) \simeq \Th_{\THH(R)}(\xi^\rep_0) = \Th_{\THH(R)}(\xi^{\gp,\cy}_0)\end{equation}
  of $\THH(R)$-modules. To analyze the right-hand term, we note that Lemma~\ref{lem:Bcy-on-FunIS_cPic} and the use of the comparison map $|F(-)| \to F(|-|)$ in the proof of Proposition~\ref{prop:Bcy-symm-mon-nat}(ii) imply that for each $q$ we have a commutative diagram in which the vertical maps are induced by the inclusions $\eta_q$ of the $q$-simplices of the respective cyclic bar constructions into their realizations, that is, the colimits over $\Delta^{\op}$: 
  \begin{equation}\label{eq:incl-q-simplices} \xymatrix@-1pc{ B^\cy_q(\<x^{\pm 1}\>;0) \ar[rr] \ar[d]_{\eta_q} && B^\cy_q(\Pic_R) \ar[rr] \ar[d]_{\eta_q} && \Pic_{\THH_q(R)} \ar[d]^{\Pic_{\eta_q}} \\  B^\cy(\<x^{\pm 1}\>) \ar[rr] && B^\cy(\Pic_R) \ar[rr]  && \Pic_{\THH(R)}}\rlap{\,.}
  \end{equation}
  By the proof of Proposition~\ref{prop:BrepZge0} (see e.g.~\cite{Rog09}*{Prop.~3.21}), the inclusion of
  \[ \xymatrix{
(x^{-1}, x) \ar@<1ex>[r]^-{d_0} \ar@<-1ex>[r]_-{d_1} & (1)
}\]
into the $\Delta^\op$-diagram $B^\rep_\bullet(\<x\>; 0) = B^\cy_\bullet(\<x^{\pm 1}\>; 0)$ induces a homotopy equivalence $\Delta^1/\partial\Delta^1 \to B^\rep(\<x\>; 0)$ upon  passage to colimits. Composing with the maps to $\Pic_{\THH(R)}$  induced by $\xi^\gp_*$ that are specified in~\eqref{eq:incl-q-simplices}, and passing to Thom $\THH(R)$-modules, the coequalizer diagram above gets mapped to a diagram of the form
 \begin{equation}\label{eq:THHR-coequalizer-x-1x}\xymatrix{
\THH(R) \otimes_{R \otimes R} P^* \otimes P
	\ar@<1ex>[r]^-{d_0} \ar@<-1ex>[r]_-{d_1}
	& \THH(R) \otimes_R R \,.
      } 
    \end{equation}

Since $\eta_1 \: R \otimes R \simeq \THH_1(R) \to \THH(R)$ factors up to homotopy as $\eta_0 \circ d_0$, the face operator $d_0$ induces equivalences
\begin{multline}\label{eq:THHRRP-starP-identification}
\THH(R) \otimes_{R \otimes R} P^* \otimes P
	\simeq \THH(R) \otimes_R R \otimes_{R \otimes R} P^* \otimes P \\
	\simeq \THH(R) \otimes_R (P^* \otimes_R P)
	\simeq \THH(R) \otimes_R R = \THH(R) 
\end{multline}
identifying the first term in~\eqref{eq:THHR-coequalizer-x-1x}. (This equivalence is not a canonical one because we could also argue with $d_1$ and $\eta_1 \simeq \eta_0 \circ d_1$.)
   
To identify the coequalizer~\eqref{eq:THHR-coequalizer-x-1x}, we note that while the $1$-simplex $(x^{-1},x)$ represents a closed loop in $B^\cy_\bullet(\<x^{\pm 1}\>; 0)$, its image in $B^\cy_\bullet(\Pic_R)$ is an edge from $P^*\otimes_R P$ to $P\otimes_R P^*$. By the discussion of coherence homotopies before the lemma, the closed loop in $\Pic_R$ represented by the image of~$(x^{-1}, x)$ is composed of the edge $\lambda_{x^{-1}, x}$ (in reverse), the above-mentioned edge from $P^* \otimes_R P$ to $P \otimes_R P^*$, and the edge $\lambda_{x, x^{-1}}$. Upon passage to Thom spectra, this means that the maps $d_0$ and $d_1$ in~\eqref{eq:THHR-coequalizer-x-1x} are induced by $\lambda_{x^{-1},x}$ and $\lambda_{x,x^{-1}}\circ \tau$. Our assumption that $k\geq 2$ therefore implies that $d_0$ and $d_1$ are homotopic.

Using~\eqref{eq:logTHHThRxi0-as-ThTHHR}, the coequalizer of~\eqref{eq:THHR-coequalizer-x-1x} is equivalent to $\THH(\Th_R(\xi_*), \xi; 0)$. This equivalence is compatible with the corresponding identification for $(1) \to B^\cy(\<x\>; 0)$, and yields the following diagram of $\THH(R)$-modules, with vertical coequalizers and horizontal cofiber sequences. The upper right-hand entry is characterized by its suspension being the cofiber $\Sigma \THH^\xi(R)$ of $\rho_0$, hence is equivalent to the repletion fiber $\THH^\xi(R)$:
\begin{equation} \label{eq:fibrho}
\xymatrix@-1pc{
0 \ar@<-1ex>[d] \ar@<1ex>[d] \ar[rr]
	&& \THH(R) \otimes_{R \otimes R} P^* \otimes P
	\ar[rr]^-{\simeq} \ar@<-1ex>[d]_-{d_0} \ar@<1ex>[d]^-{d_1}
	&& \THH^\xi(R)
	\ar@<-1ex>[d] \ar@<1ex>[d] \\
\THH(R) \ar[d]_-{\simeq} \ar[rr]^-{\simeq}
	&& \THH(R) \otimes_R R \ar[rr] \ar[d]
	&& 0 \ar[d] \\
\THH(\Th_R(\xi_*); 0) \ar[rr]^-{\rho_0}
	&& \THH(\Th_R(\xi_*), \xi; 0) \ar[rr]^-{\res_0}
	&& \Sigma \THH^\xi(R)\rlap{\,.}
}
\end{equation}
With the homotopy $d_0\simeq d_1$, the diagram implies the claimed identification of~$\rho_0$. The statement about $\THH^\xi(R)$ follows by forming the cofiber of~\eqref{eq:model-for-0-repl}, or by combining the upper right-hand horizontal equivalence in~\eqref{eq:fibrho} with the combined equivalence in~\eqref{eq:THHRRP-starP-identification}.
\end{proof}

\subsection{Identifying the cyclotomic structure map on the repletion fiber}\label{subsec:identifying-Froebenius} We now study the cyclotomic structure map $\varphi^\xi_p \: \THH^\xi(R) \to \THH^\xi(R)^{tC_p}$ of the repletion fiber. Recall that in the formulation of~\cite{NS18}*{\S\S{}III.2--3}, the cyclotomic
structure map $\varphi_p \: \THH(A) \to \THH(A)^{tC_p}$ for an
$\bE_1$-ring~$A$ is constructed as a colimit of the Tate diagonal maps
$\Delta_p \: X \to (X^{\otimes p})^{tC_p}$, making squares
\begin{equation} \label{eq:Deltap-varphip}
\xymatrix@-1pc{
A^{\otimes 1+q} \ar[rr]^-{\eta_q} \ar[d]_-{\Delta_p}
	&& \THH(A) \ar[d]^-{\varphi_p} \\
(A^{\otimes p(1+q)})^{tC_p} \ar[rr]^-{\eta_{pq+p-1}^{tC_p}}
	&& \THH(A)^{tC_p}
}
\end{equation}
commute for all $q\ge0$.

We begin with the following refinement of Lemma~\ref{lem:repletion-fiber-as-THHR-module}:

\begin{proposition} \label{prop:noneq-phi'}
Non-equivariantly, i.e., ignoring the $\bT$-action, the cyclotomic
structure map $\varphi^\xi_p \: \THH^\xi(R) \to \THH^\xi(R)^{tC_p}$ is equivalent as a $\THH(R)$-module map to the standard cyclotomic structure map $\varphi_p \: \THH(R)
\to \THH(R)^{tC_p}$, for each prime~$p$.
\end{proposition}

\begin{proof}
The $C_p$-action on topological Hochschild homology is realized by a
simplicial action on its $p$-fold edgewise subdivision, see~\cite{BHM93}*{\S1},
\cite{LNR11}*{Lem.~3.4}, or \cite{NS18}*{\S{}III.2}.  This leads to the
following $C_p$-equivariant enrichment of~\eqref{eq:fibrho}:
\begin{equation} \label{eq:sdp-fibrho}
\xymatrix@-1pc{
0 \ar@<-1ex>[d] \ar@<1ex>[d] \ar[rr]
	&& \THH(R) \otimes_{R^{\otimes 2p}} (P^* \otimes P)^{\otimes p}
	\ar[rr]^-{\simeq} \ar@<-1ex>[d]_-{d_0} \ar@<1ex>[d]^-{d_1}
	&& \THH^\xi(R)
	\ar@<-1ex>[d] \ar@<1ex>[d] \\
\THH(R) \ar[d]_-{\simeq} \ar[rr]^-{\simeq}
	&& \THH(R) \otimes_{R^{\otimes p}} R^{\otimes p} \ar[rr] \ar[d]
	&& 0 \ar[d] \\
\THH(\Th_R(\xi_*); 0) \ar[rr]^-{\rho_0}
	&& \THH(\Th_R(\xi_*), \xi; 0) \ar[rr]^-{\res_0}
	&& \Sigma \THH^\xi(R)\rlap{\,.}
}
\end{equation}
In particular, the middle column is a coequalizer diagram of
$\THH(R)$-modules with $C_p$-action, and therefore yields a coequalizer
diagram after passage to $C_p$-Tate fixed points.

When applied with $A = R$ or $A = \Th_R(\xi_*^\gp)$, and evaluated in
weight $i=0$, the square~\eqref{eq:Deltap-varphip} leads to the following
commutative cube:
\[
\xymatrix@R-1pc@C-2pc{
R^{\otimes 1+q} \ar[rr]^-{\eta_q} \ar[dr] \ar[dd]_-{\Delta_p}
	&& \THH(R) \ar[dr] \ar[dd]^(0.3){\varphi_p}|-{\hole} \\
& (\Th_R(\xi_*^\gp)^{\otimes 1+q})_0 \ar[rr]^(0.35){\eta_q} \ar[dd]_(0.3){\Delta_p} 
	&& \THH(\Th_R(\xi_*^\gp); 0) \ar[dd]^-{\varphi_{p;0}}  \\
(R^{\otimes p(1+q)})^{tC_p} \ar[rr]|-{\hole} \ar[dr]
	&& \THH(R)^{tC_p} \ar[dr] \\
& (\Th_R(\xi_*^\gp)^{\otimes p(1+q)})^{tC_p}_0 \ar[rr]
	&& \THH(\Th_R(\xi_*^\gp); 0)^{tC_p}\rlap{\,.}
}
\]
When restricted along $R \to \Th_R(\xi_*^\gp)_0$ (for $q=0$) and $P^*
\otimes P \to \Th_R(\xi_*^\gp)^{\otimes 2}_0$ (for $q=1$), substituting
$\THH(\Th_R(\xi_*), \xi; 0)$ for $\THH(\Th_R(\xi_*^\gp); 0)$, this leads to
the following horizontal map of vertical coequalizer diagrams:
\[
\xymatrix@-1pc{
\THH(R) \otimes_{R^{\otimes 2}} P^* \otimes P
	\ar[rrrr]^-{\lambda_1(\varphi_p \otimes_{\Delta_p} \Delta_p)}
	\ar@<-1ex>[d]_-{d_0} \ar@<1ex>[d]^-{d_1}
	&&&& (\THH(R) \otimes_{R^{\otimes 2p}}
		(P^* \otimes P)^{\otimes p})^{tC_p}
	\ar@<-1ex>[d]_-{d_0} \ar@<1ex>[d]^-{d_1} \\
\THH(R) \otimes_R R
	\ar[rrrr]^-{\lambda_0(\varphi_p \otimes_{\Delta_p} \Delta_p)} \ar[d]
	&&&& (\THH(R) \otimes_{R^{\otimes p}} R^{\otimes p})^{tC_p} \ar[d] \\
\THH(\Th_R(\xi_*), \xi; 0)
	\ar[rrrr]^-{\varphi_{p;0}}
	&&&& \THH(\Th_R(\xi_*), \xi; 0)^{tC_p}\,.
}
\]
Here
\begin{align*}
\lambda_0 \: \THH(R)^{tC_p} \otimes_{(R^{\otimes p})^{tC_p}} (R^{\otimes p})^{tC_p}
	&\xrightarrow{\simeq} (\THH(R) \otimes_{R^{\otimes p}}
		R^{\otimes p})^{tC_p} \qquad\text{and}\\
\lambda_1 \: \THH(R)^{tC_p} \otimes_{(R^{\otimes 2p})^{tC_p}}
		((P^* \otimes P)^{\otimes p})^{tC_p}
	&\longto (\THH(R) \otimes_{R^{\otimes 2p}}
		(P^* \otimes P)^{\otimes p})^{tC_p}
\end{align*}
are instances of the lax structure maps for $C_p$-Tate fixed points.

In analogy with~\eqref{eq:THHRRP-starP-identification}, the subdivided face operator~$d_0$ induces equivalences
\begin{multline*}
\THH(R) \otimes_{R^{\otimes 2p}} (P^* \otimes P)^{\otimes p}
\simeq \THH(R) \otimes_{R^{\otimes p}} R^{\otimes p}
	\otimes_{R^{\otimes 2p}} (P^* \otimes P)^{\otimes p} \\
\simeq \THH(R) \otimes_{R^{\otimes p}} (P^* \otimes_R P)^{\otimes p}
\simeq \THH(R) \otimes_{R^{\otimes p}} R^{\otimes p} = \THH(R)
\end{multline*}
that respect the $C_p$-actions, hence also induce equivalences of
$C_p$-Tate fixed points.  Taken together, we obtain a commutative diagram
\[
\xymatrix@-1pc{
\THH^\xi(R) \ar[d]^-{\varphi^\xi_p}
	&& \THH(R) \otimes_{R^{\otimes 2}} (P^* \otimes P)
	\ar[ll]_-{\simeq} \ar[rr]^-{\simeq}
	\ar[d]^-{\lambda_1 (\varphi_p \otimes_{\Delta_p} \Delta_p)}
	&& \THH(R) \ar[d]^-{\varphi_p} \\
\THH^\xi(R)^{tC_p}
	&& (\THH(R) \otimes_{R^{\otimes 2p}} (P^* \otimes P)^{\otimes p})^{tC_p}
	\ar[ll]_-{\simeq} \ar[rr]^-{\simeq}
	&& \THH(R)^{tC_p}
}
\]
of $\THH(R)$-modules, exhibiting the commensurability of~$\varphi^\xi_p$
with $\varphi_p$.
\end{proof}
We will use the following consequence to combine this equivalence with the $\bT$-equivariant analysis of the repletion fiber to be addressed next. 
\begin{corollary} \label{cor:Frobenius-p-completion}
If $R$ is such that $\varphi_p \: \THH(R) \to \THH(R)^{tC_p}$ is
$p$-completion, then $\varphi^\xi_p \: \THH^\xi(R) \to \THH^\xi(R)^{tC_p}$
is also $p$-completion.  
\end{corollary}

\subsection{Identifying the repletion fiber as a \texorpdfstring{$\bT$}{T}-equivariant \texorpdfstring{$\THH(R)$}{THH(R)}-module}
We now study $\THH(\Th_R(\xi_*), \xi; 0) = \THH(\Th_R(\xi^\gp_*); 0)$ as a $\THH(R)$-module with $\bT$-action, using a similar strategy as in~\cite{RSS15}*{\S7}. 
\begin{definition}
For a map $f \: Y_\circ \to B^\cy_\circ(\<x^{\pm1}\>)$ of cyclic
sets, we define 
\[ Y \odot \Th_R(\xi^\gp) := \Th_{\THH(R)}(\xi^{\gp,\cy}_*\circ |f| \: |Y_\circ| \to B^\cy(\<x^{\pm1}\>) \to c^*\Pic_{\THH(R)})\,,\]
viewed as a $\THH(R)$-module with $\bT$-action. 
\end{definition}
By general properties of the Thom spectrum functor, $-\odot \Th_R(\xi^\gp)$ sends weak equivalences of cyclic sets to equivalences and (homotopy) colimits of cyclic sets over $B^\cy_\circ(\<x^{\pm1}\>)$ to (homotopy) colimits. It recovers~$\THH(\Th_R(\xi_*^\gp); i)$ when applied to the inclusion
$f \: B^\cy_\circ(\<x^{\pm1}\>; i) \to B^\cy_\circ(\<x^{\pm1}\>)\,$.

\begin{definition}
Let $S^\infty_\circ = B^\cy_\circ(\cG)$ be the cyclic bar
construction~\cite{Wal79}*{(2.3)} of the groupoid~$\cG$ with objects
$\up$ and $\down$ and non-identity morphisms $+1 : \down \to \up$
and $-1 : \up \to \down$.  Let $S^0_\circ = B^\cy_\circ(\cG^\delta)$
be the cyclic bar construction of its discrete subgroupoid $\cG^\delta$,
with objects $\up$ and $\down$ and only identity morphisms.
\[
\cG: \quad
\xymatrix{ \up \ar@(l,l)[d]_-{-1} \\ \down \ar@(r,r)[u]_-{+1} }
\qquad \qquad \qquad \qquad
\cG^\delta: \qquad
\xymatrix{ \up \\ \down }
\]
These have geometric realizations $S^\infty$ and $S^0$, respectively,
where the former is contractible and the latter has trivial $\bT$-action.
\end{definition}

\begin{definition}\label{def:cyclic-map-f}
Let $f \: S^\infty_\circ \to B^\cy_\circ(\<x^{\pm1}\>; 0)$ be the cyclic
map taking the $q$-simplex
\[
\xymatrix{
g_0 \ar[r]^-{i_1} & g_1 \ar[r] & \dots \ar[r] & g_{q-1} \ar[r]^-{i_q} & g_q
	\ar@(d,d)[llll]^-{i_0}
}
\]
to the $q$-simplex $(x^{i_0}, x^{i_1}, \dots, x^{i_q})$, with $g_0, \dots,
g_q \in \{\up, \down\}$, $i_0 = g_0 - g_q$, and $i_s = g_s - g_{s-1}$
for $1 \le s \le q$.
\end{definition}
\begin{lemma}
The cyclic map $f$ takes $S^0_\circ$ to $(1)$, and induces an isomorphism
$\bar f$ of $S^\infty_\circ/S^0_\circ$ with the cyclic subset of
$B^\cy_\circ(\<x^{\pm1}\>; 0)$ consisting of the $(x^{i_0}, \dots,
x^{i_q})$ with $i_0 + \dots + i_q = 0$ for which the nonzero $i_s$
alternate cyclically between $+1$ and $-1$.
\end{lemma}

\begin{proof}
The integer sequence $(i_0, \dots, i_q)$ drawn from $\{+1, 0, -1\}$
determines the object sequence $(g_0, \dots, g_q)$ drawn from $\{\up,
\down\}$, except when the latter sequence is constant: The first nonzero $i_s$ in $(i_0,\dots, i_q)$ gives the morphism $\down \to \up$ in $\cG$ if $i_s = 1$ and the morphism $\up \to \down$ if $i_s = -1$. This determines the domain and codomain of all other $i_t$ in the sequence and therefore $(i_0,\dots, i_q)$. 
\end{proof}

\begin{lemma}\label{lem:homotopy-pushout-cyclic-sets}
The square of cyclic sets
\begin{equation}\label{eq:homotopy-pushout-cyclic-sets}
\xymatrix@-1pc{
S^0_\circ \ar[rr] \ar[d]
	&& S^\infty_\circ \ar[d]^-{f} \\
(1) \ar[rr]
	&& B^\cy_\circ(\<x^{\pm1}\>; 0)
}
\end{equation}
is a (homotopy) pushout, in the sense that the induced $\bT$-map
\[
\bar f \: S^\infty/S^0 \longto B^\cy(\<x^{\pm1}\>; 0)
\]
is a weak equivalence.  
\end{lemma}

\begin{proof}
The cyclic map $f$ takes each of the $0$-simplices $(g_0) = \up$ and
$(g_0) = \down$ to $(1)$.  Moreover, it takes the $1$-simplex $(g_0,
g_1) = (\down, \up)$ to $(x^{-1}, x)$, hence induces the $\bT$-map~$\bar
f$ taking a generating loop in $S^\infty/S^0 \simeq S^1$ to a generating
loop in $B^\cy(\<x^{\pm1}\>; 0) \simeq S^1(0)$, which is therefore
a weak equivalence.
\end{proof}

We now apply $- \odot \Th_R(\xi^\gp)$ to~\eqref{eq:homotopy-pushout-cyclic-sets}, making the identifications $(1) \odot \Th_R(\xi^\gp) \simeq \THH(R)$ and $S^0 \odot \Th_R(\xi^\gp) \simeq \THH(R) \vee \THH(R)$. When applying ${- \odot \Th_R(\xi^\gp)}$ to the cyclic map $f$ of Definition~\ref{def:cyclic-map-f}, Proposition~\ref{prop:Bcy-THH-comparison-over-PicR}(i) implies that the codomain of the resulting map $f_*$ can be identified with $\THH(\Th_R(\xi^\gp_*); 0)$. This produces the following (homotopy) pushout of $\THH(R)$-modules with $\bT$-action
\begin{equation}\label{eq:homotopy-pushout-THHR-modules}\xymatrix@-1pc{
\THH(R) \vee \THH(R) \ar[rrrr]^{\up_* \vee\,\, \down_*} \ar[d]_-{\id \vee \id}
        &&&& S^\infty \odot \Th_R(\xi^\gp) \ar[d]^-{f_*} \\
\THH(R) \ar[rr]^-{\sim}
&&\THH(\Th_R(\xi_*); 0) \ar[rr]^-{\rho_0} && \THH(\Th_R(\xi_*^\gp); 0)\rlap{\,,}}
\end{equation}
where the bottom composite is induced by the unit map of $\Th_R(\xi^\gp)$ and therefore factors as the the weight $0$ component of the repletion map and the equivalence of Corollary~\ref{cor:THH-ThR-xi-zero}. 

We can now prove the following equivariant refinement of Lemma~\ref{lem:repletion-fiber-as-THHR-module}:
\begin{proposition} \label{prop:THHThRxixi0-Tact}
The weight $0$ component $\rho_0$ of the repletion map is equivalent as a map of $\THH(R)$-modules with $\bT$-action to the composite
  \begin{equation}\label{eq:eq-model-for-0-repl} \THH(R) \xrightarrow{\simeq} \THH(R)\wedge \{*\}_+ \to \THH(R)\wedge S^1(0)_+
  \end{equation}
  induced by the basepoint inclusion $\{*\} \to S^1$, and the repletion fiber  $\THH^\xi(R)$ is equivalent to $\THH(R)$ as a $\THH(R)$-module with $\bT$-action.
\end{proposition}
\begin{proof}
  We claim that the two $\bT$-equivariant maps $\up_*,\down_*$ in~\eqref{eq:homotopy-pushout-THHR-modules} are homotopic through $\bT$-equivariant maps. Assuming this, we can extend $\up_* \vee\,\, \down_*$ over $\THH(R)\wedge [0,1]_+$. The extension $ \THH(R)\wedge [0,1]_+ \to  S^\infty \odot \Th_R(\xi^\gp)$ is an equivalence because the inclusion $\up \to S^\infty$ induces an equivalence when applying $-\odot \Th_R(\xi^\gp)$. The statement about $\rho_0$ then follows from the homotopy pushout~\eqref{eq:homotopy-pushout-THHR-modules}, and from this we get that the fiber of $\rho_0$ is equivalent to $\THH(R)$ as a $\THH(R)$-module with $\bT$-action.

  To verify the claim about $\up_*,\down_*$, it is sufficient to find a homotopy between the corresponding maps $\up_*,\down_*\: \{*\} \to S^\infty$ in $(\cS_{/\Pic_{\THH(R)}})^{B\bT}$ and then apply $\Th_{\THH(R)}$. Using the augmentation of the codomain, this reduces to finding such a homotopy in spaces with $\bT$-action. Precomposing with the $\bT$-equivariant weak equivalence $E\bT \to \{*\}$, this in turn becomes a statement about $\bT$-equivariant maps $E\bT \to S^\infty$. Since the space of such maps is equivalent to the contractible space $(S^\infty)^{h\bT}$, any two such maps are $\bT$-equivariantly homotopic. 
\end{proof}
\subsection{Criteria for cyclotomically good pairs}\label{subsec:p-cyc-good-criteria}
Given a $\THH(R)$-module~$E$ we can generalize~\cite{Bou79} slightly
and say that a $\THH(R)$-module~$Z$ is \emph{$E_*$-acyclic} if $E
\otimes_{\THH(R)} Z \simeq 0$, a map~$X \to Y$ of $\THH(R)$-modules
is an \emph{$E_*$-equivalence} if its \mbox{(co-)fiber} is $E_*$-acyclic, and
$Y$ is \emph{$E_*$-local} if $\Map_{\THH(R)}(Z, Y) \simeq *$ for each
$E_*$-acyclic~$Z$.  The map $X \to Y$ is an \emph{$E_*$-localization}
of~$X$ if it is an $E_*$-equivalence and~$Y$ is $E_*$-local. For $\THH(R)$-modules with $\bT$-action, the notions $E_*$-acyclic and $E_*$-equiva\-lence are defined in terms of the underlying non-equivariant $\THH(R)$-modules. 

If $Z$ and~$Y$ are $\THH(R)$-modules with $\bT$-action, then the
$\bT$-equivariant mapping space from~$Z$ to~$Y$ has the homotopy type
of $\Map_{\THH(R)}(Z, Y)^{h\bT}$.  This space is contractible for all
$E_*$-acyclic~$Z$ with $\bT$-action (in which case we say that $Y$
is $\bT$-equivariantly $E_*$-local), if and only if the underlying
non-equivariant $\THH(R)$-module of~$Y$ is $E_*$-local.  A $\bT$-map $X
\to Y$ of $\THH(R)$-modules is thus a $\bT$-equivariant $E_*$-localization
if and only if the underlying map of non-equivariant $\THH(R)$-modules
is an $E_*$-localization.

When $E = \THH(R) \otimes \bS/p$, being $E_*$-local is the same as being
$p$-complete, and $E_*$-localization is the same as $p$-completion.

\begin{proof}[Proof of Proposition~\ref{prop:sufficient-criteria-cycl-good}]
  By Proposition~\ref{prop:THHThRxixi0-Tact}, $\THH(R) \simeq\THH^\xi(R)$ as $\THH(R)$-modules with $\bT$-action. Corollary~\ref{cor:Frobenius-p-completion} and the above discussion therefore imply that it is sufficient to show that $\varphi_p \: \THH(R) \to \THH(R)^{tC_p}$ is $p$-completion for every prime $p$ under the assumption of the proposition. This hypothesis is satisfied for $R = \bS$ by the theorems of Lin~\cite{Lin80} and Gunawardena~\cite{Gun81}, and for $R = MU$ by~\cite{LNR11}.

  A similar argument applies if $\varphi_p \: \THH(R) \to \THH(R)^{tC_p}$ is a smashing localization in the category of $\THH(R)$-modules with $\bT$-action. This applies when $R = H\bF$ with $\bF$ any perfect field of characteristic~$p$, since by~\cite{HM03}*{\S5} the $\bT$-map $\varphi_p \: \THH(\bF) \to \THH(\bF)^{tC_p}$ then inverts any choice of $\bT$-equivariant lift in $\pi_2 \TC^{-}(\bF) \cong W(\bF)$ of the degree~$2$ B{\"o}kstedt periodicity element~$\mu_0$ in $\pi_* \THH(\bF) \cong \bF[\mu_0]$.  Moreover, $\varphi_q$ is trivial for each prime $q \ne p$.  A similar argument applies for $\bF$ any perfectoid ring, as in~\cite{BMS19}*{Def.~4.18, Prop.~6.2}. 
\end{proof}

\section{The repletion--residue sequence in the multiple generator case}\label{sec:multiple-generator-case}
The purpose of this section is to generalize the repletion--residue sequence constructed in the previous section to $R$-based prelog $\bE_k$-rings where the monoid is a free commutative monoid of rank $r$. So throughout this section, we let $(\xi^i\: \<x_i\> \to \Pic_R)_{1\leq i\leq r}$ be a family of $\bE_k$-maps with $k\geq 2$ and let
\begin{equation}\label{eq:xi-as-product} \xi = \xi^1\cdots \xi^r\: \<x_1,\dots, x_r\> \to \Pic_R 
\end{equation}
be their product, viewed as a $\bZ_{\geq 0}^r$-graded object with the canonical grading. Prelog structures involving such a $\xi$ arise for example from Corollary~\ref{cor:prelog-from-generators}.

The next lemma shows that the condition that $\xi$ is a product can always be enforced, albeit up to losing one level of commutativity.
\begin{lemma} Let $\zeta \: \<x_1,\dots, x_r\> \to \Pic_R$ be an $\bE_k$-map. Then the composite \[\zeta|_{\<x_1\>} \cdots \zeta|_{\<x_r\>} \:  \<x_1,\dots, x_r\> = \<x_1\>\times \dots \times \<x_r\> \to \Pic_R\times \dots \times \Pic_R \to \Pic_R\]
  is homotopic to $\zeta$ through a $\bE_{k-1}$-maps.
\end{lemma}
\begin{proof}
  The underlying $\bE_1$-structure on $\zeta$ gives equivalences
  \[ \zeta|_{\<x_1\>}(x_1^{i_1}) \otimes_R \dots \otimes_R \zeta|_{\<x_r\>}(x_r^{i_r}) \xrightarrow{\simeq} \zeta(x_1^{i_1}\cdots x_r^{i_r})
  \]
  that define a homotopy from $\zeta^1\cdots \zeta^r$ to~$\zeta$ through
$\bE_0$-maps. The $\bE_k$-case follows by Dunn additivity. 
\end{proof}
Returning to the situation of~\eqref{eq:xi-as-product}, we can use the external product~\eqref{eq:external-product} to decompose the (canonical log) $\THH$ of $\Th_R(\xi^1 \cdots \xi^r)$: 
\begin{lemma}\label{lem:external-product-THH-log-THH-xi-i}
  For a family $(\xi^i\: \<x_i\> \to \Pic_R)_{1\leq i\leq r}$  of $\bE_k$-maps with $k\geq 2$, the symmetric monoidal structures of $\THH$ and $\Th_R$ induce compatible equivalences
  \begin{multline}\label{eq:external-product-THH-xi-i} \THH(\Th_R(\xi^1_*)) \otimes_{\THH(R)} \dots  \otimes_{\THH(R)}\THH(\Th_R(\xi^r_*))\\  \xrightarrow{\simeq} \THH(\Th_R(\xi^1_* \cdots \xi^r_*)) \qquad\text{and}
  \end{multline}
  \begin{multline}\label{eq:external-product-logTHH-xi-i}
    \THH(\Th_R(\xi^1_*),\xi^1) \otimes_{\THH(R)} \dots  \otimes_{\THH(R)}\THH(\Th_R(\xi^r_*),\xi^r) \\ \xrightarrow{\simeq} \THH(\Th_R(\xi^1_* \cdots \xi^r_*), \xi^1\cdots \xi^r) \qquad  \qquad
  \end{multline}
  of $\bZ_{\geq 0}^r$-graded twisted cyclotomic $\THH(R)$-modules. 
\end{lemma}
\begin{proof}
For the first equivalence, we apply the symmetric monoidal structure of $\Th_R\: \cS_{/ \Pic_R}\to \Mod_R$ to $\xi^1_*,\dots,\xi^r_*$ and use Lemma~\ref{lem:THH-commutes-w-external-product}. For the second, we argue analogously with $(\xi^1_*)^\gp,\dots,(\xi^r_*)^\gp$ and restrict the resulting $\bZ^r$-graded object to $\bZ_{\geq 0}^r$. As in~\eqref{eq:rho_for_canonical}, we get compatible repletion maps. 
\end{proof}
We will use the following notion for the bookkeeping of the cofiber sequences we are after: 
\begin{definition}
  Let $\cC$ be a symmetric monoidal stable $\infty$-category and let $[2]$ denote the poset $\{0 < 1 < 2\}$. An \emph{$r$-dimensional cube of cofiber sequences} in $\cC$ is a functor $C(-)\: [2]^r \to \cC$ such that for any sequence $(a_1,\dots, a_{i-1},a_{i+1},\dots ,a_r)$ in $\{0,1,2\}^{r-1}$ the diagram
  \[ C(a_0,\dots,a_{i-1},-,a_{i+1},\dots, a_r) \: [2] \to \cC\] is a cofiber sequence in $\cC$. 
\end{definition}
\begin{construction} We write
\begin{equation}\label{eq:THH-ThR-xi1-xin-cube} \THH(\Th_R(\xi_*);-) \: [2]^r \to \TwCyc\Fun(\bZ^r_{\geq 0}, \Mod_{\THH(R)})
\end{equation}
for the $r$-fold external product of the cofiber sequences of $\bZ_{\geq 0}$-graded twisted cyclotomic $\THH(R)$-modules
\[ \THH(\Th_R(\xi^i_*)) \xrightarrow{\rho} \THH(\Th_R(\xi^i_*),\xi^i) \xrightarrow{\res} \Sigma \THH^{\xi^i}(R) \]
from Proposition~\ref{prop:rep-res-can}. Via the $\bE_k$-equivalence~\eqref{eq:external-product-THH-xi-i} and the $\THH(\Th_R(\xi^i_*))$-module structures from Proposition~\ref{prop:rep-res-can}, we can view $\THH(\Th_R(\xi_*);-)$ as a diagram of $\bZ_{\geq 0}^r$-graded twisted cyclotomic $\THH(\Th_R(\xi_*))$-modules.
\end{construction}

To describe $\THH(\Th_R(\xi_*);-)$ in more detail, we introduce the notation 
\begin{equation}
S_j = S_j(a_1,\dots, a_r) = \{ i \mid a_i = j\} \subseteq \{1,\dots, r\}
\end{equation}
for $0\leq j \leq 2$ and an index $(a_1,\dots, a_r) \in \{0,1,2\}^r$. Then $\THH(\Th_R(\xi_*);a_1,\dots, a_r)$ is equivalent to
  \begin{equation}\label{eq:external-product-canonical}
    \textstyle \bigotimes_{i \in S_0} \THH(\Th_R(\xi^i_*)) \otimes \bigotimes_{i \in S_1} \THH(\Th_R(\xi^i_*),\xi^i) \otimes \bigotimes_{i \in S_2} \Sigma \THH^{\xi^i}\!(R) 
  \end{equation}
  where all $\otimes$-products are formed over $\THH(R)$.

\begin{definition}
For an $R$-based prelog $\bE_k$-ring
\[ (A,\xi=\xi^1\cdots\xi^r \: \<x_1,\dots,x_r\> \to \Pic_R, \bar\alpha \: \Th_R(\xi) \to A)\]
with $\xi$ as in \eqref{eq:xi-as-product}, we let \[\THH(A,\xi,\bar\alpha;-) \: [2]^r \to \Cyc\Mod_{\THH(A)}\]
be the diagram of cyclotomic $\THH(A)$-modules obtained from (the underlying ungraded object of) the external product~\eqref{eq:THH-ThR-xi1-xin-cube} by base change along the $\bE_{k-1}$-map $\THH(\bar\alpha) \: \THH(\Th_R(\xi)) \to \THH(A)$.
\end{definition}
For an $S \subseteq \{1,\dots,r\}$, we define $\bar\alpha_{S,\emptyset}$ to be the composite $\bE_k$-map 
\[ \Th_R(\textstyle\prod_{i\in S} \xi^i) \xrightarrow{\incl} \Th_R(\xi^1\cdots\xi^r) \xrightarrow{\bar\alpha} A.\]
Similarly to~\eqref{eq:Amodbaralphax}, we write
\[ \modmod{A}{\bar\alpha(S)} = A \otimes_{\Th_R\textstyle(\prod_{i \in S}\xi^i)} R \]
for the $\bE_{k-1}$-algebra formed using $\bar\alpha_{S,\emptyset}$ and the collapse map (where we understand $\modmod{A}{\bar\alpha(\emptyset)} = A$). For subsets $S, T \subseteq \{1,\dots,r\}$, we define $\bar\alpha_{S,T}$ to be the composite
\[ \Th_R(\textstyle\prod_{i\in S} \xi^i) \to A\to \modmod{A}{\bar\alpha(T)}\]
of $\bar\alpha_{S,\emptyset}$ and the canonical map. 

\begin{theorem}\label{thm:identification-of-hty-cofiber-cube}
  The diagram $\THH(A,\xi,\bar\alpha;-)$ is an $n$-dimensional cube of cofiber sequences of cyclotomic $\THH(A)$-modules. For an object $(a_1,\dots,a_r)$ of~$[2]^r$, the term $\THH(A,\xi,\bar\alpha;a_1,\dots,a_r)$ is equivalent as a $\THH(A)$-module with $\bT$-action to $\Sigma^{|S_2|} \THH(\modmod{A}{\bar\alpha(S_2)})$ if $S_1(a_1,\dots,a_r) = \emptyset$, and to
  \begin{equation}\label{eq:gen-term-in-hty-cofiber-cube} \textstyle  \Sigma^{|S_2|} \THH(\modmod{A}{\bar\alpha(S_2)})\otimes_{\THH(\Th_R(\prod_{i\in S_1}\xi^i))} \THH(\Th_R(\prod_{i\in S_1}\xi^i),\prod_{i\in S_1}\xi^i)
  \end{equation}
  otherwise. 
If $(R,\xi^i_*)$ is in addition cyclotomically good for all $1\leq i \leq r$, then these are equivalences of cyclotomic $\THH(A)$-modules.
\end{theorem}
\begin{remark}\label{rem:log-THH-for-E1-quotients}
When $k\geq 3$, then $\modmod{A}{\bar\alpha(S_2)}$ is an $\bE_2$-ring and~\eqref{eq:gen-term-in-hty-cofiber-cube} is equivalent to $\Sigma^{|S_2|}\THH(\modmod{A}{\bar\alpha(S_2)},\prod_{i \in S_1}\xi^i,\bar\alpha_{S_1,S_2})$ in the sense of Definition~\ref{def:logTHHAxi}. When $k=2$ and $S_2\neq \emptyset$, then $\modmod{A}{\bar\alpha(S_2)}$ is a priori only $\bE_1$. We view~\eqref{eq:gen-term-in-hty-cofiber-cube} as an extension of the definition of log $\THH$ to this situation, and use the same notation.  
\end{remark}

If the $(R,\xi^i_*)$ are all cyclotomically good, then Corollary~\ref{cor:rep-res-TC} generalizes and we obtain a similar $r$-dimensional cube of cofiber sequences of $\TC$-terms.  

\begin{proof}[Proof of Theorem~\ref{thm:identification-of-hty-cofiber-cube}]
We begin by observing that if $S,T \subseteq \{1,\dots,r\}$ are disjoint subsets with $S \cup T =  \{1,\dots,r\}$, then using (symmetric) monoidality of $\Th_R$ we get an equivalence
\begin{equation}\label{eq:subset-collapse}
  \textstyle A \otimes_{\Th_R(\xi^1\cdots \xi^r)} \Th_R(\prod_{i\in S}\xi^i) \simeq  A \otimes_{\Th_R(\prod_{i\in T}\xi^i)} R.
  \end{equation}
  We first assume $S_1\neq \emptyset$. Applying $\THH(A)\otimes_{\THH(\Th_R(\xi))}-$ to~\eqref{eq:external-product-canonical} and using the equivalences $\THH(R) \simeq \THH^{\xi^i}(R)$ from Proposition~\ref{prop:THHThRxixi0-Tact} as well as~\eqref{eq:external-product-THH-xi-i} and~\eqref{eq:external-product-logTHH-xi-i}, the term is equivalent to 
  \[
\textstyle \Sigma^{|S_2|} \THH(A)\otimes_{\THH(\Th_R(\xi))}\left( \bigotimes_{i \in S_0} \THH(\Th_R(\xi^i)) \otimes \bigotimes_{i \in S_1} \THH(\Th_R(\xi^i),\xi^i)\right). 
\]
Replacing $\THH(\Th_R(\xi^i),\xi^i)$ by $\THH(\Th_R(\xi^i))\otimes_{\THH(\Th_R(\xi^i))}\THH(\Th_R(\xi^i),\xi^i)$, Lemma~\ref{lem:external-product-THH-log-THH-xi-i} and the equivalence~\eqref{eq:subset-collapse} with the complementary subsets $S_0\cup S_1$ and $S_2$ imply that the latter term is equivalent to the claimed one. An analogous argument applies when $S_1=\emptyset$.

When in addition each $(R,\xi^i_*)$ is cyclotomically good, we can argue with the definition of cyclotomically good instead of Proposition~\ref{prop:THHThRxixi0-Tact}. 
\end{proof}
The maps in the cubical diagram provided by the theorem are induced by the repletion and residue maps from Theorem~\ref{thm:rep-res-Rbas}. 
\begin{example}\label{ex:general-two-gen-repl-res-seq} We let $r=2$ and assume the $(R,\xi^i_*)$ to be cyclotomically good. Abbreviating all indices $\{i\}$ to $i$, the theorem provides a diagram
\[
\xymatrix@-.8pc{
\THH(A) \ar[r]^-{\rho} \ar[d]_-{\rho}
	& \THH(A, \xi^1, \bar\alpha_{1,\emptyset}) \ar[r]^-{\res} \ar[d]_-{\rho}
	& \Sigma \THH(\modmod{A}{\bar\alpha(1)}) \ar[d]_-{\Sigma \rho} \\
\THH(A, \xi^2, \bar\alpha_{2,\emptyset}) \ar[r]^-{\rho} \ar[d]_-{\res}
	& \THH(A, \xi^1 \cdot \xi^2, \bar\alpha) \ar[r]^-{\res} \ar[d]_-{\res}
	& \Sigma \THH(\modmod{A}{\bar\alpha(1)}, \xi^2, \bar\alpha_{2,1})
		\ar[d]_-{\Sigma \res} \\
\Sigma \THH(\modmod{A}{\bar\alpha(2)}) \ar[r]^-{\Sigma \rho}
	& \Sigma \THH(\modmod{A}{\bar\alpha(2)}, \xi^1, \bar\alpha_{1,2})
		\ar[r]^-{\Sigma \res}
	& \Sigma^2 \THH(\modmod{A}{\bar\alpha})
}
\]
of horizontal and vertical cofiber sequences of cyclotomic $\THH(A)$-modules, where the lower right-hand square commutes up to a sign $-1$ and the three other squares commute. Here the sign comes from two different ways to identify a smash product of suspensions with an iterated suspension, and the two terms mapping to $\Sigma^2 \THH(\modmod{A}{\bar\alpha(1,2)})$ are extending Definition~\ref{def:logTHHAxi} to the case where the underlying $\bE_k$-ring is only $\bE_1$ in the way explained in Remark~\ref{rem:log-THH-for-E1-quotients}.
Passing to $\TC$, we get an analogous diagram of horizontal and vertical cofiber sequences of $\TC(A)_p$-modules \[
\xymatrix@-.8pc{
\TC(A)_p \ar[r]^-{\rho} \ar[d]_-{\rho}
	& \TC(A, \xi^1, \bar\alpha_{1,\emptyset})_p \ar[r]^-{\res} \ar[d]_-{\rho}
	& \Sigma \TC(\modmod{A}{\bar\alpha(1)})_p \ar[d]_-{\Sigma \rho} \\
\TC(A, \xi^2, \bar\alpha_{2,\emptyset})_p \ar[r]^-{\rho} \ar[d]_-{\res}
	& \TC(A, \xi^1 \cdot \xi^2, \bar\alpha)_p \ar[r]^-{\res} \ar[d]_-{\res}
	& \Sigma \TC(\modmod{A}{\bar\alpha(1)}, \xi^2, \bar\alpha_{2,1})_p
		\ar[d]_-{\Sigma \res} \\
\Sigma \TC(\modmod{A}{\bar\alpha(2)})_p \ar[r]^-{\Sigma \rho}
	& \Sigma \TC(\modmod{A}{\bar\alpha(2)}, \xi^1, \bar\alpha_{1,2})_p
		\ar[r]^-{\Sigma \res}
	& \Sigma^2 \TC(\modmod{A}{\bar\alpha(1,2)})_p\rlap{\,.}
}
\]
\end{example}

\begin{example}
In the specific case where the prelog ring is $(\ell,\< p, v_1\>)$ from Definition~\ref{def:ku-ell-BPn}, Example~\ref{ex:general-two-gen-repl-res-seq} specializes to diagram~\eqref{eq:fraction-field-3x3-diagram-intro} from the introduction. The collapse of the complements of the faces $\< p \>, \< v_1\>$, and $\< p, v_1\>$ appearing in that diagram was described in Examples~\ref{ex:collapse-of-p-or-v1-in-ell} and~\ref{ex:collapse-of-p-and-v1-in-ell}. 
\end{example}

\section{Even-periodic sphere spectra}
The purpose of this section is to identify the terms in the log $\THH$ and log $\TC$ repletion--residue sequences for even-periodic sphere spectra.

Let $d \geq 0$. In what follows, we write $\xi = \xi_{2d} \: \bZ_{\geq 0} \to \Pic_\bS$ for the $r = 1$ instance of the $\bE_2$-map~\eqref{eq:map-defining-Sx1xr} and $\bS[x] = \Th_\bS(\xi)$ for the associated $\bE_2$-ring. Moreover, we write $\xi^\gp\: \bZ \to \Pic_\bS$ for the extension of $\xi$ arising from Construction~\ref{cons:xigp} and $\bS[x^{\pm 1}] = \Th_\bS(\xi)$ for the associated $\bE_2$-ring. 
Then
\[  \textstyle \bS[x^{\pm1}] \simeq \bigvee_{i\in\bZ} \bS^{2di} \qqandqq \bS[x]\simeq \bigvee_{i\ge0} \bS^{2di}, \]
meaning that these $\bE_2$-rings are models for the $2d$-periodic and non-negative $2d$-periodic sphere spectra. They come with canonical $\bZ_{\geq 0}$- and $\bZ$-gradings, and give rise to the canonical prelog $\bE_2$-ring $(\bS[x],\< x \>) := (\bS[x],\xi,\id)$ with $\xi$ the underlying total object of the $\bZ_{\geq 0}$-graded $\xi_*$ (see Example~\ref{ex:non-neg-even-periodic-sphere-sp}).

Before analyzing log $\THH$ and log $\TC$ in this example, we return to the question whether $\bS[x]$ admits an $\bE_3$-structure (compare Remark~\ref{rem:E3-structure-on-Sx}).

\begin{proposition}\label{prop:obstructions-E3-structures-Sx} If $p$ is odd and $p\nmid d$, then neither $\bS[x]$ nor its $p$-localization admit an $\bE_3$-ring structure extending the given $\bE_1$-ring structure.  
\end{proposition}
\begin{proof} We set $n = 2d$. Writing $\bE_kX = \bigvee_{j\ge0} C_{k,j} X$ for the free $\bE_k$-ring on $X$, we get canonical $\bE_1$- and $\bE_3$-maps $\bE_1S^{n}\to \bE_3 S^{n}\to \bE_\infty S^{n}$ given in operadic weight $j= 1$ by the identity map of $S^n$ and in operadic weight $j = p$ by
  \[ S^{np} = C_{1,p} S^{n} \to C_{3,p} S^{n} \to C_{\infty, p} S^{n} = (S^{n})^{\otimes p}_{h\Sigma_p}. \] 
  If  the free $\bE_1$-structure on $\bE_1 S^{n}$ extends to an $\bE_3$-structure, then $\bE_1S^{n}\to \bE_3 S^{n}$ admits a retraction $r$.

To exclude the existence of this $r$, we argue with Cohen's mod $p$ homology operations defined in~\cite{CLM76}*{\S{}IV.1}. Let $e \in H_n(S^n)$ be a generator.  Then $Q^{d}(e) = e^p$ generates $H_{pn}(-)$ at the bottom of operadic weight $p$ in each case.  The ``top'' operation $\xi_2$ is defined for $\bE_3$-algebras, and $ \xi_2(e)  \in  H_{pn+2p-2}( C_{3,p} S^n )$. If $P^1_* \xi_2(e)  \in  H_{pn}( C_{3,p} S^n )$ is nonzero, then no retraction $r$ can exist, since $P^1_*$ commutes with spectrum maps and $P^1_* = 0$ acting on $H_*( C_{1,p} S^n )$. To see that $P^1_* \xi_2(e) \ne 0$, we can use that $\xi_2(e) = Q^{d+1}(e)$ for
$\bE_\infty$-algebras, and instead show that $Q^{d+1}(e)  \in  H_{pn+2p-2}( C_{\infty, p} S^n )$ has nonzero image $P^1_* Q^{d+1}(e)  \in  H_{pn}( C_{\infty, p} S^n)$. When $p$ does not divide $d$, this now follows from the original Nishida relation $P^1_* Q^{d+1} = d Q^{d} P^0_*$, see e.g.~\cite{CLM76}*{Thm.~IV.1.1(6)}, since $P^0_*(e) = e$ and $Q^{d}(e) = e^p$.
\end{proof}

\subsection{Log \texorpdfstring{$\THH$}{THH} of even-periodic sphere spectra}
We get maps
\begin{equation}\label{eq:THHSP}
\THH(\bS[x]) \xrightarrow{\rho} \THH(\bS[x],\<x\>)
\xrightarrow{\tilde\gamma} \THH(\bS[x^{\pm 1}] )
\end{equation}
between (the total objects of) $\bZ_{\ge0}$-, $\bZ_{\ge0}$- and $\bZ$-graded $L_p$-twisted cyclotomic $\bE_1$-rings, where $\rho$ is the repletion map and $\tilde\gamma$ is induced by the inclusion of the non-negatively graded part of $ \THH(\bS[x^{\pm 1}])$. 
\begin{proposition}\label{prop:THHSP}
As a sequence of graded spectra with $\bT$-action,~\eqref{eq:THHSP} is equivalent to
\begin{multline*}
\textstyle\bS \vee \bigvee_{i>0} \bT_+ \wedge_{C_i} (\bS^{2d})^{\otimes i}
	\longto \bS[S^1(0)] \vee
	\bigvee_{i>0} \bT_+ \wedge_{C_i} (\bS^{2d})^{\otimes i} \\
	\longto
\textstyle	\bigvee_{i<0} \bT_+ \wedge_{C_{-i}} (\bS^{-2d})^{\otimes -i}
	\vee \bS[S^1(0)] \vee
	\bigvee_{i>0} \bT_+ \wedge_{C_i} (\bS^{2d})^{\otimes i} \,.
\end{multline*}
In each case, $C_{|i|} \subset \bT$ acts by cyclically permuting the
$|i|$ smash factors.
\end{proposition}
\begin{proof}
  Propositions~\ref{prop:Bcy-THH-comparison-over-PicS} and~\ref{prop:THHThRxistarxi} assert that the three terms in~\eqref{eq:THHSP} are equivalent as $\bZ_{\ge0}$-, $\bZ_{\ge0}$- and $\bZ$-graded $\bE_1$-rings with $\bT$-action to the Thom spectra over $B^\cy(\bZ_{\ge0}) \simeq (0) \sqcup \coprod_{i>0} S^1(i)$, $B^\rep(\bZ_{\ge0}) \simeq \coprod_{i\ge0} S^1(i)$ and~$B^\cy(\bZ) \simeq \coprod_{i \in \bZ} S^1(i)$, respectively. Here we use the notation~$S^1(i)$ from Proposition~\ref{prop:BrepZge0}.

  In weight $i > 0$, we start with the $C_i$-equivariant map $\bS[x]^{\otimes i} \to \THH(\bS[x])$ obtained as the inclusion of the $0$-simplices in the $i$-fold edgewise subdivision~\cite{BHM93}*{\S1}. The restriction along $(\bS^{2d})^{\otimes i} \to \bS[x]^{\otimes i}$ factors through the inclusion  $ \THH(\bS[x];i) \to \THH(\bS[x])$, thus giving a $C_i$-map $(\bS^{2d})^{\otimes i} \to \THH(\bS[x];i)$ that is adjoint to a $\bT$-equivariant map 
  \begin{equation}\label{eq:THHbSxwi-identification} \bT_+ \wedge_{C_i} (\bS^{2d})^{\otimes i} \to \THH(\bS[x];i).\end{equation}
  Writing $\sd_i$ for the $i$-fold edgewise subdivision, we can obtain the map~\eqref{eq:THHbSxwi-identification} under the equivalence of Proposition~\ref{prop:Bcy-THH-comparison-over-PicS} by extending the $C_i$-equivariant inclusion $\{(1,\dots,1)\} \to \sd_i(B^\cy(\bZ_{\geq 0};i))$ to a $\bT$-equivariant map $\bT \times_{C_i}\{(1,\dots,1)\} \to \sd_i(B^\cy(\bZ_{\geq 0};i))$ and then applying $\Th_\bS$. Since the latter map is the weak equivalence from Proposition~\ref{prop:BrepZge0}, it follows that~\eqref{eq:THHbSxwi-identification} is an equivalence. 

In weight $i < 0$, we invoke the group isomorphism $\bZ \to \bZ, 1 \mapsto -1$, which reverses the $\bZ$-grading.  Its composite with $\xi_*^\gp \: \bZ \to \Pic_\bS$ gives an $\bE_2$-map $\bar\xi_*^\gp \: \bZ \to \Pic_\bS$ sending $1$ to $\bS^{-2d}$.  Let $\bar\xi_* \: \bZ_{\ge0} \to \Pic_\bS$ be the restricted map.  Then
\begin{align*}
  \THH(\bS[x^{\pm 1}]; i) &= \THH(\Th_\bS(\xi_*^\gp); i)
	\simeq \THH(\Th_\bS(\bar\xi_*^\gp); -i) \\
&\simeq \THH(\Th_\bS(\bar\xi_*); -i)
	\simeq \bT_+ \wedge_{C_{-i}} (\bS^{-2d})^{\otimes -i}
\end{align*}
as a spectrum with $\bT$-action by the weight $i > 0$ case above.

In weight $i = 0$, the assertion follows from Proposition~\ref{prop:THHThRxixi0-Tact}.
\end{proof}

The following is a homological incarnation of Proposition~\ref{prop:THHSP}, capturing an algebraic image of the $\bE_1$-ring structures and $\bT$-actions.

\begin{proposition}
There are B{\"o}kstedt spectral sequences
\begin{align*}
E^2_{**} = \HH_*(H_*(\bS[x])) &= \HH_*(\bZ[x])
	= \bZ[x] \otimes \Lambda(dx) \\
	&\Longrightarrow H_*(\THH(\bS[x])) \\
E^2_{**} = \HH_*(H_*(\bS[x]), \<x\>) &= \HH_*(\bZ[x], \<x\>)
	= \bZ[x] \otimes \Lambda(\dlog x) \\
	&\Longrightarrow H_*(\THH(\bS[x],\<x\>)) \\
E^2_{**} = \HH_*(H_*(\bS[x^{\pm 1}])) &= \HH_*(\bZ[x^{\pm1}])
	= \bZ[x^{\pm1}] \otimes \Lambda(\dlog x) \\
	&\Longrightarrow H_*(\THH(\bS[x^{\pm 1}]))
\end{align*}
with classes of (Hochschild, homological) bidegree $\|x\| = (0,2d), \|dx\|=(1,2d)$, and $\|\dlog x\|=(1,0)$. 
These are algebra spectral sequences that collapse at the $E^2 = E^\infty$-terms, with no additive or multiplicative extensions.  The repletion map~$\rho$ takes $dx$ to $x \cdot \dlog x$.  The $\bT$-actions induce the derivations satisfying $x \mapsto dx$, $dx \mapsto 0$, $x \mapsto x \cdot \dlog x$ and $\dlog x \mapsto 0$.
\end{proposition}
\begin{proof}
Applying homology to the skeleton filtration of the geometric realization defining $\THH(\bS[x^{\pm 1}])$ gives the $(E^1, d^1)$-term in the third case, which we identify with the normalized Hochschild complex for~$H_*(\bS[x^{\pm 1}]) = \bZ[x^{\pm1}]$.  Hence its $E^2$-term is given by the stated Hochschild homology, with $\dlog x \in \HH_1(\bZ[x^{\pm1}])$ the class of the Hochschild $1$-cycle $x^{-1} \otimes x$.  The $\bE_2$-ring structure on $\bS[x^{\pm 1}]$ gives a map $\bS[x^{\pm 1}] \otimes \bS[x^{\pm 1}] \to \bS[x^{\pm 1}]$ of $\bE_1$-rings, inducing a pairing of spectral sequences given at the $E^1$-terms by the shuffle product and at the $E^2$-terms by the usual product in the Hochschild homology of a (graded) commutative ring.  This shows that we have an algebra spectral sequence.

There is no room for $d^r$-differentials with $r\ge2$, so $E^2 = E^\infty$.  Since the $E^\infty$-term is free abelian in each bidegree, there is no room for additive extensions.  The only possible multiplicative extension concerns the value of $\dlog x \cdot \dlog x$ in $\bZ\{x\}$ when $d=1$. But since $\THH(\bS[x^{\pm 1}])$ is a $\bZ$-graded $\bE_1$-ring, $(\dlog x)^2$ must lie in the weight~$0$ summand, while $x$ lies in the weight~$1$ summand.

Restricting to weights~$\ge0$ gives the second B{\"o}kstedt spectral sequence, with the logarithmic Hochschild homology $\HH_*(\bZ[x], \<x\>) \cong \bZ[x] \otimes \Lambda(\dlog x)$ defined as in~\cite{Rog09}*{Def.~3.23}.

Repeating the argument with $\bS[x]$ in place of~$\bS[x^{\pm 1}]$ gives the first B{\"o}kstedt spectral sequence, now with $\HH_1(\bZ[x]) = \bZ[x] \{dx\}$ generated by the homology class of the $1$-cycle $1 \otimes x$.  The chain level product $x \cdot (x^{-1} \otimes x) = 1 \otimes x$ then implies the relation $x \cdot \dlog x = \rho(dx)$ in $\HH_1(\bZ[x], \<x\>)$ and $\HH_1(\bZ[x^{\pm1}])$.

The $\bT$-action on $\THH(\bS[x^{\pm 1}])$ arising from the cyclic structure takes each stage of the skeleton filtration into the next, sending the Hochschild $0$-cycle $x^i \in E^1_{0,2di}$ to the $1$-cycle $1 \otimes x^i \in E^1_{1,2di}$, and induces the Connes $B$-operator on the $E^2$-terms.  In particular, this operator sends~$x$ to $x \cdot \dlog x$, and squares to zero.  It also acts as a derivation, e.g.~because the filtered pairing $\THH(\bS[x^{\pm 1}]) \otimes \THH(\bS[x^{\pm 1}]) \to \THH(\bS[x^{\pm 1}])$ takes the diagonal $\bT$-action on the source to the given $\bT$-action on the target.  It follows that the $B$-operator agrees with the algebraic differential~$d$, mapping the generator $x^i \in H_{2di}(\THH(\bS[x^{\pm 1}]))$ to $d(x^i) = i x^{i-1} \, dx = i x^i \, \dlog x$, which is $i$ times the generator $x^i \, \dlog x \in H_{2di+1}(\THH(\bS[x^{\pm 1}]))$, and mapping $x^i \, \dlog x$ to $0$ for all~$i \in \bZ$.
\end{proof}

\subsection{Log \texorpdfstring{$\TC$}{TC} of even-periodic sphere spectra} The next lemma is the first step towards the identification of $\TC(\bS[x])_p$ and $\TC(\bS[x],\<x\>)_p$. 
\begin{lemma}\label{lem:varphi_pi_p-equivalence}
For $i > 0$, the cyclotomic structure map $\varphi_{p;i} \: \THH(\bS[x];i)_p \to \THH(\bS[x];pi)^{tC_p}_p$ is an equivalence. 
\end{lemma}
\begin{proof}
Restricting diagram~\eqref{eq:Deltap-varphip} with $A = \Th_\bS(\xi)$ and $q = i-1$ along the map $\bS^{2d} \to \Th_\bS(\xi)$, we obtain the $C_i$-equivariant commutative outer rectangle in the diagram 
\[
\xymatrix@-1pc{
(\bS^{2d})^{\otimes i} \ar[rr] \ar[d]_-{\Delta_p}
	&& \bT_+ \wedge_{C_i} (\bS^{2d})^{\otimes i} \ar[rr]^{\simeq}
		\ar[d]_-{\bar\Delta_p}
	&& \THH(\Th_\bS(\xi); i) \ar[d]^-{\varphi_{p;i}} \\
((\bS^{2d})^{\otimes pi})^{tC_p} \ar[rr]
	&& (\bT_+ \wedge_{C_{pi}} (\bS^{2d})^{\otimes pi})^{tC_p} \ar[rr]^{\simeq}
	&& \THH(\Th_\bS(\xi); pi)^{tC_p}\,.
}
\]
The two left-hand horizontal arrows are the evident inclusions, and the two right-hand horizontal equivalences arise like the equivalence~\eqref{eq:THHbSxwi-identification} in the proof of Proposition~\ref{prop:THHSP}. The two maps $\Delta_p $ and $\bar\Delta_p $ are $p$-completions by the affirmed Segal
conjecture for $C_p$~\cite{Lin80},~\cite{Gun81}. It follows that $\varphi_{p;i}$ is a $p$-equivalence. 
\end{proof}

The pair $(\bS,\xi_*)$ is cyclotomically good by Proposition~\ref{prop:sufficient-criteria-cycl-good}, and the collapse $\modmod{\bS[x]}{\<x\>}$ is equivalent to $\bS$ (see Example~\ref{ex:collapse-canonical}). Hence the $\TC$ repletion--residue sequence for $(\bS[x],\xi,\id)$ from Corollary~\ref{cor:rep-res-TC} takes the form
\[ \TC(\bS)_p \xrightarrow{\trf} \TC(\bS[x])_p
	\xrightarrow{\rho} \TC(\bS[x],\<x\>)_p
	\xrightarrow{\res} \Sigma \TC(\bS)_p. 
      \]      
\begin{theorem}\label{thm:TCSP}
Let $|x|=2d$ with $d > 0$. There are $p$-equivalences
\begin{align*}
\TC(\bS[x])_p &\simeq_p \TC(\bS) \vee
	\bigvee_{i>0} \Sigma ((\bS^{2d})^{\otimes i})_{hC_i}  \qquad\text{and}\\ 
\TC(\bS[x], \<x\>)_p &\simeq_p \TC(\bS)[S^1] \vee
	\bigvee_{i>0} \Sigma ((\bS^{2d})^{\otimes i})_{hC_i} \,,
\end{align*}
and the map $\trf$ is null-homotopic as a map of $\TC(\bS)_p$-modules. 
\end{theorem}
\begin{proof}
Let $X_* = \THH(\bS[x])$ or $X_* = \THH(\bS[x],\<x\>)$, viewed as $\bZ_{\ge0}$-graded $L_p$-twisted cyclotomic spectra.  As spectra with $\bT$-action, $X_i \simeq \bT_+ \wedge_{C_i} (\bS^{2d})^{\otimes i}$ for each $i>0$ by Proposition~\ref{prop:THHSP}, while $X_0 \simeq \bS$ or $X_0 \simeq \bS[S^1(0)]$ according to the case.  Hence $X_i$ is $2di$-connective for each $i\ge0$, so that $X_*$ is properly connective in the sense of Definition~\ref{def:propconn}, and there is a natural equivalence
\begin{equation} \label{eq:TCX-as-product}
\textstyle \TC(X)_p \simeq \TC(X_0)_p \times \prod_{p \nmid i > 0} \TC(X_{[i]_p})_p
\end{equation}
by Proposition~\ref{prop:TC-prodIp}.

First consider the two cases with $i=0$.  When $X_0 = \THH(\bS[x]); 0) = \THH(\bS)$ we immediately have $\TC(X_0)_p = \TC(\bS)_p$.  When $X_0 = \THH(\bS[x],\<x\>; 0)$, the weight $0$ part of Proposition~\ref{prop:rep-res-can} leads to a repletion--residue sequence
\begin{equation} \label{eq:TCS-rep-res}
\TC(\bS)_p \xrightarrow{\trf_0} \TC(\bS)_p
        \xrightarrow{\rho_0} \TC(\bS[x], \<x\>; 0)_p
        \xrightarrow{\res_0} \Sigma \TC(\bS)_p
\end{equation}
of $\TC(\bS)_p$-modules. In the commutative square 
\[
\xymatrix@-1pc{
\TC(\bS)_p \ar[rr]^-{\trf_0} \ar[d]_-{\beta}
	&& \TC(\bS)_p \ar[d]^-{\beta} \\
\THH(\bS)_p \ar[rr]^-{\trf_0}
	&& \THH(\bS)_p
}
\]
with $\beta$ as in~\eqref{eq:beta-TC-THH}, both maps $\pi_0(\beta) \: \bZ_p \to \bZ_p$ are surjective because $\bS \to \TC(\bS)_p \to \THH(\bS)_p$ is a $p$-equivalence, and injective because $\pi_0$ of the homotopy fiber vanishes,  cf.~\cite{BHM93}*{Thm.~5.17}, \cite{Rog02}*{Cor.~1.21}. Lemma~\ref{lem:repletion-fiber-as-THHR-module} shows that $\pi_0(\trf_0)$ is the zero homomorphism for $\THH$, hence also for $\TC$. This implies that~$\trf_0$ in~\eqref{eq:TCS-rep-res} is null-homotopic as
a map of $\TC(\bS)_p$-modules, so that
\[
\TC(\bS[x],\<x\>; 0)_p \simeq \TC(\bS)_p \wedge S^1_+
	= \TC(\bS)_p [S^1] \,.
\]

Next, fix an $i>0$ with $p \nmid i$, so that $[i]_p = \{i, pi, p^2 i,
\dots\} = \{p^e i \mid e\ge0\}$.  Then $\TC(X_{[i]_p})_p$ is the equalizer
of the two maps in~\eqref{eq:TCXipp}, which equals the limit of the
$p$-complete diagram in Figure~\ref{fig:lim-system}.
\begin{figure}\[\xymatrix{
& ((\bT_+ \wedge_{C_i} (\bS^{2d})^{\otimes i})^{tC_p})^{h\bT}_p = 0 \\
(\bT_+ \wedge_{C_i} (\bS^{2d})^{\otimes i})^{h\bT}_p
	\ar[r]^-{\can_i}
	\ar[dr]_-{\varphi_{p;i}^{h\bT}}^-{\simeq}
& (\bT_+ \wedge_{C_i} (\bS^{2d})^{\otimes i})^{t\bT}_p
	\ar[u]^-{G_i}_{\simeq} \\
& ((\bT_+ \wedge_{C_{pi}} (\bS^{2d})^{\otimes pi})^{tC_p})^{h\bT}_p \\
(\bT_+ \wedge_{C_{pi}} (\bS^{2d})^{\otimes pi})^{h\bT}_p
	\ar[r]^-{\can_{pi}}
	\ar[dr]_-{\varphi_{p;pi}^{h\bT}}^-{\simeq}
& (\bT_+ \wedge_{C_{pi}} (\bS^{2d})^{\otimes pi})^{t\bT}_p
	\ar[u]^{G_{pi}}_-{\simeq} \\
& ((\bT_+ \wedge_{C_{p^2 i}} (\bS^{2d})^{\otimes p^2 i})^{tC_p})^{h\bT}_p  \\
{\dots} \ar[r]^-{\can_{p^2 i}}
& (\bT_+ \wedge_{C_{p^2 i}} (\bS^{2d})^{\otimes p^2 i})^{t\bT}_p
	\ar[u]^-{G_{p^2 i}}_-{\simeq}
}
\]
\caption{The limit system for $\TC(X_{[i]_p})_p$}\label{fig:lim-system}
\end{figure}
For brevity, let $j = p^e i$.  The maps $G_j$ in the diagram are equivalences by \cite{BBLNR14}*{Prop.~3.8} or \cite{NS18}*{Lem.~II.4.1}.  The maps $\varphi_{p; j}$ in the diagram are equivalences by Lemma~\ref{lem:varphi_pi_p-equivalence}, hence so are the~$\varphi_{p; j}^{h\bT}$.  For $p \nmid i$ the $C_p$-action on $\bT_+ \wedge_{C_i} (\bS^{2d})^{\otimes i}$ is free, and its $C_p$-Tate fixed point spectrum is zero.  The cofiber sequence~\eqref{eq:T-NRseq} for $X_j$ takes the form
\begin{multline*}
\Sigma ((\bS^{2d})^{\otimes j})_{hC_j}
  \simeq \Sigma (\bT_+ \wedge_{C_j} (\bS^{2d})^{\otimes j})_{h\bT} \\
  \xrightarrow{N^h} (\bT_+ \wedge_{C_j} (\bS^{2d})^{\otimes j})^{h\bT}
  \xrightarrow{R^h} (\bT_+ \wedge_{C_j} (\bS^{2d})^{\otimes j})^{t\bT}
\,.
\end{multline*}
By Lemma~\ref{lem:section} below, with $Y = (S^{2d})^{\otimes j}$, the map
$R^h = \can_j$ admits a section.  Hence the limit in question agrees
with the limit of a sequence
\[
\dots \longto (\bT_+ \wedge_{C_{p^e i}} (\bS^{2d})^{\otimes p^e i})^{h\bT}_p
  \longto (\bT_+ \wedge_{C_{p^{e-1} i}} (\bS^{2d})^{\otimes p^{e-1} i})^{h\bT}_p
  \longto \dots \longto 0
\]
where each map admits a section. The limit is therefore the product
\[
\TC(X_{[i]_p})_p \simeq \prod_{e\ge0}
	(\Sigma ((\bS^{2d})^{\otimes p^e i})_{hC_{p^e i}})_p
\]
of the fibers of these maps. 

We have now identified the factors in~\eqref{eq:TCX-as-product}.
Each $j>0$ can be written as $p^e i$ for unique $e\ge0$ and
$p \nmid i>0$, so we can rewrite this equivalence as
\[
\TC(X) \simeq_p \TC(X_0) \times
	\prod_{j>0} \Sigma ((\bS^{2d})^{\otimes j})_{hC_j} \,.
\]
Since $\Sigma ((\bS^{2d})^{\otimes j})_{hC_j}$ is $(2dj+1)$-connective
for each $j>0$, this product is equivalent to the asserted sum.
\end{proof}

The following lemma extracts a consequence of the Segal--tom\,Dieck splitting and the proven Segal conjecture, and was used in the proof above.

\begin{lemma} \label{lem:section}
Let $j = p^e i$ with $e\ge0$ and $p \nmid i>0$, and let $Y$ be a
finite based $C_j$-CW-space.  Suppose that $C_i$ acts trivially on
$\pi_*(\Sigma^\infty Y \otimes \bS/p)$.  The $p$-completed canonical map
\[
R^h = \can_j \: (\bT_+ \wedge_{C_j} \Sigma^\infty Y)^{h\bT}_p
	\longto (\bT_+ \wedge_{C_j} \Sigma^\infty Y)^{t\bT}_p
\]
admits a section.
\end{lemma}

\begin{proof}
Using the Wirthm{\"u}ller equivalence~\cite{Wir74}*{Thm.~2.1}, \cite{LMS86}*{Thm.~II.6.2}, the cofiber sequence~\eqref{eq:T-NRseq} for $\bT_+ \wedge_{C_j} \Sigma^\infty Y$ can be rewritten in the form
\[
\Sigma (\Sigma^\infty Y)_{hC_j}
	\xrightarrow{N^h} \Sigma (\Sigma^\infty Y)^{hC_j}
	\xrightarrow{R^h} \Sigma (\Sigma^\infty Y)^{tC_j}
\,.
\]
Hence it suffices to prove that the homotopy norm map
\[
N^h \: (\Sigma^\infty Y)_{hC_j} \longto (\Sigma^\infty Y)^{hC_j}
\]
admits a retraction after $p$-completion.  For this, we momentarily
work with genuinely equivariant spectra.  The Segal--tom\,Dieck splitting
\cite{Seg71}, \cite{tD75}*{Satz~2}, \cite{LMS86}*{Ch.~V} for the based
$C_{p^e}$-CW-space~$Y$ is an equivalence
\[
\textstyle \bigvee_{0 \le d \le e} (\Sigma^\infty Y^{C_{p^d}})_{hC_{p^{e-d}}}
	\xrightarrow{\simeq} (\Sigma^\infty Y)^{C_{p^e}} \,.
\]
Its $d=0$ component is given by the norm map
\[
N \: (\Sigma^\infty Y)_{hC_{p^e}} \longto (\Sigma^\infty Y)^{C_{p^e}} \,,
\]
whose composite with the comparison map $\Gamma_e \: (\Sigma^\infty
Y)^{C_{p^e}} \to (\Sigma^\infty Y)^{hC_{p^e}}$ is~$N^h$ for~$C_{p^e}$.
A choice of retraction~$r_0$ to the $d=0$ component of the Segal--tom
Dieck splitting then produces the central part of the following
homotopy commutative diagram.
\[
\xymatrix{
(\Sigma^\infty Y)_{hC_j} \ar[d]_-{N^h} \ar[r]^-{\trf}_-{\simeq_p}
	\ar@(u,u)[rrr]^-{i \cdot \id}_-{\simeq_p}
	& (\Sigma^\infty Y)_{hC_{p^e}} \ar[d]_-{N^h} \ar[dr]^-{N}
	\ar[r]^-{\id}
	& (\Sigma^\infty Y)_{hC_{p^e}} \ar[r]^-{\pr}_-{\simeq_p}
	& (\Sigma^\infty Y)_{hC_j} \\
(\Sigma^\infty Y)^{hC_j} \ar[r]^-{\res}_-{\simeq_p}
	& (\Sigma^\infty Y)^{hC_{p^e}}
	& (\Sigma^\infty Y)^{C_{p^e}} \ar[l]_-{\Gamma_e}^-{\simeq_p}
	\ar[u]_-{r_0}
}
\]
Here $\pr$, $\res$ and $\trf$ are the projection, restriction and transfer
maps associated to the inclusion $C_{p^e} \subset C_j$.  The first two
are $p$-equivalences, since homotopy fixed point and homotopy orbit
spectral sequence arguments using the hypothesis that $C_i$ acts
trivially on $\pi_*(\Sigma^\infty Y \otimes \bS/p)$ show that the $C_{p^e}$-maps
$(\Sigma^\infty Y)^{hC_i} \to \Sigma^\infty Y \to (\Sigma^\infty
Y)_{hC_i}$ are $p$-equivalences.  The composite map $\pr \circ \trf$
is multiplication by the index~$i$ of $C_{p^e}$ in $C_j$, hence it (and
$\trf$) is also a $p$-equivalence.  Finally, the comparison map $\Gamma_e$
is a $p$-equivalence by the proven Segal conjecture~\cite{Car84},
\cite{Rav84} for~$C_{p^e}$, using the finiteness of~$Y$.  Taken together,
this means that the composite
\[
(\pr \circ \trf)^{-1} \circ \pr \circ r_0 \circ (\Gamma_e)^{-1} \circ \res
\]
(in the homotopy category, with all maps implicitly $p$-completed)
provides a left inverse to the left-hand homotopy norm map.
\end{proof}

\appendix

\section{Comparison with point-set level Thom spectra}\label{app:Thom-comparison}
In~\cite{SS19}*{\S2}, the second and third authors consider the topologically enriched symmetric monoidal category $\cW$ given by Quillen's localization construction applied to the category $\cV^{\mathrm{iso}}$ of Euclidean vector spaces $\mathbb R^n, n\geq 0,$ and linear isometric isomorphisms. The classifying space of $\cW$ is $B\cW \simeq \mathbb Z \times BO$ and there is a strong symmetric monoidal functor $\cV^{\mathrm{iso}} \to \cW$. Its induced map of classifying spaces $B\cV^{\mathrm{iso}} \to B\cW$ is an $\bE_{\infty}$-map exhibiting $B\cW$ as the group completion of $B\cV^{\mathrm{iso}} \simeq \coprod_{n\geq 0}BO(n)$. The continuous functor \begin{equation}\label{eq:fun-F_VS^V} F_{-}(S^{-})\: \cW^{\op} \to \Sp^{O}\end{equation} to orthogonal spectra from~\cite{SS19}*{Lem.~4.3} extends the functor $(\cV^{\mathrm{iso}})^\op \to \Sp^{O}$ sending $V$ to $\Sigma^{\infty}S^V$. Passing to underlying $\infty$-categories, the latter functor co-restricts to an $\bE_{\infty}$-map $\coprod_{n\geq 0}BO(n) \to \Pic_{\bS}$, and $F_{-}(S^{-})\: \cW^{\op} \to \Sp^{O}$ induces an $\bE_{\infty}$ $J$-map $J \: B\cW^{\op} \to  \Pic_{\bS}$, the extension of the former $\bE_{\infty}$-map over the group completion $\coprod_{n\geq 0}BO(n) \to \mathbb Z \times BO$. 

\begin{definition}\phantomsection\label{def:ThW}\begin{enumerate}[(i)]
  \item We let $J_! \: \cS_{/ B\cW^{\op}} \to  \cS_{/ \Pic_{\bS}}$ be the symmetric monoidal functor obtained from applying Corollary~\ref{cor:construction-Fiota} to  $y_{\Pic_{\bS}}\circ J \:  B\cW^{\op} \to \cS_{/ \Pic_{\bS}}$.
  \item We let $\Th_{\cW} \: \cS_{/ B\cW^{\op}} \to \Sp$ be the symmetric monoidal functor obtained from applying Corollary~\ref{cor:construction-Fiota} to $\iota \circ J \: B\cW^{\op}\to \Sp$ (with $\iota\:\Pic_\bS\to\Sp$ the inclusion). 
\end{enumerate}
\end{definition}
Under the symmetric monoidal equivalences $\Fun(P^{\op},\cS) \simeq \cS_{/ P}$ for $P=B\cW^{\op}$ and $P= \Pic_{\bS}$, the functor $J_!$ can be identified with the left Kan extension along $J^\op$.
\begin{corollary}\label{cor:ThW-as-composite}
  The functors $\Th_{\bS}\circ J_!$ and $\Th_{\cW}$ are equivalent as colimit-preserving symmetric monoidal functors $\cS_{/ B\cW^{\op}} \to \Sp$.
\end{corollary}
\begin{proof}
  By construction, we have $\Th_{\bS}\circ J_! \circ y_{B\cW^{\op}} \simeq \Th_{\bS} \circ y_{\Pic_{\bS}} \circ J \simeq \iota \circ J$ and $\Th_{\cW} \circ y_{B\cW^{\op}} \simeq \iota \circ J$. The claim follows from Proposition~\ref{prop:extenstion-to-cocompletion}.  
\end{proof}

The $1$-categorical functor category $\cS^{\cW}= \Fun(\cW,\cS)$ admits a \emph{$\cW$-model structure} making it Quillen equivalent to the model category $\cS/(\mathbb Z \times BO)$, that is, the $1$-categorical comma category with the overcategory model structure (see~\cite{SS19}*{\S3 and \S6}). Via the compatibility of the model categorical overcategory with the slice category construction~\cite{L:HTT}*{Lem.~6.1.3.13}, $\cS^{\cW}$ is therefore a symmetric monoidal model category representing  $\Fun(B\cW,\cS) \simeq \cS_{/ B\cW^{\op}}$. By construction, the strong symmetric monoidal left Quillen functor $\bS^{\cW}\: \cS^{\cW} \to \Sp^O$ of~\cite{SS19}*{\S4.2} extends the strong symmetric monoidal functor $F_{-}(S^{-})\: \cW^{\op} \to \Sp^{O}$ of~\eqref{eq:fun-F_VS^V}. With this observation, Definition~\ref{def:ThW} implies:
\begin{corollary}\label{cor:SW-induces-ThW}
The symmetric monoidal functor of $\infty$-categories induced by $\bS^{\cW}$ is equivalent to $\Th_{\cW}$.\qed
\end{corollary}
In~\cite{SS19}, $\bS^{\cW}$ was constructed as one model for a graded Thom spectrum functor for virtual real vector bundles. The combination of Corollaries~\ref{cor:ThW-as-composite} and~\ref{cor:SW-induces-ThW} implies that this point-set level construction is symmetric monoidally equivalent to the $\infty$-categorical one. 

An analogous discussion applies to the $\cJ$-spaces considered in~\cite{SS12}: Passing to classifying spaces, the category $\cJ$ of~\cite{SS12}*{\S4} models the group completion $QS^0$ of $\coprod_{n\geq 0}B\Sigma_n$. There is a strong symmetric monoidal functor $F_{-}(S^{-})\: \cJ^{\op} \to \Sp^{\Sigma}$ to symmetric spectra that extends to a strong symmetric monoidal left Quillen functor $\bS^{\cJ} \: \cS^{\cJ} = \Fun(\cJ, \cS) \to \Sp^{\Sigma}$, see~\cite{SS12}*{Lem.~4.22 and~(4.5)}. The map of classifying spaces
$B\cJ^{\op} \to \Pic_{\bS}$
induced by $F_{-}(S^{-})\: \cJ^{\op} \to \Sp^{\Sigma}$ factors as the composite of the map $\rho\: B\cJ^{\op} \to B\cW^{\op}$ induced by the inclusions $\Sigma_n \to O(n)$ (or the unit $\bS \to \mathrm{ko}$) and the above $J \: B\cW^{\op} \to \Pic_{\bS}$. 

Arguing as in the case of Corollaries~\ref{cor:ThW-as-composite} and~\ref{cor:SW-induces-ThW}, we obtain:
\begin{corollary}\label{cor:J-space-comparison}
The symmetric monoidal functor $\Th_{\cJ}\: \cS_{/ B\cJ^{\op}} \to \Sp$ of $\infty$-ca\-te\-gories induced by $\bS^{\cJ}$ is equivalent to the composite \[ \cS_{/ B\cJ^{\op}} \xrightarrow{(J\circ\rho)_!}  \cS_{/ \Pic_{\bS}} \xrightarrow{\Th_{\bS}} \Sp.\qedhere\]
\end{corollary}

Applying Proposition~\ref{prop:additive-adjunction} and Corollary~\ref{cor:Giota-lax-sym-mon}, it follows like in Definition~\ref{def:ThR} that the functors considered above participate in adjunctions
\begin{equation} (\Th_{\cW},\Omega^{\cW}), \quad (\Th_{\cJ},\Omega^{\cJ}), \quad (J_!,J^*) \text{ and}\qquad ((J\circ\rho)_!, (J\circ\rho)^*)
\end{equation}
in each of which the right adjoint is a lax symmetric monoidal functor. Moreover, using~\cite{L:HA}*{Rem.~7.3.2.13} as in the case of~\eqref{eq:R-mod-E-k-adjunction}, each of these adjunctions induces an adjunction on the associated categories of $\bE_k$-algebras for $1 \leq k \leq \infty$. 

\begin{proof}[Proof of Proposition~\ref{prop:Omega-S-multiplicative-monoid}]
For the model categorical counterpart of $\Omega^\cJ$, the statement of Proposition~\ref{prop:Omega-S-multiplicative-monoid} is verified in~\cite{SS12}*{Prop.~4.24}. Since $J\circ \rho\: QS^0 \to \Pic_\bS$ is an isomorphism on $\pi_*$ for $*=0,1$, the equivalence $(J\circ\rho)^* \circ \Omega^\bS \simeq \Omega^{\cJ}$ resulting from the above discussion induces a bijection between the path components of the relevant homotopy fibers, which respects the multiplication and the $\{\pm 1\}$-action. 
\end{proof}

\section{Comparison with point-set level log \texorpdfstring{$\THH$}{THH}}\label{app:pointsetlevel-logTHH-comparision}
The purpose of this appendix is to compare the notion of log $\THH$ introduced in Section~\ref{sec:logTHH} with that of \cite{RSS15}*{Def.~4.6} and \cite{RSS18}*{Def.~4.1} when $M$ is repetitive in the sense of \cite{RSS15}*{Def.~6.4}.

One building block for the prelog ring spectra of the cited paper are commutative $\cJ$-space monoids~\cite{SS12}*{\S{}4}. When equipped with one of the positive $\cJ$-model structures, the category of commutative $\cJ$-space monoids $\cC\cS^{\cJ}$ models the $\infty$-category $\CAlg(\cS_{ /B\cJ^{\op}})$. It also admits a group completion model structure where the fibrant replacement is a model for the group completion~\cite{S16}*{Thm.~1.6}. This is compatible with the group completions discussed in Subsection~\ref{subsec:weight-graded-Thom}:
\begin{lemma}\label{lem:compatible-group-completions}
If $\gamma \: M \to M^\gp$ is a group completion of commutative $\cJ$-space monoids, then $(J\circ \rho)_!$ sends the map in $\CAlg(\cS_{/ B\cJ^{\op}})$ represented by $\gamma$ to a group completion in $\CAlg(\cS_{/ \Pic_{\bS}})$. 
\end{lemma}
\begin{proof}
Since the $\bE_\infty$-spaces $B\cJ^{\op}$ and $\Pic_{\bS}$ are grouplike, group completions of $\bE_\infty$-algebras in the slice categories $\cS_{/ B\cJ^{\op}}$ and $\cS_{/ \Pic_{\bS}}$ are detected on their underlying $\bE_{\infty}$-spaces (compare Lemma~\ref{lem:grouplike-in-slice}). 
\end{proof}

A $\cJ$-space prelog ring spectrum (as considered in~\cites{SS12,S14,RSS15,RSS18}) consists of an $M \in\cC\cS^{\cJ}$, a commutative symmetric ring spectrum $A \in \cC\Sp^{\Sigma}$, and a structure map $\alpha \: M \to \Omega^{\cJ}(A)$ in $ \cC\cS^{\cJ}$ or, equivalently, a structure map $\bar{\alpha}\: \bS^{\cJ}[M] \to A$ in $\cC\Sp^{\Sigma}$. The commutative $\cJ$-space monoid $M$ represents an object in the $\infty$-category $\CAlg(\Fun(B\cJ,\cS)) \simeq \CAlg(\cS_{/ B\cJ^{\op}})$. We write $\xi(M)$ for its image under $(J\circ \rho)_!\: \CAlg(\cS_{/ B\cJ^{\op}}) \to  \CAlg(\cS_{/ \Pic_{\bS}}) $ and view $\xi(M)$ as a graded object with the canonical $\pi_0$-grading of Definition~\ref{def:canonical-pi0grading}. As in \emph{op.\,cit.}, we write $M_{h\cJ} = \hocolim_{\cJ}(M)$ for the homotopy colimit and note that this $\bE_\infty$-space represents the total space of $\xi(M)$. So in this case, the canonical $\pi_0$-grading is a $\pi_0(M_{h\cJ})$-grading. By Corollary~\ref{cor:J-space-comparison}, $\bar{\alpha}$ induces a structure map $\bar{\alpha} \: \Th_{\bS}(\xi(M)) \to A$ that allows us to view $(A,\xi(M),\bar{\alpha})$ as a prelog $\bE_{\infty}$-ring.  An analogous discussion applies to the $\cW$-space prelog ring spectra considered in~\cite{SS19}. 

A commutative $\cJ$-space monoid $M$ is \emph{repetitive} in the sense of~\cite{RSS15}*{Def.~6.4} if $M$ is not concentrated in $\cJ$-space degree $0$ and if the group completion map induces an equivalence $M\to (M^\gp)_{\geq 0}$. If $M$ is repetitive, $M_{hJ} \to \pi_0(M_{h\cJ})\times_{\pi_0(M^\gp_{h\cJ})} M^\gp_{h\cJ}$ is an equivalence, meaning that $M$ is $\pi_0$-replete in the sense of Definition~\ref{def:zero-replete}. 

\begin{proposition}\label{prop:graded-THH-of-Thomcompatible}
Let $M$ be a repetitive commutative $\cJ$-space monoid that is flat as a $\cJ$-space. Then $\bS^{\cJ}[B^\rep(M)]$ represents $\THH(\Th_{\bS}(\xi(M)),\xi(M))$ in the $\infty$-category $\CAlg(\Sp)$.
\end{proposition}
\begin{proof}
Let $ B^\rep(M)$ be  the replete bar construction of $M$ in commutative $\cJ$-space monoids. The observation before the proposition and Lemma~\ref{lem:compatible-group-completions} imply that $(J\circ\rho)_!\: \CAlg(\cS_{/ B\cJ^{\op}}) \to  \CAlg(\cS_{/ \Pic_{\bS}})$ sends the object represented by $B^\rep(M)$ to $B^\zerorep((J\circ\rho)_!(M))$. With this, the claim follows from Corollary~\ref{cor:J-space-comparison} and Proposition~\ref{prop:Bcy-THH-comparison-over-PicS}. 
  \end{proof}
\begin{corollary}\label{cor:abstract-comparison-with-old-prelog}
  Let $(A,M,\alpha)$ be $\cJ$-space prelog ring spectrum with $M$ repetitive. Then the log $\THH$ spectrum $\THH(A,M,\alpha)$ considered in~\cite{RSS15}*{Def.~4.6}, \cite{RSS18}*{Def.~4.1} represents the log $\THH$ of the associated  prelog $\bE_{\infty}$-ring as an $\bE_{\infty}$-ring.\qed 
\end{corollary}

\subsection{Comparison with direct image log structures}
Let $A$ be an $\bE_{\infty}$-ring, let $d > 0$, and let $x \in \pi_{2d}(A)$ be a homotopy class such that the localization map $j\: A \to A[1/x]$ exhibits $A$ as the connective cover of $A[1/x]$. In this situation,~\cite{S14}*{Con.~4.2} provides a $\cJ$-space prelog ring spectrum $(A,D(x),\alpha^x)$ that we can view as a prelog $\bE_\infty$-ring $(A,\xi(D(x)),\bar\alpha^x)$ by the above construction. As in Construction~\ref{constr:map-to-direct-image}, the adjoint of the composite $\bS^\cJ[D(x)] \simeq \Th_\bS(\xi(D(x))) \to A \to A[1/x]$  factors through $\iota\: \GL_1^\bS(A[1/x]) \to \Omega^\bS(A[1/x])$ and therefore induces a canonical map $(A,\xi(D(x)),\bar\alpha^x) \to (A,j_* \xi^{\GL_1}, j_*\bar\iota)$ to the direct image log $\bE_\infty$-ring determined by~$j$ (compare Example~\ref{ex:directimagelog}).

\begin{lemma}
The canonical map $(A,\xi(D(x)),\bar\alpha^x) \to (A,j_*\GL_1^{\bS}(A[1/x]), \alpha)$ is an $\bE_\infty$-logification in the sense of Definition~\ref{def:E_k-logification}.   
\end{lemma}  
\begin{proof}
This is analogous to~\cite{S14}*{Lem.~4.7}. 
\end{proof}
Suppose in addition that  $\pi_{i(2d+2)+2d-1}(A) \cong 0$ for all $i \geq 1$. 
Then by Remark~\ref{rem:vanishing-for-one-gen-case}, we get a prelog $\bE_2$-ring $(A,\xi_{2d},\bar\alpha)$ determined by $x$. Construction~\ref{constr:map-to-direct-image} and Lemma~\ref{lem:map-to-direct-image} show that in this situation the canonical map $(A,\xi_{2d},\bar\alpha) \to (A,j_* \xi^{\GL_1}, j_*\bar\iota)$ is an $\bE_2$-logification. Combining this with the last lemma and applying Theorem~\ref{thm:logification-invariance}, we obtain: 
\begin{corollary}\label{cor:Dx-comparison}
  The maps $ (A,\xi_{2d},\bar\alpha) \to (A,f_* \xi^{\GL_1}, j_*\bar\iota) \leftarrow (A,\xi(D(x)),\bar\alpha^x)$ induce equivalences when applying log $\THH$. \qed
\end{corollary}
In the corollary, the log $\THH$ of each outer term is formed with respect to the canonical $\pi_0$-grading, which is a grading by the (additive) monoid $\bZ_{\geq 0}$ in these two cases. The middle term is graded by $\pi_0(f_* \GL_1^\bS(A[1/x])) \cong (\pi_0(A)^\times/\{\pm 1\})[x]$. Our connective cover assumptions imply that the group completion of this monoid is $(\pi_0(A)^\times/\{\pm 1\})[x^{\pm 1}]$. Hence the map from $\bZ_{\geq 0}$ to $\pi_0(f_* \GL_1^\bS(A[1/x]))$ is exact. Lemma~\ref{lem:change-grading-of-canonical} implies that we can change the implicit grading of the outer terms to a $(\pi_0(A)^\times/\{\pm 1\})[x]$-grading without altering the log $\THH$. 

Using Corollary~\ref{cor:abstract-comparison-with-old-prelog}, we deduce that the log $\THH$ of the $\bE_2$-prelog rings $(ku,\< u \>)$, $(\kup,\< u \>)$, and $(\ell,\< v_1\>)$ from Definition~\ref{def:ku-ell-BPn} is equivalent to that defined and studied in~\cites{RSS15,RSS18}.

\end{document}